\documentclass[10pt]{amsart}
\usepackage[cp1251]{inputenc}
\usepackage[german,english]{babel}
\usepackage{amsmath}
\usepackage{amssymb}
\usepackage{amsfonts}

\begin{document}
\renewcommand{\refname}{References}

\thispagestyle{empty}

\title[Application of the Method of Approximation
of Iterated 
It\^{o} Stochastic Integrals]
{Application of the Method of Approximation of Iterated 
It\^{o} Stochastic Integrals Based on Generalized 
Multiple Fourier Series to the High-Order Strong Numerical Methods for 
Non-Commutative Semilinear Stochastic Partial Differential Equations}
\author[D.F. Kuznetsov]{Dmitriy F. Kuznetsov}
\address{Dmitriy Feliksovich Kuznetsov
\newline\hphantom{iii} Peter the Great Saint-Petersburg Polytechnic University,
\newline\hphantom{iii} Institute of Applied Mathematics and Mechanics,
\newline\hphantom{iii} Department of Mathematics,
\newline\hphantom{iii} Polytechnicheskaya ul., 29,
\newline\hphantom{iii} 195251, Saint-Petersburg, Russia}
\email{sde\_kuznetsov@inbox.ru}
\thanks{\sc Mathematics Subject Classification: 60H05, 60H10, 42B05, 42C10}
\thanks{\sc Keywords: Non-commutative semilinear
stochastic partial differential equation,
Infinite-dimensional $Q$-Wiener process, 
Iterated It\^{o} stochastic integral,
Generalized multiple Fourier series,
Multiple Fourier--Legendre series,
Legendre polynomial, Mean-square approximation, Expansion, Trace class noise.}

\maketitle {\small
\begin{quote}
\noindent{\sc Abstract.} 
We consider a method for the approximation of iterated stochastic
integrals of arbitrary multiplicity $k$ $(k\in \mathbb{N})$
with respect to the 
infinite-dimensional $Q$-Wiener process using the mean-square
approximation method of iterated It\^{o} stochastic integrals with 
respect to the scalar standard Wiener processes based on 
generalized multiple Fourier series. The case of multiple Fourier--Legendre 
series is considered in details.
The results of the article can be applied to construction of 
high-order strong numerical methods (with respect to the temporal 
discretization) for the approximation of mild solution for non-commutative semilinear 
stochastic partial differential equations with multiplicative trace class noise.
\medskip
\end{quote}
}

\vspace{3mm}

\section{Introduction}

There exists a lot of publications on the subject of 
numerical integration of stochastic partial differential
equations (SPDEs) (see, for example, \cite{1}-\cite{26}).
One of the perspective approaches to the construction of high-order
strong numerical methods (with respect to the temporal 
discretization) for SPDEs is based
on the Taylor formula in Banach spaces and exponential formula for 
the mild solution of 
SPDEs \cite{12} (2015), 
\cite{13} (2016).
As shown in \cite{12} (2015) and \cite{18} (2007) the exponential Milstein type 
approximation method has the strong order of convergence 
$1.0-\varepsilon$ (where $\varepsilon$
is an arbitrary small posilive real number) \cite{12} or $1.0$ \cite{18}. 
In \cite{13} the exponential Wagner--Platen type numerical 
method 
for SPDEs
with strong order $1.5-\varepsilon$ (where $\varepsilon$
is an arbitrary small posilive real number) has been considered.
An important feature of these numerical methods is a presence 
in them of the so-called iterated stochastic integrals with respect 
to the infinite-dimensional $Q$-Wiener process \cite{20}.
Approximation of these stochastic integrals is a complex problem.
This problem can be significantly simplified 
if special commutativity conditions be fulfilled
\cite{12}, \cite{13}.
In \cite{26} (2019) two methods of the mean-square approximation of simplest
iterated (double) 
stochastic integrals with respect to the infinite-dimensional 
$Q$-Wiener process are considered and theorems on the convergence of 
these methods are given (the basic idea about Karhunen--Loeve
expansion of the Brownian bridge process was taken from monograph
\cite{27} (1988, In Russian)).
It is important to note that the approximation of iterated 
stochastic integrals with respect to the infinite-dimensional 
$Q$-Wiener process can be reduced to the approximation of iterated 
It\^{o} stochastic integrals with respect to the scalar standard 
Wiener processes. In a lot of author's publications  
\cite{eee}-\cite{1003} the effective methods for the mean-square
approximation of 
iterated 
It\^{o} and Stratonovich stochastic integrals with respect to the scalar standard Wiener 
processes were proposed and developed. 
One of these methods \cite{30} (also see \cite{30a}-\cite{1003})
is based on generalized multiple Fourier series,
in particular, on multiple Fourier--Legendre series. 
The purpose of this article is an adaptation of the method 
\cite{30}-\cite{1003}
for the mean-square approximation of iterated stochastic integrals 
of multiplicity $k$ $(k\in \mathbb{N})$
with respect to the finite-dimensional approximation 
of the infinite-dimensional $Q$-Wiener process.

Let $U, H$ be separable $\mathbb{R}$-Hilbert spaces and
$L_{HS}(U,H)$ be a space of Hilbert--Schmidt operators mapping from $U$ to $H.$
Let $(\Omega, {\bf F}, \sf{P})$ be a probability space 
with a normal filtration $\{{\bf F}_t, t\in [0, \bar{T}]\}$
\cite{20}, let ${\bf W}_t$ be an $U$-valued $Q$-Wiener process 
with respect to $\{{\bf F}_t, t\in [0, \bar{T}]\},$  
which has a covariance trace class operator $Q\in L(U)$. 
Here $L(U)$ denotes all bounded linear operators
mapping from $U$ to $U$. 
Consider the semilinear parabolic 
SPDE

\vspace{-2mm}
\begin{equation}
\label{xx1}
dX_t = \left(A X_t +F(X_t)\right)dt + B(X_t)d{\bf W}_t,\ \ \
X_0=\xi,\ \ \ t\in [0, \bar{T}],
\end{equation}

\vspace{3mm}
\noindent
where nonlinear operators $F,$ $B$ ($F:$ $H\rightarrow H$, 
$B:$ $H\rightarrow L_{HS}(U_0,H)$), linear operator
$A:$ $D(A)\subset H\rightarrow H$
as well as the initial value $\xi$ 
are assumed to satisfy the conditions
of existence and uniqueness of
the SPDE (\ref{xx1}) mild solution \cite{23}
(see also \cite{12}, \cite{13}).
Here $U_0$ is an $\mathbb{R}$-Hilbert space 
defined by $U_0=Q^{1/2}(U).$ The scalar product in $U_0$ is defined as follows
$\left\langle u, w\right\rangle_{U_0}=\left\langle 
Q^{-1/2}u, Q^{-1/2}w\right\rangle_{U}$ for all $u, w\in U_0.$

As it is known, strong numerical methods with high-orders of accuracy 
(with respect to the temporal discretization)
for approximating the mild solution of the SPDE (\ref{xx1}),
which are based on the Taylor formula in Banach spaces and an 
exponential formula for the mild solution of 
SPDEs, contain iterated stochastic integrals with respect to the $Q$-Wiener 
process \cite{8}, \cite{10}-\cite{13}, \cite{18}.

Note that the exponential Milstein type numerical scheme 
\cite{12}, \cite{18}, \cite{25} 
and exponential Wagner--Platen type numerical scheme \cite{13}
contain, for example, the following 
iterated stochastic integrals

\begin{equation}
\label{xx2}
\int\limits_{t}^{T}B(Z) d{\bf W}_{t_1},\ \ 
\int\limits_{t}^{T}B'(Z) \left(
\int\limits_{t}^{t_2}B(Z) d{\bf W}_{t_1} \right) d{\bf W}_{t_2},\
\end{equation}

\begin{equation}
\label{xx3}
\int\limits_{t}^{T}F'(Z)
\left(\int\limits_{t}^{t_2} B(Z) d{\bf W}_{t_1} \right) dt_2,\ \ 
\int\limits_{t}^{T}B'(Z) \left(
\int\limits_{t}^{t_3}B'(Z) \left(
\int\limits_{t}^{t_2}B(Z)
d{\bf W}_{t_1} \right) d{\bf W}_{t_2} \right) d{\bf W}_{t_3},\
\end{equation}

\begin{equation}
\label{xx3aaa}
\int\limits_{t}^{T}B'(Z) \left(
\int\limits_{t}^{t_2}F(Z) dt_1 \right) d{\bf W}_{t_2},\ \
\int\limits_{t}^{T}B''(Z) \left(
\int\limits_{t}^{t_2}B(Z) d{\bf W}_{t_1},
\int\limits_{t}^{t_2}B(Z) d{\bf W}_{t_1} 
\right) d{\bf W}_{t_2},\
\end{equation}

\vspace{4mm}
\noindent
where $0\le t<T \le \bar{T},$\ 
$Z: \Omega \rightarrow H$ is an ${\bf F}_t/{\mathcal{B}}(H)$-measurable mapping
and $F',$ $B',$ $B''$ denote 
Fr\^{e}chet derivatives.
At that, the exponential Milstein type scheme \cite{12} contains integrals 
(\ref{xx2}) while the exponential Wagner--Platen 
type scheme \cite{13} 
contains integrals (\ref{xx2})--(\ref{xx3aaa}).
It is easy to notice that the numerical schemes for SPDEs with higher 
orders of convergence (with respect to the temporal discretization) 
in contrast with numerical schemes from 
\cite{12}, \cite{13} will include iterated stochastic 
integrals (with respect to the $Q$-Wiener process) 
with multiplicities $k>3$ \cite{22} (2012). 
So, this work is partially devoted to the approximation of iterated 
stochastic integrals of the form

\begin{equation}
\label{xx4}
I[\Phi^{(k)}(Z)]_{T,t}=\int\limits_{t}^{T}\Phi_k(Z) \left( \ldots \left(
\int\limits_{t}^{t_3}\Phi_2(Z) \left(
\int\limits_{t}^{t_2}\Phi_1(Z)
d{\bf W}_{t_1} \right)
d{\bf W}_{t_2} \right) \ldots  \right) d{\bf W}_{t_k},\
\end{equation}

\vspace{4mm}
\noindent
where 
$Z: \Omega \rightarrow H$ is an ${\bf F}_t/{\mathcal{B}}(H)$-measurable 
mapping, 
$\Phi_k(v)(\ \ldots (\Phi_2(v)(\Phi_1(v)) \ldots\ ))$ 
is a $k$-linear Hilbert--Schmidt operator mapping from
$\underbrace{U_0\times \ldots \times U_0}_{\small{\hbox{$k$ times}}}$ to $H$
for all $v\in H,$
and $0\le t<T \le \bar{T}.$ 
In Sect.~5 we consider the approximation of more general
iterated stochastic integrals than (\ref{xx4}). 
In Sect.~6,~7 some other
types of ite\-ra\-ted stochastic integrals of multiplicities
2--4 with respect to the $Q$-Wiener process will be considered.
In this paper, in all the integrals mentioned above, the
infinite-dimensional $Q$-Wiener process will be replaced by its 
finite-dimensional approximation. 
In \cite{48aa}-\cite{new-new-1}, 
(also see \cite{37a}, Chapter 7)
one can find a continuation of the studies begun in this work.
In \cite{37a}, \cite{48aa}-\cite{new-new-1} we consider the approximation of 
iterated stochastic integrals (\ref{xx2})--(\ref{xx3aaa})
with respect to the infinite-dimensional $Q$-Wiener process.

Note that the second stochastic integral in (\ref{xx3aaa}) is not a special case 
of the stochastic integral (\ref{xx4}) for $k=3$.
Nevertheless, the expanded representation of the approximation of stochastic 
integral (\ref{xx3aaa}) has a 
close structure to (\ref{xx405}) for $k=3$ (see below). Moreover,
the mentioned
representation of stochastic integral (\ref{xx3aaa})
contains the same iterated It\^{o} stochastic integrals 
of third multiplicity as in (\ref{xx405}) for $k=3$
(see Sect.~6).
These conclusions mean that the main result of this article (Theorem 4, Sect.~5) 
for $k=3$ can be reformulated naturally for the stochastic integral 
(\ref{xx3aaa}) (see Sect.~6).

It should be noted that by developing an approach from the work \cite{13}, 
which uses the Taylor formula in Banach spaces and a formula for 
the mild solution of 
the SPDE (\ref{xx1}),
we obviously obtain a number of other iterated stochastic 
integrals with respect to the $Q$-Wiener process. 
For example, the following stochastic integrals

\vspace{2mm}
$$
\int\limits_{t}^{T}B'''(Z) \left(
\int\limits_{t}^{t_2}B(Z) d{\bf W}_{t_1},
\int\limits_{t}^{t_2}B(Z) d{\bf W}_{t_1},
\int\limits_{t}^{t_2}B(Z) d{\bf W}_{t_1}
\right) d{\bf W}_{t_2},\
$$

\vspace{1mm}
$$
\int\limits_{t}^{T}B'(Z) \left(
\int\limits_{t}^{t_3}B''(Z) \left(
\int\limits_{t}^{t_2}B(Z) d{\bf W}_{t_1},
\int\limits_{t}^{t_2}B(Z) d{\bf W}_{t_1} 
\right) d{\bf W}_{t_2}\right)d{\bf W}_{t_3},
$$

\vspace{1mm}
$$
\int\limits_{t}^{T}B''(Z) \left(
\int\limits_{t}^{t_3}B(Z)d{\bf W}_{t_1},
\int\limits_{t}^{t_3}B'(Z)\left(
\int\limits_{t}^{t_2}B(Z) d{\bf W}_{t_1}
\right) d{\bf W}_{t_2}\right)d{\bf W}_{t_3},\
$$

\vspace{1mm}
$$
\int\limits_{t}^{T}F'(Z)\left(
\int\limits_{t}^{t_3}B'(Z) \left(
\int\limits_{t}^{t_2}B(Z) d{\bf W}_{t_1}
\right) d{\bf W}_{t_2} 
\right) dt_3,
$$

\vspace{1mm}
$$
\int\limits_{t}^{T}F''(Z) \left(
\int\limits_{t}^{t_2}B(Z) d{\bf W}_{t_1},
\int\limits_{t}^{t_2}B(Z) d{\bf W}_{t_1} 
\right) dt_2,
$$

\vspace{1mm}
$$
\int\limits_{t}^{T}B''(Z) \left(
\int\limits_{t}^{t_2}F(Z) dt_1,
\int\limits_{t}^{t_2}B(Z) d{\bf W}_{t_1} 
\right) d{\bf W}_{t_2}
$$

\vspace{5mm}
\noindent
will be considered in Sect.~7. 
Here 
$Z: \Omega \rightarrow H$ is an ${\bf F}_t/{\mathcal{B}}(H)$-measurable 
mapping and 
$B'$, $B'',$ $B''',$ $F'$, $F''$ are
Fr\^{e}chet derivatives.

Consider eigenvalues $\lambda_i$ and 
eigenfunctions $e_i(x)$ of the covariance operator $Q,$ where
$i=(i_1,\ldots,i_d)$ $\in$ $J,$
$x=(x_1,\ldots,x_d),$
and 
$J=\{i:\ i\in \mathbb{N}^d,\ \hbox{and}\ \lambda_i>0\}$.

The series representation of 
the $Q$-Wiener process has the following form \cite{20}

\vspace{1mm}
$$
{\bf W}(t,x)=\sum\limits_{i\in J} e_i(x)\sqrt{\lambda_i}
{\bf w}_t^{(i)},\ \ \ t\in[0, \bar{T}],
$$

\vspace{2mm}
\noindent
or in the shorter notations

$$
{\bf W}_t=\sum\limits_{i\in J} e_i\sqrt{\lambda_i}
{\bf w}_t^{(i)},\ \ \ t\in[0, \bar{T}],
$$

\vspace{2mm}
\noindent
where ${\bf w}_t^{(i)},$ $i\in J$ are independent standard
Wiener processes.
Note that eigenfunctions $e_i,$ $i\in J$ form
an orthonormal basis of $U$ \cite{20}.

Consider the finite-dimensional approximation of ${\bf W}_t$ \cite{20}

\begin{equation}
\label{yy.1}
{\bf W}^M_t=\sum\limits_{i\in J_M} e_i\sqrt{\lambda_i}
{\bf w}_t^{(i)},\ \ \ t\in[0, \bar{T}],
\end{equation}

\vspace{2mm}
\noindent
where $J_M=\{i:\ 1\le i_1,\ldots,i_d\le M,\ \hbox{and}\ 
\lambda_i>0\}$.

Using (\ref{yy.1}) and the relation \cite{20}

\begin{equation}
\label{xx.1}
{\bf w}_t^{(i)}=\frac{1}{\sqrt{\lambda_i}}\langle e_i,{\bf W}_t\rangle_U,\ \ \
i\in J,
\end{equation}

\vspace{2mm}
\noindent
we obtain

\begin{equation}
\label{yy.12}
{\bf W}^M_t=\sum\limits_{i\in J_M} 
e_i\ \langle e_i,{\bf W}_t\rangle_U,\ \ \ t\in[0, \bar{T}],
\end{equation}

\vspace{2mm}
\noindent
where $\langle \cdot,\cdot\rangle_U$ is a scalar product in $U.$

Taking into account (\ref{xx.1}), (\ref{yy.12}), we note that
the approximation $I[\Phi^{(k)}(Z)]_{T,t}^M$ of ite\-ra\-ted stochastic integral
$I[\Phi^{(k)}(Z)]_{T,t}$ (see (\ref{xx4}))
can be rewritten with probability 1 (further w. p. 1) in the following form

\vspace{2mm}
$$
I[\Phi^{(k)}(Z)]_{T,t}^M=\int\limits_{t}^{T}\Phi_k(Z) \left( \ldots \left(
\int\limits_{t}^{t_3}\Phi_2(Z) \left(
\int\limits_{t}^{t_2}\Phi_1(Z)
d{\bf W}_{t_1}^M\ \right) 
d{\bf W}_{t_2}^M\ \right) \ldots  \right) d{\bf W}_{t_k}^M=
$$

\vspace{4mm}
$$
=
\sum_{r_1,\ldots,r_k\in J_M}
\Phi_k(Z)\left(\ldots \left(\Phi_2(Z) \left(\Phi_1(Z)
e_{r_1} \right)
e_{r_2} \right) \ldots \right) e_{r_k} \times
$$

$$
\times
\int\limits_t^T \ldots \int\limits_t^{t_{3}}\int\limits_t^{t_{2}}
d\langle e_{r_1},{\bf W}_{t_1}
\rangle_U\
d\langle e_{r_2},{\bf W}_{t_2}
\rangle_U\ \ldots \
d\langle e_{r_k},{\bf W}_{t_k}
\rangle_U =
$$

\vspace{4mm}
$$
=\sum_{r_1,\ldots,r_k\in J_M}
\Phi_k(Z)\left(\ldots \left(\Phi_2(Z) \left(\Phi_1(Z)
e_{r_1} \right) 
e_{r_2} \right) \ldots \right) e_{r_k}
\sqrt{\lambda_{r_1}\lambda_{r_2}\ldots \lambda_{r_k}}\times
$$

\begin{equation}
\label{xx405}
\times
\int\limits_t^T \ldots \int\limits_t^{t_3}\int\limits_t^{t_{2}}
d{\bf w}_{t_1}^{(r_1)}d{\bf w}_{t_2}^{(r_2)}\ldots
d{\bf w}_{t_k}^{(r_k)},
\end{equation}

\vspace{5mm}
\noindent
where $0\le t<T \le \bar{T}.$

\vspace{2mm}

{\bf Remark 1.}\ {\it Obviously, without the loss of generality
we can write $J_M=\{1, 2,\ldots, M\}.$}

\vspace{2mm}

When special conditions of commutativity 
for SPDEs in the form (\ref{xx1}) be fulfilled 
it is 
proposed to simulate 
numerically
the stochastic integrals (\ref{xx2})--(\ref{xx3aaa}) using 
the simple formulas \cite{12}, \cite{13}. 
In this case, the numerical simulation of mentioned stochastic 
integrals requires the use of increments of the $Q$-Wiener process only.
However, if these commutativity conditions are not fulfilled
(which often corresponds to SPDEs in numerous applications), the 
numerical simulation of stochastic integrals (\ref{xx2})--(\ref{xx3aaa})
becomes 
much more difficult. 
In \cite{26} two methods for the mean-square approximation of simplest
iterated (double) 
stochastic integrals with respect to 
the $Q$-Wiener process are proposed. In this article, we consider 
a substantially more general and effective method
for the mean-square approximation of 
iterated stochastic integrals of multiplicity $k$ $(k\in\mathbb{N})$
with respect to 
the $Q$-Wiener process.
The convergence analysis in the transition from $J_M$ to $J$, i.e.
from ${\bf W}^M_t$ to ${\bf W}t$
is carried out in \cite {37a} (Sect.7.4.2), \cite {37aa} (Sect.7.4.2), 
\cite{37aaxx}, \cite{48aa}, \cite{48aaa}
for stochastic integrals of multiplicity $k$ ($k=1, 2, 3$) with respect to the 
$Q$-Wiener process (the cases $k=1, 2$ is considered in Theorem 1 from \cite{26}).

The monographs \cite{37} (Chapters 5 and 6) and \cite{37a} or \cite{37aa},
\cite{37aaxx} (Chapters 1, 2, and 5) (also
see \cite{30}-\cite{37pred}, \cite{38}-\cite{48a})
are devoted to constructing of efficient methods
of the mean-square approximation of iterated It\^{o} stochastic
integrals with respect to the scalar standard
Wiener processes. 
These results are adapted for iterated Stratonovich stochastic integrals
\cite{eee}-\cite{48a}.
Below (Sect.~2--4) we consider a very short review of results from 
monographs \cite{37} (Chapters 5 and 6) and \cite{37a} or \cite{37aa},
\cite{37aaxx} (Chapters 1, 2, and 5) 
and some new results (Sect.~5--7).

\vspace{3mm}

\section{Method of Approximation of Iterated 
It\^{o} Stochastic integrals Based on Generalized Multiple Fourier Series}

\vspace{3mm}

Consider more general 
iterated It\^{o} stochastic integrals than in (\ref{xx405})

\begin{equation}
\label{sodom20}
J[\psi^{(k)}]_{T,t}^{(i_1\ldots i_k)}
=\int\limits_t^T\psi_k(t_k) \ldots \int\limits_t^{t_{2}}
\psi_1(t_1) d{\bf w}_{t_1}^{(i_1)}\ldots
d{\bf w}_{t_k}^{(i_k)},
\end{equation}

\vspace{2mm}
\noindent
where $0\le t<T \le \bar{T}$ and every $\psi_l(\tau)\ (l=1,\ldots,k)$ is 
a continuous non-random function on $[t, T]$;
${\bf w}_{\tau}^{(i)}$ ($i=1,\ldots,m)$ are independent
standard Wiener processes (see Sect.~1) and
${\bf w}_{\tau}^{(0)}=\tau;$\
$i_1,\ldots,i_k=0,\ 1,\ldots,m.$ The case $\psi_1(\tau),\ldots,\psi_k(\tau)\in L_2([t, T])$
will be considered in Theorem~2 (see below).

Suppose that $\{\phi_j(x)\}_{j=0}^{\infty}$
is a complete orthonormal system of functions in 
$L_2([t, T])$.
Define the following function on the hypercube $[t, T]^k$

\begin{equation}
\label{ppp}
K(t_1,\ldots,t_k)=
\begin{cases}
\psi_1(t_1)\ldots \psi_k(t_k),\ t_1<\ldots<t_k\\
~\\
~\\
0,\ \hbox{\rm otherwise}
\end{cases}
=\ \
\prod\limits_{l=1}^k
\psi_l(t_l)\ \prod\limits_{l=1}^{k-1}{\bf 1}_{\{t_l<t_{l+1}\}},\ 
\end{equation}

\vspace{3mm}
\noindent
where $t_1,\ldots,t_k\in [t, T]$ for $k\ge 2$ and 
$K(t_1)\equiv\psi_1(t_1)$ for $t_1\in[t, T].$ Here 
${\bf 1}_A$ is the indicator of the set $A$.

The function $K(t_1,\ldots,t_k)$ is piecewise continuous on the 
hypercube $[t, T]^k.$
At this situation it is well known that the generalized 
multiple Fourier series 
of $K(t_1,\ldots,t_k)\in L_2([t, T]^k)$ converges
to $K(t_1,\ldots,t_k)$ in the hypercube $[t, T]^k$ in 
the mean-square sense, i.e.

\begin{equation}
\label{sos1z}
\hbox{\vtop{\offinterlineskip\halign{
\hfil#\hfil\cr
{\rm lim}\cr
$\stackrel{}{{}_{p_1,\ldots,p_k\to \infty}}$\cr
}} }\Biggl\Vert
K(t_1,\ldots,t_k)-
\sum_{j_1=0}^{p_1}\ldots \sum_{j_k=0}^{p_k}
C_{j_k\ldots j_1}\prod_{l=1}^{k} \phi_{j_l}(t_l)\Biggr
\Vert_{L_2([t, T]^k)}=0,
\end{equation}

\vspace{2mm}
\noindent
where
\begin{equation}
\label{ppppa}
C_{j_k\ldots j_1}=\int\limits_{[t,T]^k}
K(t_1,\ldots,t_k)\prod_{l=1}^{k}\phi_{j_l}(t_l)dt_1\ldots dt_k
\end{equation}

\vspace{3mm}
\noindent
is the Fourier coefficient and
$$
\left\Vert f\right\Vert_{L_2([t, T]^k)}=\left(\int\limits_{[t,T]^k}
f^2(t_1,\ldots,t_k)dt_1\ldots dt_k\right)^{1/2}.
$$

\vspace{4mm}

Consider the discretization $\{\tau_j\}_{j=0}^N$ of $[t,T]$ such that

\vspace{-1mm}
\begin{equation}
\label{1111}
t=\tau_0<\ldots <\tau_N=T,\ \ \ \
\Delta_N=
\hbox{\vtop{\offinterlineskip\halign{
\hfil#\hfil\cr
{\rm max}\cr
$\stackrel{}{{}_{0\le j\le N-1}}$\cr
}} }\Delta\tau_j\to 0\ \ \hbox{if}\ \ N\to \infty,\ \ \ \
\Delta\tau_j=\tau_{j+1}-\tau_j.
\end{equation}

\vspace{4mm}

{\bf Theorem 1} \cite{30} (2006) (also see \cite{30a}-\cite{48aaa}). {\it 
Suppose that
every $\psi_l(\tau)$ $(l=1,\ldots, k)$ is a conti\-nu\-ous non-random function 
on $[t, T]$ and
$\{\phi_j(x)\}_{j=0}^{\infty}$ is a complete orthonormal system  
of continuous functions in $L_2([t,T]).$ Then

\vspace{1mm}
$$
J[\psi^{(k)}]_{T,t}^{(i_1\ldots i_k)}=
\hbox{\vtop{\offinterlineskip\halign{
\hfil#\hfil\cr
{\rm l.i.m.}\cr
$\stackrel{}{{}_{p_1,\ldots,p_k\to \infty}}$\cr
}} }\sum_{j_1=0}^{p_1}\ldots\sum_{j_k=0}^{p_k}
C_{j_k\ldots j_1}\Biggl(
\prod_{l=1}^k\zeta_{j_l}^{(i_l)}-
\Biggr.
$$

\vspace{3mm}
\begin{equation}
\label{tyyy}
-\ \Biggl.
\hbox{\vtop{\offinterlineskip\halign{
\hfil#\hfil\cr
{\rm l.i.m.}\cr
$\stackrel{}{{}_{N\to \infty}}$\cr
}} }\sum_{(l_1,\ldots,l_k)\in {\rm G}_k}
\phi_{j_{1}}(\tau_{l_1})
\Delta{\bf w}_{\tau_{l_1}}^{(i_1)}\ldots
\phi_{j_{k}}(\tau_{l_k})
\Delta{\bf w}_{\tau_{l_k}}^{(i_k)}\Biggr),
\end{equation}

\vspace{3mm}
\noindent
where

\vspace{-1mm}
$$
{\rm G}_k={\rm H}_k\backslash{\rm L}_k;\
{\rm H}_k=\{(l_1,\ldots,l_k):\ l_1,\ldots,l_k=0,\ 1,\ldots,N-1\},
$$

\vspace{-2mm}
$$
{\rm L}_k=\{(l_1,\ldots,l_k):\ l_1,\ldots,l_k=0,\ 1,\ldots,N-1;\
l_g\ne l_r\ (g\ne r);\ g, r=1,\ldots,k\},
$$

\vspace{3mm}
\noindent
${\rm l.i.m.}$ is a limit in the mean-square sense{\rm ,} 
$i_1,\ldots,i_k=0,1,\ldots,m,$

\begin{equation}
\label{rr23}
\zeta_{j}^{(i)}=
\int\limits_t^T \phi_{j}(s) d{\bf w}_s^{(i)}
\end{equation} 

\vspace{2mm}
\noindent
are independent standard Gaussian random variables
for various
$i$\ or $j$ {\rm(}if $i\ne 0${\rm),}
$C_{j_k\ldots j_1}$ is the Fourier coefficient {\rm(\ref{ppppa}),}
$\Delta{\bf w}_{\tau_{j}}^{(i)}=
{\bf w}_{\tau_{j+1}}^{(i)}-{\bf w}_{\tau_{j}}^{(i)}$
$(i=0,\ 1,\ldots,m),$\
$\left\{\tau_{j}\right\}_{j=0}^{N}$ 
is the discretization of
$[t,T],$ which satisfies the condition {\rm (\ref{1111})}.
}

\vspace{2mm}

Note that in \cite{30}-\cite{48} the version of Theorem 1 for
systems of Haar and Rademacher--Walsh 
functions has been considered. Another modifications and generalizations 
of Theorem 1 can be found in the monographs \cite{37a}-\cite{37aaxx}
(also see Theorem~2 below).

It is not difficult to see that for the case of pairwise different numbers
$i_1,\ldots,i_k=1,\ldots,m$ from Theorem 1 we obtain
$$
J[\psi^{(k)}]_{T,t}^{(i_1\ldots i_k)}=
\hbox{\vtop{\offinterlineskip\halign{
\hfil#\hfil\cr
{\rm l.i.m.}\cr
$\stackrel{}{{}_{p_1,\ldots,p_k\to \infty}}$\cr
}} }\sum_{j_1=0}^{p_1}\ldots\sum_{j_k=0}^{p_k}
C_{j_k\ldots j_1}\zeta_{j_1}^{(i_1)}\ldots \zeta_{j_k}^{(i_k)}.
$$

\vspace{3mm}

In order to evaluate the significance of Theorem 1 for practice we will
demonstrate its transformed particular cases for 
$k=1,\ldots,6$ \cite{30}-\cite{48a}
(the cases $k=7$ and $k>7$
can be found in \cite{32}, \cite{34}, \cite{37}-\cite{37aaxx})
\begin{equation}
\label{a1}
J[\psi^{(1)}]_{T,t}^{(i_1)}
=\hbox{\vtop{\offinterlineskip\halign{
\hfil#\hfil\cr
{\rm l.i.m.}\cr
$\stackrel{}{{}_{p_1\to \infty}}$\cr
}} }\sum_{j_1=0}^{p_1}
C_{j_1}\zeta_{j_1}^{(i_1)},
\end{equation}

\vspace{2mm}
\begin{equation}
\label{a2}
J[\psi^{(2)}]_{T,t}^{(i_1 i_2)}
=\hbox{\vtop{\offinterlineskip\halign{
\hfil#\hfil\cr
{\rm l.i.m.}\cr
$\stackrel{}{{}_{p_1,p_2\to \infty}}$\cr
}} }\sum_{j_1=0}^{p_1}\sum_{j_2=0}^{p_2}
C_{j_2j_1}\Biggl(\zeta_{j_1}^{(i_1)}\zeta_{j_2}^{(i_2)}
-{\bf 1}_{\{i_1=i_2\ne 0\}}
{\bf 1}_{\{j_1=j_2\}}\Biggr),
\end{equation}

\vspace{5mm}
$$
J[\psi^{(3)}]_{T,t}^{(i_1 i_2 i_3)}=
\hbox{\vtop{\offinterlineskip\halign{
\hfil#\hfil\cr
{\rm l.i.m.}\cr
$\stackrel{}{{}_{p_1,p_2,p_3\to \infty}}$\cr
}} }\sum_{j_1=0}^{p_1}\sum_{j_2=0}^{p_2}\sum_{j_3=0}^{p_3}
C_{j_3j_2j_1}\Biggl(
\zeta_{j_1}^{(i_1)}\zeta_{j_2}^{(i_2)}\zeta_{j_3}^{(i_3)}
-\Biggr.
$$
\begin{equation}
\label{a3}
\Biggl.
-{\bf 1}_{\{i_1=i_2\ne 0\}}
{\bf 1}_{\{j_1=j_2\}}
\zeta_{j_3}^{(i_3)}
-{\bf 1}_{\{i_2=i_3\ne 0\}}
{\bf 1}_{\{j_2=j_3\}}
\zeta_{j_1}^{(i_1)}-{\bf 1}_{\{i_1=i_3\ne 0\}}
{\bf 1}_{\{j_1=j_3\}}
\zeta_{j_2}^{(i_2)}\Biggr),
\end{equation}

\vspace{5mm}

$$
J[\psi^{(4)}]_{T,t}^{(i_1\ldots i_4)}
=
\hbox{\vtop{\offinterlineskip\halign{
\hfil#\hfil\cr
{\rm l.i.m.}\cr
$\stackrel{}{{}_{p_1,\ldots,p_4\to \infty}}$\cr
}} }\sum_{j_1=0}^{p_1}\ldots\sum_{j_4=0}^{p_4}
C_{j_4\ldots j_1}\Biggl(
\prod_{l=1}^4\zeta_{j_l}^{(i_l)}
\Biggr.
-
$$
$$
-
{\bf 1}_{\{i_1=i_2\ne 0\}}
{\bf 1}_{\{j_1=j_2\}}
\zeta_{j_3}^{(i_3)}
\zeta_{j_4}^{(i_4)}
-
{\bf 1}_{\{i_1=i_3\ne 0\}}
{\bf 1}_{\{j_1=j_3\}}
\zeta_{j_2}^{(i_2)}
\zeta_{j_4}^{(i_4)}-
$$
$$
-
{\bf 1}_{\{i_1=i_4\ne 0\}}
{\bf 1}_{\{j_1=j_4\}}
\zeta_{j_2}^{(i_2)}
\zeta_{j_3}^{(i_3)}
-
{\bf 1}_{\{i_2=i_3\ne 0\}}
{\bf 1}_{\{j_2=j_3\}}
\zeta_{j_1}^{(i_1)}
\zeta_{j_4}^{(i_4)}-
$$
$$
-
{\bf 1}_{\{i_2=i_4\ne 0\}}
{\bf 1}_{\{j_2=j_4\}}
\zeta_{j_1}^{(i_1)}
\zeta_{j_3}^{(i_3)}
-
{\bf 1}_{\{i_3=i_4\ne 0\}}
{\bf 1}_{\{j_3=j_4\}}
\zeta_{j_1}^{(i_1)}
\zeta_{j_2}^{(i_2)}+
$$
$$
+
{\bf 1}_{\{i_1=i_2\ne 0\}}
{\bf 1}_{\{j_1=j_2\}}
{\bf 1}_{\{i_3=i_4\ne 0\}}
{\bf 1}_{\{j_3=j_4\}}
+
{\bf 1}_{\{i_1=i_3\ne 0\}}
{\bf 1}_{\{j_1=j_3\}}
{\bf 1}_{\{i_2=i_4\ne 0\}}
{\bf 1}_{\{j_2=j_4\}}+
$$
\begin{equation}
\label{a4}
+\Biggl.
{\bf 1}_{\{i_1=i_4\ne 0\}}
{\bf 1}_{\{j_1=j_4\}}
{\bf 1}_{\{i_2=i_3\ne 0\}}
{\bf 1}_{\{j_2=j_3\}}\Biggr),
\end{equation}

\vspace{5mm}

$$
J[\psi^{(5)}]_{T,t}^{(i_1\ldots i_5)}
=\hbox{\vtop{\offinterlineskip\halign{
\hfil#\hfil\cr
{\rm l.i.m.}\cr
$\stackrel{}{{}_{p_1,\ldots,p_5\to \infty}}$\cr
}} }\sum_{j_1=0}^{p_1}\ldots\sum_{j_5=0}^{p_5}
C_{j_5\ldots j_1}\Biggl(
\prod_{l=1}^5\zeta_{j_l}^{(i_l)}
-\Biggr.
$$
$$
-
{\bf 1}_{\{i_1=i_2\ne 0\}}
{\bf 1}_{\{j_1=j_2\}}
\zeta_{j_3}^{(i_3)}
\zeta_{j_4}^{(i_4)}
\zeta_{j_5}^{(i_5)}-
{\bf 1}_{\{i_1=i_3\ne 0\}}
{\bf 1}_{\{j_1=j_3\}}
\zeta_{j_2}^{(i_2)}
\zeta_{j_4}^{(i_4)}
\zeta_{j_5}^{(i_5)}-
$$
$$
-
{\bf 1}_{\{i_1=i_4\ne 0\}}
{\bf 1}_{\{j_1=j_4\}}
\zeta_{j_2}^{(i_2)}
\zeta_{j_3}^{(i_3)}
\zeta_{j_5}^{(i_5)}-
{\bf 1}_{\{i_1=i_5\ne 0\}}
{\bf 1}_{\{j_1=j_5\}}
\zeta_{j_2}^{(i_2)}
\zeta_{j_3}^{(i_3)}
\zeta_{j_4}^{(i_4)}-
$$
$$
-
{\bf 1}_{\{i_2=i_3\ne 0\}}
{\bf 1}_{\{j_2=j_3\}}
\zeta_{j_1}^{(i_1)}
\zeta_{j_4}^{(i_4)}
\zeta_{j_5}^{(i_5)}-
{\bf 1}_{\{i_2=i_4\ne 0\}}
{\bf 1}_{\{j_2=j_4\}}
\zeta_{j_1}^{(i_1)}
\zeta_{j_3}^{(i_3)}
\zeta_{j_5}^{(i_5)}-
$$
$$
-
{\bf 1}_{\{i_2=i_5\ne 0\}}
{\bf 1}_{\{j_2=j_5\}}
\zeta_{j_1}^{(i_1)}
\zeta_{j_3}^{(i_3)}
\zeta_{j_4}^{(i_4)}
-{\bf 1}_{\{i_3=i_4\ne 0\}}
{\bf 1}_{\{j_3=j_4\}}
\zeta_{j_1}^{(i_1)}
\zeta_{j_2}^{(i_2)}
\zeta_{j_5}^{(i_5)}-
$$
$$
-
{\bf 1}_{\{i_3=i_5\ne 0\}}
{\bf 1}_{\{j_3=j_5\}}
\zeta_{j_1}^{(i_1)}
\zeta_{j_2}^{(i_2)}
\zeta_{j_4}^{(i_4)}
-{\bf 1}_{\{i_4=i_5\ne 0\}}
{\bf 1}_{\{j_4=j_5\}}
\zeta_{j_1}^{(i_1)}
\zeta_{j_2}^{(i_2)}
\zeta_{j_3}^{(i_3)}+
$$
$$
+
{\bf 1}_{\{i_1=i_2\ne 0\}}
{\bf 1}_{\{j_1=j_2\}}
{\bf 1}_{\{i_3=i_4\ne 0\}}
{\bf 1}_{\{j_3=j_4\}}\zeta_{j_5}^{(i_5)}+
{\bf 1}_{\{i_1=i_2\ne 0\}}
{\bf 1}_{\{j_1=j_2\}}
{\bf 1}_{\{i_3=i_5\ne 0\}}
{\bf 1}_{\{j_3=j_5\}}\zeta_{j_4}^{(i_4)}+
$$
$$
+
{\bf 1}_{\{i_1=i_2\ne 0\}}
{\bf 1}_{\{j_1=j_2\}}
{\bf 1}_{\{i_4=i_5\ne 0\}}
{\bf 1}_{\{j_4=j_5\}}\zeta_{j_3}^{(i_3)}+
{\bf 1}_{\{i_1=i_3\ne 0\}}
{\bf 1}_{\{j_1=j_3\}}
{\bf 1}_{\{i_2=i_4\ne 0\}}
{\bf 1}_{\{j_2=j_4\}}\zeta_{j_5}^{(i_5)}+
$$
$$
+
{\bf 1}_{\{i_1=i_3\ne 0\}}
{\bf 1}_{\{j_1=j_3\}}
{\bf 1}_{\{i_2=i_5\ne 0\}}
{\bf 1}_{\{j_2=j_5\}}\zeta_{j_4}^{(i_4)}+
{\bf 1}_{\{i_1=i_3\ne 0\}}
{\bf 1}_{\{j_1=j_3\}}
{\bf 1}_{\{i_4=i_5\ne 0\}}
{\bf 1}_{\{j_4=j_5\}}\zeta_{j_2}^{(i_2)}+
$$
$$
+
{\bf 1}_{\{i_1=i_4\ne 0\}}
{\bf 1}_{\{j_1=j_4\}}
{\bf 1}_{\{i_2=i_3\ne 0\}}
{\bf 1}_{\{j_2=j_3\}}\zeta_{j_5}^{(i_5)}+
{\bf 1}_{\{i_1=i_4\ne 0\}}
{\bf 1}_{\{j_1=j_4\}}
{\bf 1}_{\{i_2=i_5\ne 0\}}
{\bf 1}_{\{j_2=j_5\}}\zeta_{j_3}^{(i_3)}+
$$
$$
+
{\bf 1}_{\{i_1=i_4\ne 0\}}
{\bf 1}_{\{j_1=j_4\}}
{\bf 1}_{\{i_3=i_5\ne 0\}}
{\bf 1}_{\{j_3=j_5\}}\zeta_{j_2}^{(i_2)}+
{\bf 1}_{\{i_1=i_5\ne 0\}}
{\bf 1}_{\{j_1=j_5\}}
{\bf 1}_{\{i_2=i_3\ne 0\}}
{\bf 1}_{\{j_2=j_3\}}\zeta_{j_4}^{(i_4)}+
$$
$$
+
{\bf 1}_{\{i_1=i_5\ne 0\}}
{\bf 1}_{\{j_1=j_5\}}
{\bf 1}_{\{i_2=i_4\ne 0\}}
{\bf 1}_{\{j_2=j_4\}}\zeta_{j_3}^{(i_3)}+
{\bf 1}_{\{i_1=i_5\ne 0\}}
{\bf 1}_{\{j_1=j_5\}}
{\bf 1}_{\{i_3=i_4\ne 0\}}
{\bf 1}_{\{j_3=j_4\}}\zeta_{j_2}^{(i_2)}+
$$
$$
+
{\bf 1}_{\{i_2=i_3\ne 0\}}
{\bf 1}_{\{j_2=j_3\}}
{\bf 1}_{\{i_4=i_5\ne 0\}}
{\bf 1}_{\{j_4=j_5\}}\zeta_{j_1}^{(i_1)}+
{\bf 1}_{\{i_2=i_4\ne 0\}}
{\bf 1}_{\{j_2=j_4\}}
{\bf 1}_{\{i_3=i_5\ne 0\}}
{\bf 1}_{\{j_3=j_5\}}\zeta_{j_1}^{(i_1)}+
$$
\begin{equation}
\label{a5}
+\Biggl.
{\bf 1}_{\{i_2=i_5\ne 0\}}
{\bf 1}_{\{j_2=j_5\}}
{\bf 1}_{\{i_3=i_4\ne 0\}}
{\bf 1}_{\{j_3=j_4\}}\zeta_{j_1}^{(i_1)}\Biggr),
\end{equation}

\vspace{5mm}
$$
J[\psi^{(6)}]_{T,t}^{(i_1\ldots i_6)}
=\hbox{\vtop{\offinterlineskip\halign{
\hfil#\hfil\cr
{\rm l.i.m.}\cr
$\stackrel{}{{}_{p_1,\ldots,p_6\to \infty}}$\cr
}} }\sum_{j_1=0}^{p_1}\ldots\sum_{j_6=0}^{p_6}
C_{j_6\ldots j_1}\Biggl(
\prod_{l=1}^6
\zeta_{j_l}^{(i_l)}
-\Biggr.
$$
$$
-
{\bf 1}_{\{i_1=i_6\ne 0\}}
{\bf 1}_{\{j_1=j_6\}}
\zeta_{j_2}^{(i_2)}
\zeta_{j_3}^{(i_3)}
\zeta_{j_4}^{(i_4)}
\zeta_{j_5}^{(i_5)}-
{\bf 1}_{\{i_2=i_6\ne 0\}}
{\bf 1}_{\{j_2=j_6\}}
\zeta_{j_1}^{(i_1)}
\zeta_{j_3}^{(i_3)}
\zeta_{j_4}^{(i_4)}
\zeta_{j_5}^{(i_5)}-
$$
$$
-
{\bf 1}_{\{i_3=i_6\ne 0\}}
{\bf 1}_{\{j_3=j_6\}}
\zeta_{j_1}^{(i_1)}
\zeta_{j_2}^{(i_2)}
\zeta_{j_4}^{(i_4)}
\zeta_{j_5}^{(i_5)}-
{\bf 1}_{\{i_4=i_6\ne 0\}}
{\bf 1}_{\{j_4=j_6\}}
\zeta_{j_1}^{(i_1)}
\zeta_{j_2}^{(i_2)}
\zeta_{j_3}^{(i_3)}
\zeta_{j_5}^{(i_5)}-
$$
$$
-
{\bf 1}_{\{i_5=i_6\ne 0\}}
{\bf 1}_{\{j_5=j_6\}}
\zeta_{j_1}^{(i_1)}
\zeta_{j_2}^{(i_2)}
\zeta_{j_3}^{(i_3)}
\zeta_{j_4}^{(i_4)}-
{\bf 1}_{\{i_1=i_2\ne 0\}}
{\bf 1}_{\{j_1=j_2\}}
\zeta_{j_3}^{(i_3)}
\zeta_{j_4}^{(i_4)}
\zeta_{j_5}^{(i_5)}
\zeta_{j_6}^{(i_6)}-
$$
$$
-
{\bf 1}_{\{i_1=i_3\ne 0\}}
{\bf 1}_{\{j_1=j_3\}}
\zeta_{j_2}^{(i_2)}
\zeta_{j_4}^{(i_4)}
\zeta_{j_5}^{(i_5)}
\zeta_{j_6}^{(i_6)}-
{\bf 1}_{\{i_1=i_4\ne 0\}}
{\bf 1}_{\{j_1=j_4\}}
\zeta_{j_2}^{(i_2)}
\zeta_{j_3}^{(i_3)}
\zeta_{j_5}^{(i_5)}
\zeta_{j_6}^{(i_6)}-
$$
$$
-
{\bf 1}_{\{i_1=i_5\ne 0\}}
{\bf 1}_{\{j_1=j_5\}}
\zeta_{j_2}^{(i_2)}
\zeta_{j_3}^{(i_3)}
\zeta_{j_4}^{(i_4)}
\zeta_{j_6}^{(i_6)}-
{\bf 1}_{\{i_2=i_3\ne 0\}}
{\bf 1}_{\{j_2=j_3\}}
\zeta_{j_1}^{(i_1)}
\zeta_{j_4}^{(i_4)}
\zeta_{j_5}^{(i_5)}
\zeta_{j_6}^{(i_6)}-
$$
$$
-
{\bf 1}_{\{i_2=i_4\ne 0\}}
{\bf 1}_{\{j_2=j_4\}}
\zeta_{j_1}^{(i_1)}
\zeta_{j_3}^{(i_3)}
\zeta_{j_5}^{(i_5)}
\zeta_{j_6}^{(i_6)}-
{\bf 1}_{\{i_2=i_5\ne 0\}}
{\bf 1}_{\{j_2=j_5\}}
\zeta_{j_1}^{(i_1)}
\zeta_{j_3}^{(i_3)}
\zeta_{j_4}^{(i_4)}
\zeta_{j_6}^{(i_6)}-
$$
$$
-
{\bf 1}_{\{i_3=i_4\ne 0\}}
{\bf 1}_{\{j_3=j_4\}}
\zeta_{j_1}^{(i_1)}
\zeta_{j_2}^{(i_2)}
\zeta_{j_5}^{(i_5)}
\zeta_{j_6}^{(i_6)}-
{\bf 1}_{\{i_3=i_5\ne 0\}}
{\bf 1}_{\{j_3=j_5\}}
\zeta_{j_1}^{(i_1)}
\zeta_{j_2}^{(i_2)}
\zeta_{j_4}^{(i_4)}
\zeta_{j_6}^{(i_6)}-
$$
$$
-
{\bf 1}_{\{i_4=i_5\ne 0\}}
{\bf 1}_{\{j_4=j_5\}}
\zeta_{j_1}^{(i_1)}
\zeta_{j_2}^{(i_2)}
\zeta_{j_3}^{(i_3)}
\zeta_{j_6}^{(i_6)}+
$$
$$
+
{\bf 1}_{\{i_1=i_2\ne 0\}}
{\bf 1}_{\{j_1=j_2\}}
{\bf 1}_{\{i_3=i_4\ne 0\}}
{\bf 1}_{\{j_3=j_4\}}
\zeta_{j_5}^{(i_5)}
\zeta_{j_6}^{(i_6)}+
{\bf 1}_{\{i_1=i_2\ne 0\}}
{\bf 1}_{\{j_1=j_2\}}
{\bf 1}_{\{i_3=i_5\ne 0\}}
{\bf 1}_{\{j_3=j_5\}}
\zeta_{j_4}^{(i_4)}
\zeta_{j_6}^{(i_6)}+
$$
$$
+
{\bf 1}_{\{i_1=i_2\ne 0\}}
{\bf 1}_{\{j_1=j_2\}}
{\bf 1}_{\{i_4=i_5\ne 0\}}
{\bf 1}_{\{j_4=j_5\}}
\zeta_{j_3}^{(i_3)}
\zeta_{j_6}^{(i_6)}
+
{\bf 1}_{\{i_1=i_3\ne 0\}}
{\bf 1}_{\{j_1=j_3\}}
{\bf 1}_{\{i_2=i_4\ne 0\}}
{\bf 1}_{\{j_2=j_4\}}
\zeta_{j_5}^{(i_5)}
\zeta_{j_6}^{(i_6)}+
$$
$$
+
{\bf 1}_{\{i_1=i_3\ne 0\}}
{\bf 1}_{\{j_1=j_3\}}
{\bf 1}_{\{i_2=i_5\ne 0\}}
{\bf 1}_{\{j_2=j_5\}}
\zeta_{j_4}^{(i_4)}
\zeta_{j_6}^{(i_6)}
+{\bf 1}_{\{i_1=i_3\ne 0\}}
{\bf 1}_{\{j_1=j_3\}}
{\bf 1}_{\{i_4=i_5\ne 0\}}
{\bf 1}_{\{j_4=j_5\}}
\zeta_{j_2}^{(i_2)}
\zeta_{j_6}^{(i_6)}+
$$
$$
+
{\bf 1}_{\{i_1=i_4\ne 0\}}
{\bf 1}_{\{j_1=j_4\}}
{\bf 1}_{\{i_2=i_3\ne 0\}}
{\bf 1}_{\{j_2=j_3\}}
\zeta_{j_5}^{(i_5)}
\zeta_{j_6}^{(i_6)}
+
{\bf 1}_{\{i_1=i_4\ne 0\}}
{\bf 1}_{\{j_1=j_4\}}
{\bf 1}_{\{i_2=i_5\ne 0\}}
{\bf 1}_{\{j_2=j_5\}}
\zeta_{j_3}^{(i_3)}
\zeta_{j_6}^{(i_6)}+
$$
$$
+
{\bf 1}_{\{i_1=i_4\ne 0\}}
{\bf 1}_{\{j_1=j_4\}}
{\bf 1}_{\{i_3=i_5\ne 0\}}
{\bf 1}_{\{j_3=j_5\}}
\zeta_{j_2}^{(i_2)}
\zeta_{j_6}^{(i_6)}
+
{\bf 1}_{\{i_1=i_5\ne 0\}}
{\bf 1}_{\{j_1=j_5\}}
{\bf 1}_{\{i_2=i_3\ne 0\}}
{\bf 1}_{\{j_2=j_3\}}
\zeta_{j_4}^{(i_4)}
\zeta_{j_6}^{(i_6)}+
$$
$$
+
{\bf 1}_{\{i_1=i_5\ne 0\}}
{\bf 1}_{\{j_1=j_5\}}
{\bf 1}_{\{i_2=i_4\ne 0\}}
{\bf 1}_{\{j_2=j_4\}}
\zeta_{j_3}^{(i_3)}
\zeta_{j_6}^{(i_6)}
+
{\bf 1}_{\{i_1=i_5\ne 0\}}
{\bf 1}_{\{j_1=j_5\}}
{\bf 1}_{\{i_3=i_4\ne 0\}}
{\bf 1}_{\{j_3=j_4\}}
\zeta_{j_2}^{(i_2)}
\zeta_{j_6}^{(i_6)}+
$$
$$
+
{\bf 1}_{\{i_2=i_3\ne 0\}}
{\bf 1}_{\{j_2=j_3\}}
{\bf 1}_{\{i_4=i_5\ne 0\}}
{\bf 1}_{\{j_4=j_5\}}
\zeta_{j_1}^{(i_1)}
\zeta_{j_6}^{(i_6)}
+
{\bf 1}_{\{i_2=i_4\ne 0\}}
{\bf 1}_{\{j_2=j_4\}}
{\bf 1}_{\{i_3=i_5\ne 0\}}
{\bf 1}_{\{j_3=j_5\}}
\zeta_{j_1}^{(i_1)}
\zeta_{j_6}^{(i_6)}+
$$
$$
+
{\bf 1}_{\{i_2=i_5\ne 0\}}
{\bf 1}_{\{j_2=j_5\}}
{\bf 1}_{\{i_3=i_4\ne 0\}}
{\bf 1}_{\{j_3=j_4\}}
\zeta_{j_1}^{(i_1)}
\zeta_{j_6}^{(i_6)}
+
{\bf 1}_{\{i_6=i_1\ne 0\}}
{\bf 1}_{\{j_6=j_1\}}
{\bf 1}_{\{i_3=i_4\ne 0\}}
{\bf 1}_{\{j_3=j_4\}}
\zeta_{j_2}^{(i_2)}
\zeta_{j_5}^{(i_5)}+
$$
$$
+
{\bf 1}_{\{i_6=i_1\ne 0\}}
{\bf 1}_{\{j_6=j_1\}}
{\bf 1}_{\{i_3=i_5\ne 0\}}
{\bf 1}_{\{j_3=j_5\}}
\zeta_{j_2}^{(i_2)}
\zeta_{j_4}^{(i_4)}
+
{\bf 1}_{\{i_6=i_1\ne 0\}}
{\bf 1}_{\{j_6=j_1\}}
{\bf 1}_{\{i_2=i_5\ne 0\}}
{\bf 1}_{\{j_2=j_5\}}
\zeta_{j_3}^{(i_3)}
\zeta_{j_4}^{(i_4)}+
$$
$$
+
{\bf 1}_{\{i_6=i_1\ne 0\}}
{\bf 1}_{\{j_6=j_1\}}
{\bf 1}_{\{i_2=i_4\ne 0\}}
{\bf 1}_{\{j_2=j_4\}}
\zeta_{j_3}^{(i_3)}
\zeta_{j_5}^{(i_5)}
+
{\bf 1}_{\{i_6=i_1\ne 0\}}
{\bf 1}_{\{j_6=j_1\}}
{\bf 1}_{\{i_4=i_5\ne 0\}}
{\bf 1}_{\{j_4=j_5\}}
\zeta_{j_2}^{(i_2)}
\zeta_{j_3}^{(i_3)}+
$$
$$
+
{\bf 1}_{\{i_6=i_1\ne 0\}}
{\bf 1}_{\{j_6=j_1\}}
{\bf 1}_{\{i_2=i_3\ne 0\}}
{\bf 1}_{\{j_2=j_3\}}
\zeta_{j_4}^{(i_4)}
\zeta_{j_5}^{(i_5)}
+
{\bf 1}_{\{i_6=i_2\ne 0\}}
{\bf 1}_{\{j_6=j_2\}}
{\bf 1}_{\{i_3=i_5\ne 0\}}
{\bf 1}_{\{j_3=j_5\}}
\zeta_{j_1}^{(i_1)}
\zeta_{j_4}^{(i_4)}+
$$
$$
+
{\bf 1}_{\{i_6=i_2\ne 0\}}
{\bf 1}_{\{j_6=j_2\}}
{\bf 1}_{\{i_4=i_5\ne 0\}}
{\bf 1}_{\{j_4=j_5\}}
\zeta_{j_1}^{(i_1)}
\zeta_{j_3}^{(i_3)}
+
{\bf 1}_{\{i_6=i_2\ne 0\}}
{\bf 1}_{\{j_6=j_2\}}
{\bf 1}_{\{i_3=i_4\ne 0\}}
{\bf 1}_{\{j_3=j_4\}}
\zeta_{j_1}^{(i_1)}
\zeta_{j_5}^{(i_5)}+
$$
$$
+
{\bf 1}_{\{i_6=i_2\ne 0\}}
{\bf 1}_{\{j_6=j_2\}}
{\bf 1}_{\{i_1=i_5\ne 0\}}
{\bf 1}_{\{j_1=j_5\}}
\zeta_{j_3}^{(i_3)}
\zeta_{j_4}^{(i_4)}
+
{\bf 1}_{\{i_6=i_2\ne 0\}}
{\bf 1}_{\{j_6=j_2\}}
{\bf 1}_{\{i_1=i_4\ne 0\}}
{\bf 1}_{\{j_1=j_4\}}
\zeta_{j_3}^{(i_3)}
\zeta_{j_5}^{(i_5)}+
$$
$$
+
{\bf 1}_{\{i_6=i_2\ne 0\}}
{\bf 1}_{\{j_6=j_2\}}
{\bf 1}_{\{i_1=i_3\ne 0\}}
{\bf 1}_{\{j_1=j_3\}}
\zeta_{j_4}^{(i_4)}
\zeta_{j_5}^{(i_5)}
+
{\bf 1}_{\{i_6=i_3\ne 0\}}
{\bf 1}_{\{j_6=j_3\}}
{\bf 1}_{\{i_2=i_5\ne 0\}}
{\bf 1}_{\{j_2=j_5\}}
\zeta_{j_1}^{(i_1)}
\zeta_{j_4}^{(i_4)}+
$$
$$
+
{\bf 1}_{\{i_6=i_3\ne 0\}}
{\bf 1}_{\{j_6=j_3\}}
{\bf 1}_{\{i_4=i_5\ne 0\}}
{\bf 1}_{\{j_4=j_5\}}
\zeta_{j_1}^{(i_1)}
\zeta_{j_2}^{(i_2)}
+
{\bf 1}_{\{i_6=i_3\ne 0\}}
{\bf 1}_{\{j_6=j_3\}}
{\bf 1}_{\{i_2=i_4\ne 0\}}
{\bf 1}_{\{j_2=j_4\}}
\zeta_{j_1}^{(i_1)}
\zeta_{j_5}^{(i_5)}+
$$
$$
+
{\bf 1}_{\{i_6=i_3\ne 0\}}
{\bf 1}_{\{j_6=j_3\}}
{\bf 1}_{\{i_1=i_5\ne 0\}}
{\bf 1}_{\{j_1=j_5\}}
\zeta_{j_2}^{(i_2)}
\zeta_{j_4}^{(i_4)}
+
{\bf 1}_{\{i_6=i_3\ne 0\}}
{\bf 1}_{\{j_6=j_3\}}
{\bf 1}_{\{i_1=i_4\ne 0\}}
{\bf 1}_{\{j_1=j_4\}}
\zeta_{j_2}^{(i_2)}
\zeta_{j_5}^{(i_5)}+
$$
$$
+
{\bf 1}_{\{i_6=i_3\ne 0\}}
{\bf 1}_{\{j_6=j_3\}}
{\bf 1}_{\{i_1=i_2\ne 0\}}
{\bf 1}_{\{j_1=j_2\}}
\zeta_{j_4}^{(i_4)}
\zeta_{j_5}^{(i_5)}
+
{\bf 1}_{\{i_6=i_4\ne 0\}}
{\bf 1}_{\{j_6=j_4\}}
{\bf 1}_{\{i_3=i_5\ne 0\}}
{\bf 1}_{\{j_3=j_5\}}
\zeta_{j_1}^{(i_1)}
\zeta_{j_2}^{(i_2)}+
$$
$$
+
{\bf 1}_{\{i_6=i_4\ne 0\}}
{\bf 1}_{\{j_6=j_4\}}
{\bf 1}_{\{i_2=i_5\ne 0\}}
{\bf 1}_{\{j_2=j_5\}}
\zeta_{j_1}^{(i_1)}
\zeta_{j_3}^{(i_3)}
+
{\bf 1}_{\{i_6=i_4\ne 0\}}
{\bf 1}_{\{j_6=j_4\}}
{\bf 1}_{\{i_2=i_3\ne 0\}}
{\bf 1}_{\{j_2=j_3\}}
\zeta_{j_1}^{(i_1)}
\zeta_{j_5}^{(i_5)}+
$$
$$
+
{\bf 1}_{\{i_6=i_4\ne 0\}}
{\bf 1}_{\{j_6=j_4\}}
{\bf 1}_{\{i_1=i_5\ne 0\}}
{\bf 1}_{\{j_1=j_5\}}
\zeta_{j_2}^{(i_2)}
\zeta_{j_3}^{(i_3)}
+
{\bf 1}_{\{i_6=i_4\ne 0\}}
{\bf 1}_{\{j_6=j_4\}}
{\bf 1}_{\{i_1=i_3\ne 0\}}
{\bf 1}_{\{j_1=j_3\}}
\zeta_{j_2}^{(i_2)}
\zeta_{j_5}^{(i_5)}+
$$
$$
+
{\bf 1}_{\{i_6=i_4\ne 0\}}
{\bf 1}_{\{j_6=j_4\}}
{\bf 1}_{\{i_1=i_2\ne 0\}}
{\bf 1}_{\{j_1=j_2\}}
\zeta_{j_3}^{(i_3)}
\zeta_{j_5}^{(i_5)}
+
{\bf 1}_{\{i_6=i_5\ne 0\}}
{\bf 1}_{\{j_6=j_5\}}
{\bf 1}_{\{i_3=i_4\ne 0\}}
{\bf 1}_{\{j_3=j_4\}}
\zeta_{j_1}^{(i_1)}
\zeta_{j_2}^{(i_2)}+
$$
$$
+
{\bf 1}_{\{i_6=i_5\ne 0\}}
{\bf 1}_{\{j_6=j_5\}}
{\bf 1}_{\{i_2=i_4\ne 0\}}
{\bf 1}_{\{j_2=j_4\}}
\zeta_{j_1}^{(i_1)}
\zeta_{j_3}^{(i_3)}
+
{\bf 1}_{\{i_6=i_5\ne 0\}}
{\bf 1}_{\{j_6=j_5\}}
{\bf 1}_{\{i_2=i_3\ne 0\}}
{\bf 1}_{\{j_2=j_3\}}
\zeta_{j_1}^{(i_1)}
\zeta_{j_4}^{(i_4)}+
$$
$$
+
{\bf 1}_{\{i_6=i_5\ne 0\}}
{\bf 1}_{\{j_6=j_5\}}
{\bf 1}_{\{i_1=i_4\ne 0\}}
{\bf 1}_{\{j_1=j_4\}}
\zeta_{j_2}^{(i_2)}
\zeta_{j_3}^{(i_3)}
+
{\bf 1}_{\{i_6=i_5\ne 0\}}
{\bf 1}_{\{j_6=j_5\}}
{\bf 1}_{\{i_1=i_3\ne 0\}}
{\bf 1}_{\{j_1=j_3\}}
\zeta_{j_2}^{(i_2)}
\zeta_{j_4}^{(i_4)}+
$$
$$
+
{\bf 1}_{\{i_6=i_5\ne 0\}}
{\bf 1}_{\{j_6=j_5\}}
{\bf 1}_{\{i_1=i_2\ne 0\}}
{\bf 1}_{\{j_1=j_2\}}
\zeta_{j_3}^{(i_3)}
\zeta_{j_4}^{(i_4)}-
$$
$$
-
{\bf 1}_{\{i_6=i_1\ne 0\}}
{\bf 1}_{\{j_6=j_1\}}
{\bf 1}_{\{i_2=i_5\ne 0\}}
{\bf 1}_{\{j_2=j_5\}}
{\bf 1}_{\{i_3=i_4\ne 0\}}
{\bf 1}_{\{j_3=j_4\}}-
$$
$$
-
{\bf 1}_{\{i_6=i_1\ne 0\}}
{\bf 1}_{\{j_6=j_1\}}
{\bf 1}_{\{i_2=i_4\ne 0\}}
{\bf 1}_{\{j_2=j_4\}}
{\bf 1}_{\{i_3=i_5\ne 0\}}
{\bf 1}_{\{j_3=j_5\}}-
$$
$$
-
{\bf 1}_{\{i_6=i_1\ne 0\}}
{\bf 1}_{\{j_6=j_1\}}
{\bf 1}_{\{i_2=i_3\ne 0\}}
{\bf 1}_{\{j_2=j_3\}}
{\bf 1}_{\{i_4=i_5\ne 0\}}
{\bf 1}_{\{j_4=j_5\}}-
$$
$$
-
{\bf 1}_{\{i_6=i_2\ne 0\}}
{\bf 1}_{\{j_6=j_2\}}
{\bf 1}_{\{i_1=i_5\ne 0\}}
{\bf 1}_{\{j_1=j_5\}}
{\bf 1}_{\{i_3=i_4\ne 0\}}
{\bf 1}_{\{j_3=j_4\}}-
$$
$$
-
{\bf 1}_{\{i_6=i_2\ne 0\}}
{\bf 1}_{\{j_6=j_2\}}
{\bf 1}_{\{i_1=i_4\ne 0\}}
{\bf 1}_{\{j_1=j_4\}}
{\bf 1}_{\{i_3=i_5\ne 0\}}
{\bf 1}_{\{j_3=j_5\}}-
$$
$$
-
{\bf 1}_{\{i_6=i_2\ne 0\}}
{\bf 1}_{\{j_6=j_2\}}
{\bf 1}_{\{i_1=i_3\ne 0\}}
{\bf 1}_{\{j_1=j_3\}}
{\bf 1}_{\{i_4=i_5\ne 0\}}
{\bf 1}_{\{j_4=j_5\}}-
$$
$$
-
{\bf 1}_{\{i_6=i_3\ne 0\}}
{\bf 1}_{\{j_6=j_3\}}
{\bf 1}_{\{i_1=i_5\ne 0\}}
{\bf 1}_{\{j_1=j_5\}}
{\bf 1}_{\{i_2=i_4\ne 0\}}
{\bf 1}_{\{j_2=j_4\}}-
$$
$$
-
{\bf 1}_{\{i_6=i_3\ne 0\}}
{\bf 1}_{\{j_6=j_3\}}
{\bf 1}_{\{i_1=i_4\ne 0\}}
{\bf 1}_{\{j_1=j_4\}}
{\bf 1}_{\{i_2=i_5\ne 0\}}
{\bf 1}_{\{j_2=j_5\}}-
$$
$$
-
{\bf 1}_{\{i_3=i_6\ne 0\}}
{\bf 1}_{\{j_3=j_6\}}
{\bf 1}_{\{i_1=i_2\ne 0\}}
{\bf 1}_{\{j_1=j_2\}}
{\bf 1}_{\{i_4=i_5\ne 0\}}
{\bf 1}_{\{j_4=j_5\}}-
$$
$$
-
{\bf 1}_{\{i_6=i_4\ne 0\}}
{\bf 1}_{\{j_6=j_4\}}
{\bf 1}_{\{i_1=i_5\ne 0\}}
{\bf 1}_{\{j_1=j_5\}}
{\bf 1}_{\{i_2=i_3\ne 0\}}
{\bf 1}_{\{j_2=j_3\}}-
$$
$$
-
{\bf 1}_{\{i_6=i_4\ne 0\}}
{\bf 1}_{\{j_6=j_4\}}
{\bf 1}_{\{i_1=i_3\ne 0\}}
{\bf 1}_{\{j_1=j_3\}}
{\bf 1}_{\{i_2=i_5\ne 0\}}
{\bf 1}_{\{j_2=j_5\}}-
$$
$$
-
{\bf 1}_{\{i_6=i_4\ne 0\}}
{\bf 1}_{\{j_6=j_4\}}
{\bf 1}_{\{i_1=i_2\ne 0\}}
{\bf 1}_{\{j_1=j_2\}}
{\bf 1}_{\{i_3=i_5\ne 0\}}
{\bf 1}_{\{j_3=j_5\}}-
$$
$$
-
{\bf 1}_{\{i_6=i_5\ne 0\}}
{\bf 1}_{\{j_6=j_5\}}
{\bf 1}_{\{i_1=i_4\ne 0\}}
{\bf 1}_{\{j_1=j_4\}}
{\bf 1}_{\{i_2=i_3\ne 0\}}
{\bf 1}_{\{j_2=j_3\}}-
$$
$$
-
{\bf 1}_{\{i_6=i_5\ne 0\}}
{\bf 1}_{\{j_6=j_5\}}
{\bf 1}_{\{i_1=i_2\ne 0\}}
{\bf 1}_{\{j_1=j_2\}}
{\bf 1}_{\{i_3=i_4\ne 0\}}
{\bf 1}_{\{j_3=j_4\}}-
$$
\begin{equation}
\label{a6}
\Biggl.-
{\bf 1}_{\{i_6=i_5\ne 0\}}
{\bf 1}_{\{j_6=j_5\}}
{\bf 1}_{\{i_1=i_3\ne 0\}}
{\bf 1}_{\{j_1=j_3\}}
{\bf 1}_{\{i_2=i_4\ne 0\}}
{\bf 1}_{\{j_2=j_4\}}\Biggr),
\end{equation}

\vspace{6mm}
\noindent
where ${\bf 1}_A$ is the indicator of the set $A$.

Consider the generalization of (\ref{a1})--(\ref{a6}) 
for the case of an arbitrary $k$ $(k\in \mathbb{N})$ 
as well as for the case of an arbitrary complete orthonormal system 
of functions $\{\phi_j(x)\}_{j=0}^{\infty}$ in the space $L_2([t,T])$
and $\psi_1(\tau),\ldots,\psi_k(\tau) \in L_2([t, T])$.

In order to do this, let us
consider the unordered
set $\{1, 2, \ldots, k\}$ 
and separate it into two parts:
the first part consists of $r$ unordered 
pairs (sequence order of these pairs is also unimportant) and the 
second one consists of the 
remaining $k-2r$ numbers.
So, we have

\vspace{-1mm}
\begin{equation}
\label{leto5007}
(\{
\underbrace{\{g_1, g_2\}, \ldots, 
\{g_{2r-1}, g_{2r}\}}_{\small{\hbox{part 1}}}
\},
\{\underbrace{q_1, \ldots, q_{k-2r}}_{\small{\hbox{part 2}}}
\}),
\end{equation}

\vspace{3mm}
\noindent
where 
$\{g_1, g_2, \ldots, 
g_{2r-1}, g_{2r}, q_1, \ldots, q_{k-2r}\}=\{1, 2, \ldots, k\},$
braces   
mean an unordered 
set, and parentheses mean an ordered set.

We will say that (\ref{leto5007}) is a partition 
and consider the sum with respect to all possible
partitions

\begin{equation}
\label{leto5008}
\sum_{\stackrel{(\{\{g_1, g_2\}, \ldots, 
\{g_{2r-1}, g_{2r}\}\}, \{q_1, \ldots, q_{k-2r}\})}
{{}_{\{g_1, g_2, \ldots, 
g_{2r-1}, g_{2r}, q_1, \ldots, q_{k-2r}\}=\{1, 2, \ldots, k\}}}}
a_{g_1 g_2, \ldots, 
g_{2r-1} g_{2r}, q_1 \ldots q_{k-2r}}.
\end{equation}

\vspace{2mm}

Below there are several examples of sums in the form (\ref{leto5008})

$$
\sum_{\stackrel{(\{g_1, g_2\})}{{}_{\{g_1, g_2\}=\{1, 2\}}}}
a_{g_1 g_2}=a_{12},
$$

\vspace{3mm}
$$
\sum_{\stackrel{(\{\{g_1, g_2\}, \{g_3, g_4\}\})}
{{}_{\{g_1, g_2, g_3, g_4\}=\{1, 2, 3, 4\}}}}
a_{g_1 g_2 g_3 g_4}=a_{1234} + a_{1324} + a_{2314},
$$

\vspace{3mm}
$$
\sum_{\stackrel{(\{g_1, g_2\}, \{q_1, q_{2}\})}
{{}_{\{g_1, g_2, q_1, q_{2}\}=\{1, 2, 3, 4\}}}}
a_{g_1 g_2, q_1 q_{2}}=a_{12,34}+a_{13,24}+a_{14,23}
+a_{23,14}+a_{24,13}+a_{34,12},
$$

\vspace{3mm}
$$
\sum_{\stackrel{(\{g_1, g_2\}, \{q_1, q_{2}, q_3\})}
{{}_{\{g_1, g_2, q_1, q_{2}, q_3\}=\{1, 2, 3, 4, 5\}}}}
a_{g_1 g_2, q_1 q_{2}q_3}
=a_{12,345}+a_{13,245}+a_{14,235}
+a_{15,234}+a_{23,145}+a_{24,135}+
$$
$$
+a_{25,134}+a_{34,125}+a_{35,124}+a_{45,123},
$$

\vspace{3mm}
$$
\sum_{\stackrel{(\{\{g_1, g_2\}, \{g_3, g_{4}\}\}, \{q_1\})}
{{}_{\{g_1, g_2, g_3, g_{4}, q_1\}=\{1, 2, 3, 4, 5\}}}}
a_{g_1 g_2, g_3 g_{4},q_1}
=
a_{12,34,5}+a_{13,24,5}+a_{14,23,5}+
a_{12,35,4}+a_{13,25,4}+a_{15,23,4}+
$$
$$
+a_{12,54,3}+a_{15,24,3}+a_{14,25,3}+a_{15,34,2}+a_{13,54,2}+a_{14,53,2}+
a_{52,34,1}+a_{53,24,1}+a_{54,23,1}.
$$

\vspace{8mm}

Now we can formulate the following generalization of Theorem 1.

\vspace{2mm}

{\bf Theorem 2} \cite{37a} (Sect.~1.11), \cite{38} (Sect.~15). {\it Suppose that
$\{\phi_j(x)\}_{j=0}^{\infty}$ is an arbitrary complete orthonormal system  
of functions in the space $L_2([t,T])$ and 
$\psi_1(\tau),\ldots,\psi_k(\tau)\in L_2([t, T]).$ 
Then the following expansion

\vspace{1mm}

$$
J[\psi^{(k)}]_{T,t}^{(i_1\ldots i_k)}=
\hbox{\vtop{\offinterlineskip\halign{
\hfil#\hfil\cr
{\rm l.i.m.}\cr
$\stackrel{}{{}_{p_1,\ldots,p_k\to \infty}}$\cr
}} }
\sum\limits_{j_1=0}^{p_1}\ldots
\sum\limits_{j_k=0}^{p_k}
C_{j_k\ldots j_1}\Biggl(
\prod_{l=1}^k\zeta_{j_l}^{(i_l)}+\sum\limits_{r=1}^{[k/2]}
(-1)^r \times
\Biggr.
$$

\vspace{2mm}
\begin{equation}
\label{leto6000}
\times
\sum_{\stackrel{(\{\{g_1, g_2\}, \ldots, 
\{g_{2r-1}, g_{2r}\}\}, \{q_1, \ldots, q_{k-2r}\})}
{{}_{\{g_1, g_2, \ldots, 
g_{2r-1}, g_{2r}, q_1, \ldots, q_{k-2r}\}=\{1, 2, \ldots, k\}}}}
\prod\limits_{s=1}^r
{\bf 1}_{\{i_{g_{{}_{2s-1}}}=~i_{g_{{}_{2s}}}\ne 0\}}
\Biggl.{\bf 1}_{\{j_{g_{{}_{2s-1}}}=~j_{g_{{}_{2s}}}\}}
\prod_{l=1}^{k-2r}\zeta_{j_{q_l}}^{(i_{q_l})}\Biggr)
\end{equation}

\vspace{6mm}
\noindent
con\-verg\-ing in the mean-square sense is valid, where $[x]$ is an integer part of 
a real number $x;$ another notations are the same as in Theorem~{\rm 1}.}

\vspace{4mm}

In particular, from (\ref{leto6000}) for $k=5$ we obtain

\vspace{3mm}

$$
J[\psi^{(5)}]_{T,t}^{(i_1\ldots i_5)}=
\sum_{j_1,\ldots,j_5=0}^{\infty}
C_{j_5\ldots j_1}\Biggl(
\prod_{l=1}^5\zeta_{j_l}^{(i_l)}-\Biggr.
\sum\limits_{\stackrel{(\{g_1, g_2\}, \{q_1, q_{2}, q_3\})}
{{}_{\{g_1, g_2, q_{1}, q_{2}, q_3\}=\{1, 2, 3, 4, 5\}}}}
{\bf 1}_{\{i_{g_{{}_{1}}}=~i_{g_{{}_{2}}}\ne 0\}}
{\bf 1}_{\{j_{g_{{}_{1}}}=~j_{g_{{}_{2}}}\}}
\prod_{l=1}^{3}\zeta_{j_{q_l}}^{(i_{q_l})}+
$$

\vspace{2mm}
$$
+
\sum_{\stackrel{(\{\{g_1, g_2\}, 
\{g_{3}, g_{4}\}\}, \{q_1\})}
{{}_{\{g_1, g_2, g_{3}, g_{4}, q_1\}=\{1, 2, 3, 4, 5\}}}}
{\bf 1}_{\{i_{g_{{}_{1}}}=~i_{g_{{}_{2}}}\ne 0\}}
{\bf 1}_{\{j_{g_{{}_{1}}}=~j_{g_{{}_{2}}}\}}
\Biggl.{\bf 1}_{\{i_{g_{{}_{3}}}=~i_{g_{{}_{4}}}\ne 0\}}
{\bf 1}_{\{j_{g_{{}_{3}}}=~j_{g_{{}_{4}}}\}}
\zeta_{j_{q_1}}^{(i_{q_1})}\Biggr).
$$

\vspace{6mm}

The last equality obviously agrees with
(\ref{a5}).
Note that the correctness of formulas (\ref{a1})--(\ref{a6}) 
can be 
verified 
by the fact that if 
$i_1=\ldots=i_6=i=1,\ldots,m$
and $\psi_1(s),\ldots,\psi_6(s)\equiv \psi(s)$,
then we can derive from (\ref{a1})--(\ref{a6}) the well known
equalities

$$
J[\psi^{(1)}]_{T,t}^{(i)}
=\frac{1}{1!}\delta_{T,t}^{(i)},
$$

\vspace{1mm}
$$
J[\psi^{(2)}]_{T,t}^{(ii)}
=\frac{1}{2!}\left(\left(\delta^{(i)}_{T,t}\right)^2-\Delta_{T,t}\right),\
$$

\vspace{1mm}
$$
J[\psi^{(3)}]_{T,t}^{(iii)}
=\frac{1}{3!}\left(\left(\delta_{T,t}^{(i)}\right)^3-
3\delta_{T,t}^{(i)}\Delta_{T,t}\right),
$$

\vspace{1mm}
$$
J[\psi^{(4)}]_{T,t}^{(iiii)}
=\frac{1}{4!}\left(\left(\delta_{T,t}^{(i)}\right)^4-
6\left(\delta_{T,t}^{(i)}\right)^2\Delta_{T,t}
+3\Delta^2_{T,t}\right),\
$$

\vspace{1mm}
$$
J[\psi^{(5)}]_{T,t}^{(iiiii)}
=\frac{1}{5!}\left(\left(\delta_{T,t}^{(i)}\right)^5-
10\left(\delta_{T,t}^{(i)}\right)^3\Delta_{T,t}
+15\delta_{T,t}^{(i)}\Delta^2_{T,t}\right),
$$

\vspace{1mm}
$$
J[\psi^{(6)}]_{T,t}^{(iiiiii)}
=\frac{1}{6!}\left(\left(\delta_{T,t}^{(i)}\right)^6-
15\left(\delta_{T,t}^{(i)}\right)^4\Delta_{T,t}
+45\left(\delta_{T,t}^{(i)}\right)^2\Delta^2_{T,t}-15\Delta_{T,t}^3\right)
$$

\vspace{2mm}
\noindent
w. p. 1 \cite{30a}-\cite{37}, where

$$
\delta_{T,t}^{(i)}=\int\limits_t^T\psi(s)d{\bf w}_s^{(i)},\ \ \
\Delta_{T,t}=\int\limits_t^T\psi^2(s)ds.
$$

\vspace{2mm}

The above equalities
can be independently  
obtained using the It\^{o} formula and Hermite polynomials \cite{16}.

\vspace{5mm} 

\section{Calculation of the Mean-Square Approximation Error
of Ite\-rated It\^{o} Stochastic Integrals in Theorems~1, 2}

\vspace{5mm}

Assume that $J[\psi^{(k)}]_{T,t}^{(i_1\ldots i_k)p_1 \ldots p_k}$ 
is an approximation 
of (\ref{sodom20}), which is
the expression on the right-hand side of {\rm (\ref{leto6000})} before passing to the limit 
$\hbox{\vtop{\offinterlineskip\halign{
\hfil#\hfil\cr
{\rm l.i.m.}\cr
$\stackrel{}{{}_{p_1,\ldots,p_k\to \infty}}$\cr
}} }$. Let us denote

$$
E^{(i_1\ldots i_k)p_1,\ldots,p_k}={\sf M}\biggl\{\biggl(
J[\psi^{(k)}]_{T,t}^{(i_1\ldots i_k)}-
J[\psi^{(k)}]_{T,t}^{(i_1\ldots i_k)p_1,\ldots,p_k}\biggr)^2\biggr\},
$$

\vspace{2mm}
$$
\left.E^{(i_1\ldots i_k)p}=E_k^{(i_1\ldots i_k)p_1,\ldots,p_k}\right|_{p_1
=\ldots=p_k=p},
$$

\vspace{2mm}
\begin{equation}
\label{g123}
I_k=\Vert K\Vert^2_{L_2([t, T]^k)}=\int\limits_{[t,T]^k}
K^2(t_1,\ldots,t_k)dt_1\ldots dt_k.
\end{equation}  

\vspace{2mm}

In \cite{34}-\cite{37aaxx}, \cite{46}-\cite{48} it was shown that

\begin{equation}
\label{star00011}
E_k^{(i_1\ldots i_k)p_1,\ldots,p_k}\le k!\left(I_k-\sum_{j_1=0}^{p_1}\ldots
\sum_{j_k=0}^{p_k}C^2_{j_k\ldots j_1}\right),
\end{equation}

\vspace{4mm}
\noindent
where $i_1,\ldots,i_k=1,\ldots,m$ for $0<T-t<\infty$ and
$i_1,\ldots,i_k=0, 1,\ldots,m$ for $0<T-t<1.$ Note that the
estimate (\ref{star00011}) is valid under the conditions of Theorem~2.

\vspace{2mm}

The exact calcutation of $E^{(i_1\ldots i_k)p}$ is presented in 
the following theorem.

\vspace{2mm}

{\bf Theorem 3} \cite{37a} (Sect.~1.12).
{\it {\it Suppose that $\{\phi_j(x)\}_{j=0}^{\infty}$ 
is an arbitrary complete orthonormal system  
of functions in the space $L_2([t,T])$ and
$\psi_1(\tau),\ldots,\psi_k(\tau)\in L_2([t, T]),$  $i_1,\ldots, i_k=1,\ldots,m$.
Then

\vspace{-1mm}
$$
E^{(i_1\ldots i_k)p}=I_k-
$$

\begin{equation}
\label{tttr11}
- \sum_{j_1,\ldots, j_k=0}^{p}
C_{j_k\ldots j_1}
{\sf M}\left\{J[\psi^{(k)}]_{T,t}^{(i_1\ldots i_k)}
\sum\limits_{(j_1,\ldots,j_k)}
\int\limits_t^T \phi_{j_k}(t_k)
\ldots
\int\limits_t^{t_{2}}\phi_{j_{1}}(t_{1})
d{\bf w}_{t_1}^{(i_1)}\ldots
d{\bf w}_{t_k}^{(i_k)}\right\},
\end{equation}

\vspace{5mm}
\noindent
where 
$J[\psi^{(k)}]_{T,t}^{(i_1\ldots i_k)p}$ is the expression on the right-hand side of
{\rm (\ref{leto6000})}
before passing to the limit 
$\hbox{\vtop{\offinterlineskip\halign{
\hfil#\hfil\cr
{\rm l.i.m.}\cr
$\stackrel{}{{}_{p_1,\ldots,p_k\to \infty}}$\cr
}} }$ for
$p_1=\ldots=p_k=p;$\
the expression 

$$
\sum\limits_{(j_1,\ldots,j_k)}
$$ 

\vspace{3mm}
\noindent
means the sum with respect to all
possible permutations
$(j_1,\ldots,j_k).$ At the same time if 
$j_r$ swapped with $j_q$ in the permutation $(j_1,\ldots,j_k)$,
then $i_r$ swapped with $i_q$ in the permutation 
$(i_1,\ldots,i_k);$
another notations are the same as in Theorems~{\rm 1, 2.}}
}

\vspace{2mm}

Note that 

\vspace{-1mm}
$$
{\sf M}\left\{J[\psi^{(k)}]_{T,t}^{(i_1\ldots i_k)}
\int\limits_t^T \phi_{j_k}(t_k)
\ldots
\int\limits_t^{t_{2}}\phi_{j_{1}}(t_{1})
d{\bf w}_{t_1}^{(i_1)}\ldots
d{\bf w}_{t_k}^{(i_k)}\right\}=C_{j_k\ldots j_1}
$$

\vspace{4mm}
\noindent
for $i_1\ldots i_k=1,\ldots,m.$

Then from Theorem 3 for $i_1,\ldots,i_k=1,\ldots,m$
we obtain \cite{35}, \cite{37pred}-\cite{37aaxx}

\vspace{1mm}

\begin{equation}
\label{12345}
E^{(i_1\ldots i_k)p}= I_k- \sum_{j_1,\ldots,j_k=0}^{p}
C_{j_k\ldots j_1}^2\ \ \ (\hbox{pairwise\ 
different}\ i_1,\ldots,i_k),
\end{equation}

\vspace{3mm}
$$
E^{(i_1 i_2)p}
=I_2
-\sum_{j_1,j_2=0}^p
C_{j_2j_1}^2-
\sum_{j_1,j_2=0}^p
C_{j_2j_1}C_{j_1j_2}\ \ \ (i_1=i_2),
$$

\vspace{3mm}
$$
E^{(i_1 i_2 i_3)p}=I_3
-\sum_{j_3,j_2,j_1=0}^p C_{j_3j_2j_1}^2-
\sum_{j_3,j_2,j_1=0}^p C_{j_3j_1j_2}C_{j_3j_2j_1}\ \ \ (i_1=i_2\ne i_3),
$$

\vspace{4mm}
$$
E^{(i_1i_2i_3i_4)p}
= I_4 -
\sum_{j_1,j_2,j_3,j_4=0}^{p}
C_{j_4j_3j_2j_1}\Biggl(\sum\limits_{(j_3,j_4)}
\Biggl(\sum\limits_{(j_1,j_2)}
C_{j_4j_3j_2j_1}\Biggr)\Biggr)\ \ \
(i_1=i_2\ne i_3=i_4),
$$

\vspace{5mm}
$$
E^{(i_1i_2i_3i_4i_5)p} = I_5 - \sum_{j_1,j_2,j_3,j_4,j_5=0}^{p}
C_{j_5j_4j_3j_2j_1}\Biggl(\sum\limits_{(j_3,j_4)}\Biggl(
\sum\limits_{(j_1,j_2,j_5)}
C_{j_5j_4j_3j_2j_1}\Biggr)\Biggr)
$$

\vspace{1mm}
$$
(i_1=i_2=i_5\ne i_3=i_4).
$$

\vspace{5mm}

\section{Some Examples of the Mean-Square Approximations of Ite\-rated 
It\^{o} Stochastic Integrals Using Legendre Polynomials}

\vspace{5mm}

Denote
$$
I_{(1)T,t}^{(i_1)}=\int\limits_t^Td{\bf w}_{t_1}^{(i_1)},\
$$

\vspace{1mm}
$$
I_{(10)T,t}^{(i_1 0)}=
\int\limits_t^T\int\limits_t^{t_2} d{\bf w}_{t_1}^{(i_1)}
dt_2,\ \ \
I_{(01)T,t}^{(0 i_2)}=\int\limits_t^T\int\limits_t^{t_2} dt_1
d{\bf w}_{t_2}^{(i_2)},
$$

\vspace{1mm}
$$
I_{(11)T,t}^{(i_1 i_2)}=\int\limits_t^T\int\limits_t^{t_2} d{\bf w}_{t_1}^{(i_1)}
d{\bf w}_{t_2}^{(i_2)},\ \ \
I_{(111)T,t}^{(i_1 i_2 i_3)}=\int\limits_t^T\int\limits_t^{t_3}\int\limits_t^{t_2} 
d{\bf w}_{t_1}^{(i_1)}
d{\bf w}_{t_2}^{(i_2)}
d{\bf w}_{t_3}^{(i_3)},
$$

\vspace{2mm}
$$
I_{(1111)T,t}^{(i_1 i_2 i_3 i_4)}=
\int\limits_t^T\int\limits_t^{t_4}\int\limits_t^{t_3}\int\limits_t^{t_2} 
d{\bf w}_{t_1}^{(i_1)}
d{\bf w}_{t_2}^{(i_2)}
d{\bf w}_{t_3}^{(i_3)}
d{\bf w}_{t_4}^{(i_4)},
$$

\vspace{3mm}
$$
I_{(11111)T,t}^{(i_1 i_2 i_3 i_4 i_5)}=
\int\limits_t^T\int\limits_t^{t_5}\int\limits_t^{t_4}
\int\limits_t^{t_3}\int\limits_t^{t_2} 
d{\bf w}_{t_1}^{(i_1)}
d{\bf w}_{t_2}^{(i_2)}
d{\bf w}_{t_3}^{(i_3)}
d{\bf w}_{t_4}^{(i_4)}
d{\bf w}_{t_5}^{(i_5)},
$$

\vspace{5mm}
\noindent
where $i_1, i_2, i_3, i_4, i_5=1,\ldots,m.$

The complete orthonormal system of Legendre polynomials in 
the space $L_2([t,T])$ looks as follows

\vspace{1mm}
\begin{equation}
\label{66}
\phi_j(x)=\sqrt{\frac{2j+1}{T-t}}P_j\biggl(\biggl(
x-\frac{T+t}{2}\biggr)\frac{2}{T-t}\biggr);\ j=0, 1, 2,\ldots,
\end{equation}

\vspace{4mm}
\noindent
where $P_j(x)$ is the Legendre polynomial. 

Using the system of 
functions (\ref{66}) and 
Theorems 1, 2 we obtain the following approximations of 
iterated 
It\^{o} stochastic integrals \cite{eee}-\cite{1003}

\vspace{2mm}
$$
I_{(1)T,t}^{(i_1)}=\sqrt{T-t}\zeta_0^{(i_1)},
$$

\vspace{1mm}
\begin{equation}
\label{opp1}
I_{(01)T,t}^{(0 i_1)}=\frac{(T-t)^{3/2}}{2}\biggl(\zeta_0^{(i_1)}+
\frac{1}{\sqrt{3}}\zeta_1^{(i_1)}\biggr),\
\end{equation}

\vspace{1mm}
\begin{equation}
\label{opp2}
I_{(10)T,t}^{(i_1 0)}=\frac{(T-t)^{3/2}}{2}\biggl(\zeta_0^{(i_1)}-
\frac{1}{\sqrt{3}}\zeta_1^{(i_1)}\biggr),
\end{equation}

\vspace{2mm}
$$
I_{(11)T,t}^{(i_1 i_2)q}=
\frac{T-t}{2}\left(\zeta_0^{(i_1)}\zeta_0^{(i_2)}+\sum_{i=1}^{q}
\frac{1}{\sqrt{4i^2-1}}\biggl(
\zeta_{i-1}^{(i_1)}\zeta_{i}^{(i_2)}-
\zeta_i^{(i_1)}\zeta_{i-1}^{(i_2)}\biggr)-{\bf 1}_{\{i_1=i_2\}}
\right),
$$

\vspace{5mm}
$$
I_{(111)T,t}^{(i_1 i_2 i_3)q_1}=
\sum_{j_1,j_2,j_3=0}^{q_1}
C_{j_3j_2j_1}
\Biggl(
\zeta_{j_1}^{(i_1)}\zeta_{j_2}^{(i_2)}\zeta_{j_3}^{(i_3)}
\biggr.-{\bf 1}_{\{i_1=i_2\}}
{\bf 1}_{\{j_1=j_2\}}
\zeta_{j_3}^{(i_3)}\Biggr.
-
$$

\begin{equation}
\label{kr1}
\Biggl.-{\bf 1}_{\{i_2=i_3\}}
{\bf 1}_{\{j_2=j_3\}}
\zeta_{j_1}^{(i_1)}-
{\bf 1}_{\{i_1=i_3\}}
{\bf 1}_{\{j_1=j_3\}}
\zeta_{j_2}^{(i_2)}\Biggr),
\end{equation}

\vspace{6mm}
$$
I_{(111)T,t}^{(i_1 i_1 i_1)}
=\frac{1}{6}(T-t)^{3/2}\left(
\biggl(\zeta_0^{(i_1)}\biggr)^3-3
\zeta_0^{(i_1)}\right),
$$

\vspace{6mm}

$$
I_{(1111)T,t}^{(i_1 i_2 i_3 i_4)q_2}
=
\sum_{j_1,j_2,j_3,j_4=0}^{q_2}
C_{j_4 j_3 j_2 j_1}\Biggl(
\zeta_{j_1}^{(i_1)}\zeta_{j_2}^{(i_2)}\zeta_{j_3}^{(i_3)}\zeta_{j_4}^{(i_4)}
\Biggr.-
$$
$$
-{\bf 1}_{\{i_1=i_2\}}
{\bf 1}_{\{j_1=j_2\}}
\zeta_{j_3}^{(i_3)}
\zeta_{j_4}^{(i_4)}
-
{\bf 1}_{\{i_1=i_3\}}
{\bf 1}_{\{j_1=j_3\}}
\zeta_{j_2}^{(i_2)}
\zeta_{j_4}^{(i_4)}
-
{\bf 1}_{\{i_1=i_4\}}
{\bf 1}_{\{j_1=j_4\}}
\zeta_{j_2}^{(i_2)}
\zeta_{j_3}^{(i_3)}
-
$$

\vspace{-2mm}
$$
-{\bf 1}_{\{i_2=i_3\}}
{\bf 1}_{\{j_2=j_3\}}
\zeta_{j_1}^{(i_1)}
\zeta_{j_4}^{(i_4)}
-
{\bf 1}_{\{i_2=i_4\}}
{\bf 1}_{\{j_2=j_4\}}
\zeta_{j_1}^{(i_1)}
\zeta_{j_3}^{(i_3)}
-
{\bf 1}_{\{i_3=i_4\}}
{\bf 1}_{\{j_3=j_4\}}
\zeta_{j_1}^{(i_1)}
\zeta_{j_2}^{(i_2)}+
$$

\vspace{-2mm}
$$
+
{\bf 1}_{\{i_1=i_2\}}
{\bf 1}_{\{j_1=j_2\}}
{\bf 1}_{\{i_3=i_4\}}
{\bf 1}_{\{j_3=j_4\}}
+
{\bf 1}_{\{i_1=i_3\}}
{\bf 1}_{\{j_1=j_3\}}
{\bf 1}_{\{i_2=i_4\}}
{\bf 1}_{\{j_2=j_4\}}+
$$
\begin{equation}
\label{rrr1}
+\Biggl.
{\bf 1}_{\{i_1=i_4\}}
{\bf 1}_{\{j_1=j_4\}}
{\bf 1}_{\{i_2=i_3\}}
{\bf 1}_{\{j_2=j_3\}}\Biggr),
\end{equation}

\vspace{7mm}

$$
I_{(1111)T,t}^{(i_1i_1i_1i_1)}=
\frac{1}{24}(T-t)^2
\left(\biggl(\zeta_0^{(i_1)}\biggr)^4-
6\biggl(\zeta_0^{(i_1)}\biggr)^2+3\right),
$$

\vspace{7mm}

$$
I_{(11111)T,t}^{(i_1 i_2 i_3 i_4 i_5)q_3}
=
\sum_{j_1,j_2,j_3,j_4,j_5=0}^{q_3}
C_{j_5j_4 j_3 j_2 j_1}\Biggl(
\zeta_{j_1}^{(i_1)}\zeta_{j_2}^{(i_2)}\zeta_{j_3}^{(i_3)}\zeta_{j_4}^{(i_4)}
\zeta_{j_5}^{(i_5)}
-\Biggr.
{\bf 1}_{\{i_1=i_2\}}
{\bf 1}_{\{j_1=j_2\}}
\zeta_{j_3}^{(i_3)}
\zeta_{j_4}^{(i_4)}
\zeta_{j_5}^{(i_5)}-
$$
$$
-
{\bf 1}_{\{i_1=i_3\}}
{\bf 1}_{\{j_1=j_3\}}
\zeta_{j_2}^{(i_2)}
\zeta_{j_4}^{(i_4)}
\zeta_{j_5}^{(i_5)}-
{\bf 1}_{\{i_1=i_4\}}
{\bf 1}_{\{j_1=j_4\}}
\zeta_{j_2}^{(i_2)}
\zeta_{j_3}^{(i_3)}
\zeta_{j_5}^{(i_5)}-
{\bf 1}_{\{i_1=i_5\}}
{\bf 1}_{\{j_1=j_5\}}
\zeta_{j_2}^{(i_2)}
\zeta_{j_3}^{(i_3)}
\zeta_{j_4}^{(i_4)}-
$$
$$
-
{\bf 1}_{\{i_2=i_3\}}
{\bf 1}_{\{j_2=j_3\}}
\zeta_{j_1}^{(i_1)}
\zeta_{j_4}^{(i_4)}
\zeta_{j_5}^{(i_5)}-
{\bf 1}_{\{i_2=i_4\}}
{\bf 1}_{\{j_2=j_4\}}
\zeta_{j_1}^{(i_1)}
\zeta_{j_3}^{(i_3)}
\zeta_{j_5}^{(i_5)}-
{\bf 1}_{\{i_2=i_5\}}
{\bf 1}_{\{j_2=j_5\}}
\zeta_{j_1}^{(i_1)}
\zeta_{j_3}^{(i_3)}
\zeta_{j_4}^{(i_4)}-
$$
$$
-{\bf 1}_{\{i_3=i_4\}}
{\bf 1}_{\{j_3=j_4\}}
\zeta_{j_1}^{(i_1)}
\zeta_{j_2}^{(i_2)}
\zeta_{j_5}^{(i_5)}-
{\bf 1}_{\{i_3=i_5\}}
{\bf 1}_{\{j_3=j_5\}}
\zeta_{j_1}^{(i_1)}
\zeta_{j_2}^{(i_2)}
\zeta_{j_4}^{(i_4)}
-{\bf 1}_{\{i_4=i_5\}}
{\bf 1}_{\{j_4=j_5\}}
\zeta_{j_1}^{(i_1)}
\zeta_{j_2}^{(i_2)}
\zeta_{j_3}^{(i_3)}+
$$
$$
+
{\bf 1}_{\{i_1=i_2\}}
{\bf 1}_{\{j_1=j_2\}}
{\bf 1}_{\{i_3=i_4\}}
{\bf 1}_{\{j_3=j_4\}}\zeta_{j_5}^{(i_5)}+
{\bf 1}_{\{i_1=i_2\}}
{\bf 1}_{\{j_1=j_2\}}
{\bf 1}_{\{i_3=i_5\}}
{\bf 1}_{\{j_3=j_5\}}\zeta_{j_4}^{(i_4)}+
$$
$$
+
{\bf 1}_{\{i_1=i_2\}}
{\bf 1}_{\{j_1=j_2\}}
{\bf 1}_{\{i_4=i_5\}}
{\bf 1}_{\{j_4=j_5\}}\zeta_{j_3}^{(i_3)}+
{\bf 1}_{\{i_1=i_3\}}
{\bf 1}_{\{j_1=j_3\}}
{\bf 1}_{\{i_2=i_4\}}
{\bf 1}_{\{j_2=j_4\}}\zeta_{j_5}^{(i_5)}+
$$
$$
+
{\bf 1}_{\{i_1=i_3\}}
{\bf 1}_{\{j_1=j_3\}}
{\bf 1}_{\{i_2=i_5\}}
{\bf 1}_{\{j_2=j_5\}}\zeta_{j_4}^{(i_4)}+
{\bf 1}_{\{i_1=i_3\}}
{\bf 1}_{\{j_1=j_3\}}
{\bf 1}_{\{i_4=i_5\}}
{\bf 1}_{\{j_4=j_5\}}\zeta_{j_2}^{(i_2)}+
$$
$$
+
{\bf 1}_{\{i_1=i_4\}}
{\bf 1}_{\{j_1=j_4\}}
{\bf 1}_{\{i_2=i_3\}}
{\bf 1}_{\{j_2=j_3\}}\zeta_{j_5}^{(i_5)}+
{\bf 1}_{\{i_1=i_4\}}
{\bf 1}_{\{j_1=j_4\}}
{\bf 1}_{\{i_2=i_5\}}
{\bf 1}_{\{j_2=j_5\}}\zeta_{j_3}^{(i_3)}+
$$
$$
+
{\bf 1}_{\{i_1=i_4\}}
{\bf 1}_{\{j_1=j_4\}}
{\bf 1}_{\{i_3=i_5\}}
{\bf 1}_{\{j_3=j_5\}}\zeta_{j_2}^{(i_2)}+
{\bf 1}_{\{i_1=i_5\}}
{\bf 1}_{\{j_1=j_5\}}
{\bf 1}_{\{i_2=i_3\}}
{\bf 1}_{\{j_2=j_3\}}\zeta_{j_4}^{(i_4)}+
$$
$$
+
{\bf 1}_{\{i_1=i_5\}}
{\bf 1}_{\{j_1=j_5\}}
{\bf 1}_{\{i_2=i_4\}}
{\bf 1}_{\{j_2=j_4\}}\zeta_{j_3}^{(i_3)}+
{\bf 1}_{\{i_1=i_5\}}
{\bf 1}_{\{j_1=j_5\}}
{\bf 1}_{\{i_3=i_4\}}
{\bf 1}_{\{j_3=j_4\}}\zeta_{j_2}^{(i_2)}+
$$
$$
+
{\bf 1}_{\{i_2=i_3\}}
{\bf 1}_{\{j_2=j_3\}}
{\bf 1}_{\{i_4=i_5\}}
{\bf 1}_{\{j_4=j_5\}}\zeta_{j_1}^{(i_1)}+
{\bf 1}_{\{i_2=i_4\}}
{\bf 1}_{\{j_2=j_4\}}
{\bf 1}_{\{i_3=i_5\}}
{\bf 1}_{\{j_3=j_5\}}\zeta_{j_1}^{(i_1)}+
$$
\begin{equation}
\label{rrr2}
+\Biggl.
{\bf 1}_{\{i_2=i_5\}}
{\bf 1}_{\{j_2=j_5\}}
{\bf 1}_{\{i_3=i_4\}}
{\bf 1}_{\{j_3=j_4\}}\zeta_{j_1}^{(i_1)}\Biggr),
\end{equation}

\vspace{5mm}

$$
I_{(11111)T,t}^{(i_1i_1i_1i_1i_1)}=
\frac{1}{120}(T-t)^{5/2}
\left(\biggl(\zeta_0^{(i_1)}\biggr)^5-
10\biggl(\zeta_0^{(i_1)}\biggr)^3+15\zeta_0^{(i_1)}\right),
$$

\vspace{7mm}

$$
C_{j_3j_2j_1}=\frac{\sqrt{(2j_1+1)(2j_2+1)(2j_3+1)}(T-t)^{3/2}}{8}\bar
C_{j_3j_2j_1},
$$

\vspace{3mm}
$$
C_{j_4j_3j_2j_1}
=\frac{\sqrt{(2j_1+1)(2j_2+1)(2j_3+1)(2j_4+1)}(T-t)^{2}}{16}\bar
C_{j_4j_3j_2j_1},
$$

\vspace{3mm}
$$
C_{j_5j_4 j_3 j_2 j_1}=
\frac{\sqrt{(2j_1+1)(2j_2+1)(2j_3+1)(2j_4+1)(2j_5+1)}(T-t)^{5/2}}{32}\bar
C_{j_5j_4 j_3 j_2 j_1},
$$

\vspace{4mm}
$$
\bar C_{j_3j_2j_1}=\int\limits_{-1}^{1}P_{j_3}(z)
\int\limits_{-1}^{z}P_{j_2}(y)
\int\limits_{-1}^{y}
P_{j_1}(x)dx dy dz,
$$

\vspace{3mm}
$$
\bar C_{j_4j_3j_2j_1}=\int\limits_{-1}^{1}P_{j_4}(u)
\int\limits_{-1}^{u}P_{j_3}(z)
\int\limits_{-1}^{z}P_{j_2}(y)
\int\limits_{-1}^{y}
P_{j_1}(x)dx dy dz du,
$$

\vspace{3mm}
$$
\bar C_{j_5j_4 j_3 j_2 j_1}=
\int\limits_{-1}^{1}P_{j_5}(v)
\int\limits_{-1}^{v}P_{j_4}(u)
\int\limits_{-1}^{u}P_{j_3}(z)
\int\limits_{-1}^{z}P_{j_2}(y)
\int\limits_{-1}^{y}
P_{j_1}(x)dx dy dz du dv,
$$

\vspace{7mm}
\noindent
random variables $\zeta_{j}^{(i)}$ are defined by 
(\ref{rr23}), and

\vspace{3mm}
$$
I_{(11)T,t}^{(i_1 i_2)}=\hbox{\vtop{\offinterlineskip\halign{
\hfil#\hfil\cr
{\rm l.i.m.}\cr
$\stackrel{}{{}_{q\to \infty}}$\cr
}} }I_{(11)T,t}^{(i_1 i_2)q},\ \ \
I_{(111)T,t}^{(i_1 i_2 i_3)}=\hbox{\vtop{\offinterlineskip\halign{
\hfil#\hfil\cr
{\rm l.i.m.}\cr
$\stackrel{}{{}_{q_1\to \infty}}$\cr
}} }I_{(111)T,t}^{(i_1 i_2 i_3)q_1},
$$

\vspace{3mm}
$$
I_{(1111)T,t}^{(i_1 i_2 i_3 i_4)}=\hbox{\vtop{\offinterlineskip\halign{
\hfil#\hfil\cr
{\rm l.i.m.}\cr
$\stackrel{}{{}_{q_2\to \infty}}$\cr
}} }I_{(1111)T,t}^{(i_1 i_2 i_3 i_4)q_2},\ \ \ 
I_{(11111)T,t}^{(i_1 i_2 i_3 i_4 i_5)}=\hbox{\vtop{\offinterlineskip\halign{
\hfil#\hfil\cr
{\rm l.i.m.}\cr
$\stackrel{}{{}_{q_3\to \infty}}$\cr
}} }I_{(11111)T,t}^{(i_1 i_2 i_3 i_4 i_5)q_3}.
$$

\vspace{6mm}

Note that $T-t\ll 1$ ($T-t$ is an 
integration step with respect to the temporal variable). Thus
$q_1\ll q$ (see Table 1 \cite{30}-\cite{34}, \cite{37pred}-\cite{37aaxx}). 
Moreover, the values $\bar C_{j_3j_2j_1},$ 
$\bar C_{j_4j_3j_2j_1},$
$\bar C_{j_5j_4j_3j_2j_1}$
do not depend on
$T-t.$ This feature is important because 
we can use a variable integration step $T-t$.
Coefficients 
$\bar C_{j_3j_2j_1},$ 
$\bar C_{j_4j_3j_2j_1},$
$\bar C_{j_5j_4j_3j_2j_1}$
are calculated once and before the start 
of the numerical scheme.
Some examples of the exact calculation of coefficients
$\bar C_{j_3j_2j_1},$ 
$\bar C_{j_4j_3j_2j_1},$
$\bar C_{j_5j_4j_3j_2j_1}$
via Python programming language can be found in Tables 2--4
(the database with 270,000 exactly
calculated Fourier--Legendre coefficients was described in \cite{1000}, \cite{1001}).

Denote

\vspace{-1mm}
$$
E^{(i_1i_2)q}={\sf M}\biggl\{\biggl(
I_{(11)T,t}^{(i_1 i_2)}-
I_{(11)T,t}^{(i_1 i_2)q}\biggr)^2\biggr\},
$$

\vspace{2mm}
$$
E^{(i_1i_2i_3)q_1}={\sf M}\biggl\{\biggl(
I_{(111)T,t}^{(i_1 i_2 i_3)}-
I_{(111)T,t}^{(i_1 i_2 i_3)q_1}\biggr)^2\biggr\},
$$

\vspace{2mm}
$$
E^{(i_1i_2i_3i_4)q_2}={\sf M}\biggl\{\biggl(
I_{(1111)T,t}^{(i_1 i_2 i_3 i_4)}-
I_{(1111)T,t}^{(i_1 i_2 i_3 i_4)q_2}\biggr)^2\biggr\},
$$ 

\vspace{2mm}
$$
E^{(i_1i_2i_3i_4i_5)q_3}={\sf M}\biggl\{\biggl(
I_{(11111)T,t}^{(i_1 i_2 i_3 i_4i_5)}-
I_{(11111)T,t}^{(i_1 i_2 i_3 i_4i_5)q_3}\biggr)^2\biggr\}.
$$ 

\vspace{5mm}

\noindent
\begin{figure}
\begin{center}
\centerline{Table 1.\ Minimal numbers\ $q,$\ $q_1$\  
such that\ $E^{(i_1i_2)q},$\ $E^{(i_1i_2i_3)q_1}\le (T-t)^4,$\ \ \ $q_1\ll q$.}
\vspace{4mm}
\tabcolsep=0.28em
\begin{tabular}{|c|c|c|c|c|c|c|}
\hline
$T-t$&$0.08222$&$0.05020$&$0.02310$&$0.01956$\\
\hline
$q$&19&51&235&328\\
\hline
$q_1$&1&2&5&6\\
\hline
\end{tabular}
\end{center}
\vspace{8mm}
\begin{center}
\centerline{Table 2.\ Coefficients $\bar C_{3jk}.$}
\vspace{4mm}
\begin{tabular}{|c|c|c|c|c|c|c|c|c|}
\hline
${}_j {}^k$&0&1&2&3&4&5&6\\
\hline
0&$0$&$\frac{2}{105}$&$0$&$-\frac{4}{315}$&$0$&$\frac{2}{693}$&0\\
\hline
1&$\frac{4}{105}$&0&$-\frac{2}{315}$&0&$-\frac{8}{3465}$&0&$\frac{10}{9009}$\\
\hline
2&$\frac{2}{35}$&$-\frac{2}{105}$&$0$&$\frac{4}{3465}$&
$0$&$-\frac{74}{45045}$&0\\
\hline
3&$\frac{2}{315}$&$0$&$-\frac{2}{3465}$&0&
$\frac{16}{45045}$&0&$-\frac{10}{9009}$\\
\hline
4&$-\frac{2}{63}$&$\frac{46}{3465}$&0&$-\frac{32}{45045}$&
0&$\frac{2}{9009}$&0\\
\hline
5&$-\frac{10}{693}$&0&$\frac{38}{9009}$&0&
$-\frac{4}{9009}$&0&$\frac{122}{765765}$\\
\hline
6&$0$&$-\frac{10}{3003}$&$0$&$\frac{20}{9009}$&$0$&$-\frac{226}{765765}$&$0$\\
\hline
\end{tabular}
\end{center}
\vspace{8mm}
\begin{center}
\centerline{Table 3.\ Coefficients $\bar C_{21kl}.$}
\vspace{4mm}
\begin{tabular}{|c|c|c|c|}
\hline
${}_k {}^l$&0&1&2\\
\hline
0&$\frac{2}{21}$&$-\frac{2}{45}$&$\frac{2}{315}$\\
\hline
1&$\frac{2}{315}$&$\frac{2}{315}$&$-\frac{2}{225}$\\
\hline
2&$-\frac{2}{105}$&$\frac{2}{225}$&$\frac{2}{1155}$\\
\hline
\end{tabular}
\end{center}
\vspace{8mm}
\begin{center}
\centerline{Table 4.\ Coefficients $\bar C_{101lr}.$}
\vspace{4mm}
\begin{tabular}{|c|c|c|}
\hline
${}_l {}^r$&0&1\\
\hline
0&$\frac{4}{315}$&$0$\\
\hline
1&$\frac{4}{315}$&$-\frac{8}{945}$\\
\hline
\end{tabular}
\end{center}
\vspace{5mm}
\end{figure}

Then for pairwise different $i_1,i_2,i_3,i_4,i_5=1,\ldots,m$
from Theorem 3 we obtain \cite{eee}-\cite{1003}

\vspace{1mm}
\begin{equation}
\label{909}
E^{(i_1i_2)q}=\frac{(T-t)^2}{2}\left(\frac{1}{2}-\sum_{i=1}^{q}
\frac{1}{4i^2-1}\right),
\end{equation}

\vspace{1mm}
\begin{equation}
\label{zzz}
E^{(i_1 i_2i_3)q_1}=
\frac{(T-t)^{3}}{6}-\sum_{j_1,j_2,j_3=0}^{q_1}
C_{j_3j_2j_1}^2,
\end{equation}

\vspace{1mm}
\begin{equation}
\label{zzz999}
E^{(i_1 i_2i_3i_4)q_2}=
\frac{(T-t)^{4}}{24}-\sum_{j_1,j_2,j_3,j_4=0}^{q_2}
C_{j_4j_3j_2j_1}^2,
\end{equation}

\vspace{1mm}
\begin{equation}
\label{zzz9999}
E^{(i_1 i_2i_3i_4i_5)q_3}=
\frac{(T-t)^{5}}{120}-\sum_{j_1,j_2,j_3,j_4,j_5=0}^{q_3}
C_{j_5j_4j_3j_2j_1}^2.
\end{equation}

\vspace{5mm}

On the basis of 
the presented 
approximations of 
iterated It\^{o} stochastic integrals we 
can see that increasing of multiplicities of these integrals 
leads to increasing 
of orders of smallness with respect to $T-t$
($T-t\ll 1$) in the mean-square sense 
for iterated It\^{o} stochastic integrals. This leads to a sharp decrease  
of member 
quantities
in the approximations of iterated It\^{o} stochastic
integrals, which are required for achieving the acceptable accuracy
of approximation ($q_1\ll q$). 

From (\ref{zzz})--(\ref{zzz9999}) we obtain
\cite{30}-\cite{34}, \cite{37pred}-\cite{37aaxx}

\vspace{1mm}
\begin{equation}
\label{ooo1}
\left.E^{(i_1i_2i_3)q_1}\right|_{q_1=6}
\approx
0.01956000(T-t)^3,
\end{equation}

\vspace{1mm}

\begin{equation}
\label{ooo2}
\left.E^{(i_1i_2i_3i_4)q_2}\right|_{q_2=2}
\approx 0.02360840(T-t)^4,
\end{equation}

\vspace{1mm}

\begin{equation}
\label{ooo3}
\left.E^{(i_1i_2i_3i_4i_5)q_3}\right|_{q_3=1}
\approx 0.00759105(T-t)^5.
\end{equation}

\vspace{4mm}

It is not difficult to see that the accuracy in (\ref{ooo2})
and (\ref{ooo3}) is significantly better
than in (\ref{ooo1}) 
$(T-t\ll 1)$ even for $q_2=2$ and $q_3=1.$
This means that in such situation in formulas 
(\ref{rrr1}), (\ref{rrr2})
the number of terms can be chosen significantly less 
than $3^4$ ($q_2=2$) and $2^5$ ($q_3=1$).
So, in practice, we can leave only few terms in these formulas.
For more details see \cite{1000}-\cite{1003}.

\vspace{5mm}

\section{Approximation of Iterated Stochastic Integrals of Multiplicity
$k$ with
Respect to the $Q$-Wiener Process}

\vspace{5mm}

Consider the iterated stochastic integral with respect to 
the $Q$-Wiener process in the form

\vspace{1mm}
$$
I[\Phi^{(k)}(Z), \psi^{(k)}]_{T,t}=
$$

\vspace{-1mm}
\begin{equation}
\label{ssss11}
=
\int\limits_{t}^{T}\Phi_k(Z) \left( \ldots \left(
\int\limits_{t}^{t_3}\Phi_2(Z) \left(
\int\limits_{t}^{t_2}\Phi_1(Z)
\psi_1(t_1)d{\bf W}_{t_1} \right)
\psi_2(t_2)d{\bf W}_{t_2} \right) \ldots  \right) \psi_k(t_k)d{\bf W}_{t_k},
\end{equation}

\vspace{5mm}
\noindent
where $Z: \Omega \rightarrow H$ is 
an ${\bf F}_t/{\mathcal{B}}(H)$-measurable mapping, 
$\Phi_k(v)(\ \ldots (\Phi_2(v)(\Phi_1(v))) \ldots\ )$ 
is a $k$-linear Hilbert--Schmidt operator
mapping from
$\underbrace{U_0\times \ldots \times U_0}_{\small{\hbox{$k$ times}}}$ to $H$
for all $v\in H$, and $\psi_1(\tau),\ldots,\psi_k(\tau)\in L_2([t, T]).$

\vspace{2mm} 

Let $I[\Phi^{(k)}(Z), \psi^{(k)}]_{T,t}^M$ be an approximation
of the stochastic integral (\ref{ssss11})

\vspace{2mm}
$$
I[\Phi^{(k)}(Z), \psi^{(k)}]_{T,t}^M=
$$

\vspace{2mm}
$$
=
\int\limits_{t}^{T}\Phi_k(Z) \left( \ldots \left(
\int\limits_{t}^{t_3}\Phi_2(Z) \left(
\int\limits_{t}^{t_2}\Phi_1(Z)
\psi_1(t_1)d{\bf W}_{t_1}^M \right)
\psi_2(t_2)d{\bf W}_{t_2}^M 
\right) \ldots  \right) \psi_k(t_k)d{\bf W}_{t_k}^M=
$$

\vspace{4mm}
$$
=\sum_{r_1,r_2,\ldots,r_k\in J_M}
\Phi_k(Z)\left(\ldots \left(\Phi_2(Z) \left(\Phi_1(Z)
e_{r_1} \right) 
e_{r_2} \right) \ldots \right) e_{r_k}\times
$$

\vspace{2mm}

\begin{equation}
\label{xx605}
\times
\sqrt{\lambda_{r_1}\lambda_{r_2}\ldots \lambda_{r_k}}\ 
J[\psi^{(k)}]_{T,t}^{(r_1 r_2\ldots r_k)},
\end{equation}

\vspace{6mm}
\noindent
where $0\le t<T\le \bar{T},$  and

\vspace{2mm}
$$
J[\psi^{(k)}]_{T,t}^{(r_1\ldots r_k)}=
\int\limits_t^T \psi_k(t_k)\ldots \int\limits_t^{t_3}
\psi_2(t_2)\int\limits_t^{t_{2}}\psi_1(t_1)
d{\bf w}_{t_1}^{(r_1)}d{\bf w}_{t_2}^{(r_2)}\ldots
d{\bf w}_{t_k}^{(r_k)}
$$

\vspace{4mm}
\noindent
is the iterated It\^{o} stochastic integral (\ref{sodom20}),
$r_1,r_2,\ldots,r_k\in J_M.$

Let $I[\Phi^{(k)}(Z), \psi^{(k)}]_{T,t}^{M,p_1\ldots,p_k}$ be an
approximation of the stochastic integral (\ref{xx605})

\vspace{5mm}
$$
I[\Phi^{(k)}(Z), \psi^{(k)}]_{T,t}^{M,p_1\ldots,p_k}
=
$$

\vspace{2mm}
$$
=\sum_{r_1,r_2,\ldots,r_k\in J_M}
\Phi_k(Z)\left(\ldots \left(\Phi_2(Z) \left(\Phi_1(Z)
e_{r_1} \right) 
e_{r_2} \right) \ldots \right) e_{r_k}\times
$$

\vspace{2mm}
\begin{equation}
\label{xx705}
\times
\sqrt{\lambda_{r_1}\lambda_{r_2}\ldots \lambda_{r_k}}\
J[\psi^{(k)}]_{T,t}^{(r_1 r_2\ldots r_k)p_1,\ldots,p_k},
\end{equation}

\vspace{6mm}
\noindent
where
$J[\psi^{(k)}]_{T,t}^{(r_1 r_2\ldots r_k)p_1,\ldots,p_k}$
is defined as a prelimit expression on the right-hand side of  (\ref{leto6000})

\vspace{2mm}
$$
J[\psi^{(k)}]_{T,t}^{(r_1\ldots r_k)p_1\ldots p_k}=
\sum\limits_{j_1=0}^{p_1}\ldots
\sum\limits_{j_k=0}^{p_k}
C_{j_k\ldots j_1}\Biggl(
\prod_{l=1}^k\zeta_{j_l}^{(r_l)}+\sum\limits_{m=1}^{[k/2]}
(-1)^m \times
\Biggr.
$$

\vspace{3mm}
\begin{equation}
\label{f1}
\times
\sum_{\stackrel{(\{\{g_1, g_2\}, \ldots, 
\{g_{2m-1}, g_{2m}\}\}, \{q_1, \ldots, q_{k-2m}\})}
{{}_{\{g_1, g_2, \ldots, 
g_{2m-1}, g_{2m}, q_1, \ldots, q_{k-2m}\}=\{1, 2, \ldots, k\}}}}
\prod\limits_{s=1}^m
{\bf 1}_{\{r_{g_{{}_{2s-1}}}=~r_{g_{{}_{2s}}}\ne 0\}}
\Biggl.{\bf 1}_{\{j_{g_{{}_{2s-1}}}=~j_{g_{{}_{2s}}}\}}
\prod_{l=1}^{k-2m}\zeta_{j_{q_l}}^{(r_{q_l})}\Biggr).
\end{equation}

\vspace{6mm}

Let $U,$ $H$ be separable
$\mathbb{R}$-Hilbert spaces,
$U_{0}=Q^{1/2}(U)$, and
$L(U, H)$ be the space of linear and bounded operators mapping
from $U$ to $H$. Let
$L(U, H)_{0}=\left\{T \vert_{U_{0}}:\
T\in L(U, H)\right\}$ (here $T \vert_{U_{0}}$
is the restriction of operator $T$ to
the space $U_0$). It is known \cite{7} that
$L(U, H)_{0}$
is a dense subset of the space of
Hilbert--Schmidt operators $L_{HS}(U_{0}, H)$.

\vspace{2mm}

{\bf Theorem 4}\ \cite{37a}-\cite{37aaxx}, \cite{48aa}, \cite{48aaa}, \cite{2000}. {\it 
Suppose that
$\{\phi_j(x)\}_{j=0}^{\infty}$ is an arbitrary complete orthonormal system  
of functions in the space $L_2([t,T])$ and
$\psi_1(\tau),\ldots,\psi_k(\tau)\in L_2([t, T]).$ 
Furthermore$,$ let the following conditions be satisfied$:$

{\rm 1}. $Q\in L(U)$ is a nonnegative and symmetric 
trace class operator {\rm (}$\lambda_i$ and $e_i$ $(i\in J)$ are
its eigenvalues and eigenfunctions {\rm (}which form
an orthonormal basis of $U${\rm )}
correspondingly{\rm )}, and ${\bf W}_{\tau},$ $\tau\in [0, \bar T]$
is an $U$-valued $Q$-Wiener process.

{\rm 2}. $Z: \Omega \rightarrow H$ is 
an ${\bf F}_t/{\mathcal{B}}(H)$-measurable mapping.

{\rm 3}. $\Phi_1\in L(U, H)_{0},$ 
$\Phi_2\in L(H,L(U,H)_0),$ and $\Phi_k(v)(\ \ldots (\Phi_2(v)(\Phi_1(v))) \ldots\ )$ 
is a $k$-linear Hilbert--Schmidt operator mapping from
$\underbrace{U_0\times \ldots \times U_0}_{\small{\hbox{$k$ times}}}$ to $H$
for all $v\in H$
such that

\vspace{1mm}

$$
\Biggl\Vert \Phi_k(Z)\left(\ldots \left(\Phi_2(Z) \left(\Phi_1(Z)
e_{r_1} \right) 
e_{r_2} \right) \ldots \right) e_{r_k}\Biggr\Vert_H^2
\le L_k<\infty
$$

\vspace{4mm}
\noindent
w.\ p.\ {\rm 1}\ for all $r_1,r_2,\ldots,r_k\in J_M$, $M\in\mathbb{N}$.
\vspace{2mm}
\noindent

Then
$$
{\sf M}\left\{
\Biggl\Vert
I[\Phi^{(k)}(Z), \psi^{(k)}]_{T,t}^M-
I[\Phi^{(k)}(Z), \psi^{(k)}]_{T,t}^{M,p_1\ldots p_k}
\Biggr\Vert_H^2\right\}\le 
$$

\vspace{-2mm}
\begin{equation}
\label{zzz1}
\le L_k (k!)^2
\left({\rm tr}\ Q\right)^k
\left(I_k-\sum_{j_1=0}^{p_1}\ldots
\sum_{j_k=0}^{p_k}C^2_{j_k\ldots j_1}\right),
\end{equation}

\vspace{5mm}
\noindent
where $I_k$ is defined by {\rm (\ref{g123})}, ${\rm tr}\ Q=\sum\limits_{i\in J}\lambda_i,$
and 

$$
C_{j_k\ldots j_1}=\int\limits_{[t,T]^k}
K(t_1,\ldots,t_k)\prod_{l=1}^{k}\phi_{j_l}(t_l)dt_1\ldots dt_k,
$$

\vspace{1mm}
$$
K(t_1,\ldots,t_k)=
\begin{cases}
\psi_1(t_1)\ldots \psi_k(t_k),\ t_1<\ldots<t_k\\
~\\
~\\
0,\ \hbox{\rm otherwise}
\end{cases}.
$$

}

\vspace{3mm}

{\bf Remark 2.}\ {\it It should be noted that the right-hand side 
of the inequality {\rm (\ref{zzz1})} is independent of $M$
and tends to zero if $p_1,\ldots,p_k\to\infty$ due to the
Parseval equality.}

\vspace{3mm}

{\bf Proof.} Using (\ref{star00011}), we obtain

$$
{\sf M}\left\{
\Biggl\Vert
I[\Phi^{(k)}(Z), \psi^{(k)}]_{T,t}^M-
I[\Phi^{(k)}(Z), \psi^{(k)}]_{T,t}^{M,p_1\ldots p_k}
\Biggr\Vert_H^2\right\}=
$$

\vspace{10mm}
$$
={\sf M}\left\{
\Biggl\Vert
\sum_{r_1,r_2,\ldots,r_k\in J_M}
\Phi_k(Z)\left(\ldots \left(\Phi_2(Z) \left(\Phi_1(Z)
e_{r_1} \right) 
e_{r_2} \right) \ldots \right) e_{r_k}
\sqrt{\lambda_{r_1}\lambda_{r_2}\ldots \lambda_{r_k}}\times\right.\Biggr.
$$
\begin{equation}
\label{zz1}
~~~~~~~~~~~~~~~~~~~~~~~~~~~~~~~~~~~~~~~~~~~~~~~~~~~\times
\left.\Biggl.\Biggl(J[\psi^{(k)}]_{T,t}^{(r_1 r_2\ldots r_k)}-
J[\psi^{(k)}]_{T,t}^{(r_1 r_2\ldots r_k)p_1,\ldots,p_k}\Biggr)
\Biggr\Vert_H^2\right\}=
\end{equation}

\vspace{10mm}
$$
=\left|
{\sf M}\Biggl\{\sum_{r_1,r_2,\ldots,r_k\in J_M}\
\sum_{(r_1^{,},r_2^{,},\ldots,r_k^{,}):\ 
\{r_1^{,},r_2^{,},\ldots,r_k^{,}\}=\{r_1,r_2,\ldots,r_k\}}
\Biggl\langle\ \Phi_k(Z)\left(\ldots \left(\Phi_2(Z) \left(\Phi_1(Z)
e_{r_1} \right) 
e_{r_2} \right) \ldots \right) e_{r_k}\ ,\Biggr.\Biggr.\right. 
$$

\vspace{5mm}
$$
~~~~~~~~~~~\Biggl.
\Phi_k(Z)\left(\ldots \left(\Phi_2(Z) \left(\Phi_1(Z)
e_{r_1^{,}} \right) 
e_{r_2^{,}} \right) \ldots \right) e_{r_k^{,}}\ \Biggr\rangle_H\
\sqrt{\lambda_{r_1}\lambda_{r_2}\ldots \lambda_{r_k}}
\sqrt{\lambda_{r_1^{,}}\lambda_{r_2^{,}}\ldots \lambda_{r_k^{,}}}\times
$$

\vspace{5mm}

$$
\times\
{\sf M}\Biggl\{\Biggl(J[\psi^{(k)}]_{T,t}^{(r_1 r_2\ldots r_k)}-
J[\psi^{(k)}]_{T,t}^{(r_1 r_2\ldots r_k)p_1,\ldots,p_k}\Biggr)\times\Biggr.
$$

\vspace{5mm}

\begin{equation}
\label{zz100}
\left.
\Biggl.
\Biggl.
~~~~~~~~~~~~~~~~~~~~~~~~~~~~~~~~~~~~~\times
\Biggl(J[\psi^{(k)}]_{T,t}^{(r_1^{,} r_2^{,}\ldots r_k^{,})}-
J[\psi^{(k)}]_{T,t}^{(r_1^{,} r_2^{,}\ldots r_k^{,})p_1,\ldots,p_k}\Biggr)
\biggl.\biggr|{\bf F}_t\Biggr\}
\Biggr\}\right|\le
\end{equation}

\vspace{13mm}
$$
\le
\hspace{-1mm}\sum_{r_1,r_2,\ldots,r_k\in J_M}\ 
\sum_{(r_1^{,},r_2^{,},\ldots,r_k^{,}):\
\{r_1^{,},r_2^{,},\ldots,r_k^{,}\}=\{r_1,r_2,\ldots,r_k\}}
\hspace{-2mm}{\sf M}\Biggl\{
\Biggl\Vert\ \Phi_k(Z)\left(\ldots \left(\Phi_2(Z) \left(\Phi_1(Z)
e_{r_1} \right) 
e_{r_2} \right) \ldots \right) e_{r_k} \Biggr\Vert_H\hspace{-1mm}\times \Biggr.
$$

\vspace{5mm}
$$
~~~~~~~~~~\times
\Biggl\Vert\ \Phi_k(Z)\left(\ldots \left(\Phi_2(Z) \left(\Phi_1(Z)
e_{r_1^{,}} \right) 
e_{r_2^{,}} \right) \ldots \right) e_{r_k^{,}}\ \Biggr\Vert_H\
\sqrt{\lambda_{r_1}\lambda_{r_2}\ldots \lambda_{r_k}}
\sqrt{\lambda_{r_1^{,}}\lambda_{r_2^{,}}\ldots \lambda_{r_k^{,}}}\times
$$

\vspace{5mm}
$$
\times
\left|{\sf M}\Biggl\{\Biggl(J[\psi^{(k)}]_{T,t}^{(r_1 r_2\ldots r_k)}-
J[\psi^{(k)}]_{T,t}^{(r_1 r_2\ldots r_k)p_1,\ldots,p_k}\Biggr)\times
\Biggr.\right.
$$

\vspace{5mm}

$$
\Biggl.\left.\Biggl.
~~~~~~~~~~~~~~~~~~~~~~~~~~~~~~~~~\times
\Biggl(J[\psi^{(k)}]_{T,t}^{(r_1^{,} r_2^{,}\ldots r_k^{,})}-
J[\psi^{(k)}]_{T,t}^{(r_1^{,} r_2^{,}\ldots r_k^{,})p_1,\ldots,p_k}\Biggr)
\biggl.\biggr|{\bf F}_t\Biggr\}
\right|
\Biggr\}\le
$$

\vspace{14mm}
$$
\le L_k
\sum_{r_1,r_2,\ldots,r_k\in J_M}\ \ 
\sum_{(r_1^{,},r_2^{,},\ldots,r_k^{,}):\ 
\{r_1^{,},r_2^{,},\ldots,r_k^{,}\}=\{r_1,r_2,\ldots,r_k\}}\ \
\sqrt{\lambda_{r_1}\lambda_{r_2}\ldots \lambda_{r_k}}
\sqrt{\lambda_{r_1^{,}}\lambda_{r_2^{,}}\ldots \lambda_{r_k^{,}}}\times
$$

\vspace{3mm}
$$
\hspace{-15mm}
\times{\sf M}\Biggl\{\left|\Biggl(J[\psi^{(k)}]_{T,t}^{(r_1 r_2\ldots r_k)}-
J[\psi^{(k)}]_{T,t}^{(r_1 r_2\ldots r_k)p_1,\ldots,p_k}\Biggr)\times
\right.\Biggr.
$$

\vspace{3mm}
$$
~~~~~~~~~~~~~~~~~~~~~~~~~~~\Biggl.\left.\times
\Biggl(J[\psi^{(k)}]_{T,t}^{(r_1^{,} r_2^{,}\ldots r_k^{,})}-
J[\psi^{(k)}]_{T,t}^{(r_1^{,} r_2^{,}\ldots r_k^{,})p_1,\ldots,p_k}\Biggr)
\right|\Biggr\}\le
$$

\vspace{14mm}
$$
\le L_k
\sum_{r_1,r_2,\ldots,r_k\in J_M}\ \ 
\sum_{(r_1^{,},r_2^{,},\ldots,r_k^{,}):\ 
\{r_1^{,},r_2^{,},\ldots,r_k^{,}\}=\{r_1,r_2,\ldots,r_k\}}\ \
\sqrt{\lambda_{r_1}\lambda_{r_2}\ldots \lambda_{r_k}}
\sqrt{\lambda_{r_1^{,}}\lambda_{r_2^{,}}\ldots \lambda_{r_k^{,}}}\times
$$

\vspace{3mm}
$$
\times
\left({\sf M}\left\{\Biggl(J[\psi^{(k)}]_{T,t}^{(r_1 r_2\ldots r_k)}-
J[\psi^{(k)}]_{T,t}^{(r_1 r_2\ldots r_k)p_1,\ldots,p_k}\Biggr)^2\right\}
\right)^{1/2}\times
$$

\vspace{3mm}
$$
~~~~~~~~~~~~~~~~~~~~~~~~~~~~~~~~~~~\times
\left({\sf M}\left\{
\Biggl(J[\psi^{(k)}]_{T,t}^{(r_1^{,} r_2^{,}\ldots r_k^{,})}-
J[\psi^{(k)}]_{T,t}^{(r_1^{,} r_2^{,}\ldots r_k^{,})p_1,\ldots,p_k}\Biggr)^2
\right\}
\right)^{1/2}\le
$$

\vspace{13mm}
$$
\le L_k
\sum_{r_1,r_2,\ldots,r_k\in J_M}\ \ 
\sum_{(r_1^{,},r_2^{,},\ldots,r_k^{,}):\ 
\{r_1^{,},r_2^{,},\ldots,r_k^{,}\}=\{r_1,r_2,\ldots,r_k\}}\ \
\sqrt{\lambda_{r_1}\lambda_{r_2}\ldots \lambda_{r_k}}
\sqrt{\lambda_{r_1^{,}}\lambda_{r_2^{,}}\ldots \lambda_{r_k^{,}}}\times
$$

\vspace{3mm}
$$
~~~~~~~~~~~~~~\times
\left(k!\left(I_k-\sum_{j_1=0}^{p_1}\ldots
\sum_{j_k=0}^{p_k}C^2_{j_k\ldots j_1}\right)\right)^{1/2}
\left(k!\left(I_k-\sum_{j_1=0}^{p_1}\ldots
\sum_{j_k=0}^{p_k}C^2_{j_k\ldots j_1}\right)\right)^{1/2}\le
$$

\vspace{12mm}
$$
\le L_k
\sum_{r_1,r_2,\ldots,r_k\in J_M} k!\ 
\lambda_{r_1}\lambda_{r_2}\ldots \lambda_{r_k}
\left(k!\left(I_k-\sum_{j_1=0}^{p_1}\ldots
\sum_{j_k=0}^{p_k}C^2_{j_k\ldots j_1}\right)\right)=
$$

\vspace{8mm}
$$
=L_k \left(k!\right)^2
\sum_{r_1,r_2,\ldots,r_k\in J_M}
\lambda_{r_1}\lambda_{r_2}\ldots \lambda_{r_k}
\left(I_k-\sum_{j_1=0}^{p_1}\ldots
\sum_{j_k=0}^{p_k}C^2_{j_k\ldots j_1}\right)\le
$$

\vspace{6mm}
$$
\le L_k \left(k!\right)^2
\left({\rm tr}\ Q\right)^k
\left(I_k-\sum_{j_1=0}^{p_1}\ldots
\sum_{j_k=0}^{p_k}C^2_{j_k\ldots j_1}\right),
$$

\vspace{6mm}
\noindent
where $\langle \cdot,\cdot\rangle_H$ is a scalar product in $H,$ and

\vspace{2mm}
$$
\sum_{(r_1^{,},r_2^{,},\ldots,r_k^{,}):\ 
\{r_1^{,},r_2^{,},\ldots,r_k^{,}\}=\{r_1,r_2,\ldots,r_k\}}
$$

\vspace{6mm}
\noindent
means the sum with respect to all
possible permutations
$(r_1^{,},r_2^{,},\ldots,r_k^{,})$ such that

\vspace{1mm}
$$
\{r_1^{,},r_2^{,},\ldots,r_k^{,}\}=\{r_1,r_2,\ldots,r_k\}.
$$

\vspace{7mm}

The transition from (\ref{zz1}) to (\ref{zz100}) 
is based on the following theorem.

\vspace{6mm}

{\bf Theorem 5}\ \cite{37a}-\cite{37aaxx}, \cite{2000}. {\it 
The following 
equality is true

\vspace{3mm}
$$
\hspace{-40mm}
{\sf M}
\Biggl\{\Biggl(
J[\psi^{(k)}]_{T,t}^{(r_1\ldots r_k)}-
J[\psi^{(k)}]_{T,t}^{(r_1\ldots r_k)p_1\ldots p_k}
\Biggr)\times\Biggr.
$$

\vspace{1mm}
\begin{equation}
\label{uuu2}
~~~~~~~~~~~~~~~~~~~~~\Biggl.\times
\Biggl(
J[\psi^{(k)}]_{T,t}^{(m_1\ldots m_k)}-
J[\psi^{(k)}]_{T,t}^{(m_1\ldots m_k)p_1\ldots p_k}
\Biggr)\biggl.\biggr|{\bf F}_t\Biggr\}=0
\end{equation}

\vspace{7mm}
\noindent
w. p. {\rm 1} for all $r_1,\ldots,r_k,m_1,\ldots,m_k\in J_M$ $(M\in\mathbb{N})$
such that $\{r_1,\ldots,r_k\}\ne 
\{m_1,\ldots,m_k\}$.}

\vspace{5mm}
{\bf Proof.} Using the standard moment properties of 
the It\^{o} stochastic integral,
we obtain

\vspace{2mm}
\begin{equation}
\label{uuu3}
{\sf M}
\biggl\{J[\psi^{(k)}]_{T,t}^{(r_1\ldots r_k)}
J[\psi^{(k)}]_{T,t}^{(m_1\ldots m_k)}\biggl.\biggr|{\bf F}_t
\biggr\}=0
\end{equation}

\vspace{5mm}
\noindent
w. p. 1\ for all $r_1,\ldots,r_k,m_1,\ldots,m_k\in J_M$ 
such that $(r_1,\ldots,r_k)\ne (m_1,\ldots,m_k),$ $M\in\mathbb{N}.$

From the
proof of Theorem {\rm 1.18} in \cite{37a} (Sect.~1.12)
it follows that

\vspace{3mm}
$$
\prod_{l=1}^k\zeta_{j_l}^{(r_l)}+\sum\limits_{m=1}^{[k/2]}
(-1)^m \times
\Biggr.
$$

\vspace{2mm}
$$
\times
\sum_{\stackrel{(\{\{g_1, g_2\}, \ldots, 
\{g_{2m-1}, g_{2m}\}\}, \{q_1, \ldots, q_{k-2m}\})}
{{}_{\{g_1, g_2, \ldots, 
g_{2m-1}, g_{2m}, q_1, \ldots, q_{k-2m}\}=\{1, 2, \ldots, k\}}}}
\prod\limits_{s=1}^m
{\bf 1}_{\{r_{g_{{}_{2s-1}}}=~r_{g_{{}_{2s}}}\ne 0\}}
\Biggl.{\bf 1}_{\{j_{g_{{}_{2s-1}}}=~j_{g_{{}_{2s}}}\}}
\prod_{l=1}^{k-2m}\zeta_{j_{q_l}}^{(r_{q_l})}=
$$

\vspace{1mm}
\begin{equation}
\label{ttt2}
=\sum\limits_{(j_1,\ldots,j_k)}
\int\limits_t^T \phi_{j_k}(t_k)
\ldots
\int\limits_t^{t_{2}}\phi_{j_{1}}(t_{1})
d{\bf w}_{t_1}^{(r_1)}\ldots
d{\bf w}_{t_k}^{(r_k)}\ \ \ {\rm w.\ p.\ 1,}
\end{equation}

\vspace{4mm}
\noindent
where 
$$
\sum\limits_{(j_1,\ldots,j_k)}
$$ 

\vspace{2mm}
\noindent
means the sum with respect to all
possible permutations
$(j_1,\ldots,j_k)$. At the same time if 
$j_l$ swapped  with $j_q$ in the permutation $(j_1,\ldots,j_k)$,
then $r_l$ swapped  with $r_q$ in the permutation
$(r_1,\ldots,r_k);$
another notations are the same as in Theorem 2.

Using (\ref{leto6000}) and (\ref{ttt2}), we get 

\vspace{2mm}
\begin{equation}
\label{est11}
J[\psi^{(k)}]_{T,t}^{(r_1\ldots r_k)p_1\ldots p_k}=
\sum\limits_{j_1=0}^{p_1}\ldots
\sum\limits_{j_k=0}^{p_k}
C_{j_k\ldots j_1}\sum\limits_{(j_1,\ldots,j_k)}
\int\limits_t^T \phi_{j_k}(t_k)
\ldots
\int\limits_t^{t_{2}}\phi_{j_{1}}(t_{1})
d{\bf w}_{t_1}^{(r_1)}\ldots
d{\bf w}_{t_k}^{(r_k)},
\end{equation}

\vspace{5.5mm}
\noindent
where notations are the same as in (\ref{ttt2}).

Then w. p. 1 

\vspace{2mm}
$$
{\sf M}
\biggl\{J[\psi^{(k)}]_{T,t}^{(m_1\ldots m_k)}
J[\psi^{(k)}]_{T,t}^{(r_1\ldots r_k)p_1\ldots p_k}
\biggl.\biggr|{\bf F}_t\biggr\}=
$$

\vspace{3mm}
$$
=\sum_{j_1=0}^{p_1}\ldots\sum_{j_k=0}^{p_k}
C_{j_k\ldots j_1}\times
$$

\vspace{3mm}
$$
\times{\sf M}\left\{J[\psi^{(k)}]_{T,t}^{(m_1\ldots m_k)}
\sum\limits_{(j_1,\ldots,j_k)}
\int\limits_t^T \phi_{j_k}(t_k)
\ldots
\int\limits_t^{t_{2}}\phi_{j_{1}}(t_{1})
d{\bf w}_{t_1}^{(r_1)}\ldots
d{\bf w}_{t_k}^{(r_k)}
\biggl.\biggr|{\bf F}_t\right\}.
$$

\vspace{8mm}

From the standard moment properties of the It\^{o} stochastic integral 
it follows that

\vspace{2mm}
$$
{\sf M}\left\{J[\psi^{(k)}]_{T,t}^{(m_1\ldots m_k)}
\sum\limits_{(j_1,\ldots,j_k)}
\int\limits_t^T \phi_{j_k}(t_k)
\ldots
\int\limits_t^{t_{2}}\phi_{j_{1}}(t_{1})
d{\bf w}_{t_1}^{(r_1)}\ldots
d{\bf w}_{t_k}^{(r_k)}
\biggl.\biggr|{\bf F}_t\right\}=0
$$

\vspace{5.5mm}
\noindent
w. p. 1\ for all $r_1,\ldots,r_k,m_1,\ldots,m_k\in J_M$ 
such that $\{r_1,\ldots,r_k\}\ne \{m_1,\ldots,m_k\},$ $M\in\mathbb{N}.$

Then 

\vspace{1mm}
\begin{equation}
\label{uuu4}
{\sf M}
\biggl\{J[\psi^{(k)}]_{T,t}^{(m_1\ldots m_k)}
J[\psi^{(k)}]_{T,t}^{(r_1\ldots r_k)p_1\ldots p_k}
\biggl.\biggr|{\bf F}_t\biggr\}=0
\end{equation}

\vspace{5.5mm}
\noindent
w. p. 1 for all $r_1,\ldots,r_k,m_1,\ldots,m_k\in J_M$ $(M\in\mathbb{N})$
such that $\{r_1,\ldots,r_k\}\ne 
\{m_1,\ldots,m_k\}$.

From (\ref{est11}) it follows that

\vspace{5mm}
$$
{\sf M}\biggl\{J[\psi^{(k)}]_{T,t}^{(r_1\ldots r_k)p_1,\ldots,p_k}
J[\psi^{(k)}]_{T,t}^{(m_1\ldots m_k)p_1,\ldots,p_k} 
\biggl.\biggr|{\bf F}_t\biggr\}
=
$$

\vspace{3mm}
$$
=\sum_{j_1=0}^{p_1}\ldots\sum_{j_k=0}^{p_k}
C_{j_k\ldots j_1}
\sum_{{q}_1=0}^{p_1}\ldots\sum_{{q}_k=0}^{p_k}
C_{{q}_k\ldots {q}_1}\times
$$

\vspace{2mm}
$$
\hspace{-40mm}\times
{\sf M}\left\{\left(
\sum\limits_{(j_1,\ldots,j_k)}
\int\limits_t^T \phi_{j_k}(t_k)
\ldots
\int\limits_t^{t_{2}}\phi_{j_{1}}(t_{1})
d{\bf w}_{t_1}^{(r_1)}\ldots
d{\bf w}_{t_k}^{(r_k)}\right) \times\right.
$$

\vspace{1mm}
\begin{equation}
\label{uuu5}
~~~~~~~~~~~~~~~~~~~~~~~~~~~~~~~~\times
\left.\left(
\sum\limits_{(q_1,\ldots,q_k)}
\int\limits_t^T \phi_{q_k}(t_k)
\ldots
\int\limits_t^{t_{2}}\phi_{q_{1}}(t_{1})
d{\bf w}_{t_1}^{(m_1)}\ldots
d{\bf w}_{t_k}^{(m_k)}\right) \biggl.\biggr|{\bf F}_t\right\}=0
\end{equation}

\vspace{8mm}
\noindent
w. p. 1 for all $r_1,\ldots,r_k,m_1,\ldots,m_k\in J_M$ $(M\in\mathbb{N})$
such that $\{r_1,\ldots,r_k\}\ne 
\{m_1,\ldots,m_k\}$.

From (\ref{uuu3}), (\ref{uuu4}), and (\ref{uuu5}) 
we obtain (\ref{uuu2}).
Theorem 5 is proved.

\vspace{1mm}

\vspace{5mm}
{\bf Corollary 1}\ \cite{37a}-\cite{37aaxx}, \cite{2000}. {\it The following equality is true

\vspace{2mm}
$$
{\sf M}
\Biggl\{\Biggl(
J[\psi^{(k)}]_{T,t}^{(r_1\ldots r_k)}-
J[\psi^{(k)}]_{T,t}^{(r_1\ldots r_k)p_1\ldots p_k}
\Biggr)
\Biggl(
J[\psi^{(l)}]_{T,t}^{(m_1\ldots m_l)}-
J[\psi^{(l)}]_{T,t}^{(m_1\ldots m_l)q_1\ldots q_l}
\Biggr)\biggl.\biggr|{\bf F}_t\Biggr\}=0
$$

\vspace{5mm}
w. p. {\rm 1}\ for all $l=1, 2, \ldots, k-1$, and 
$r_1,\ldots,r_k,m_1,\ldots,m_l\in J_M,$ $p_1,\ldots,p_k,q_1,\ldots,q_l
=0,1,2,\ldots $}

\vspace{10mm}

\section{Approximation of Some Iterated Stochastic Integrals 
of Second and Third Miltiplicity with
Respect to the $Q$-Wiener Process}

\vspace{5mm}

This section is devoted to the approximation of iterated
stochastic integrals of the following form (see Sect.~1)

\vspace{-2mm}
\begin{equation}
\label{abc1}
I_0[B(Z),F(Z)]_{T,t}^M=
\int\limits_{t}^{T}B'(Z) \left(
\int\limits_{t}^{t_2}F(Z) dt_1 \right) d{\bf W}_{t_2}^M,
\end{equation}

\begin{equation}
\label{kkk1}
I_1[B(Z),F(Z)]_{T,t}^M=\int\limits_{t}^{T}F'(Z)
\left(\int\limits_{t}^{t_2} B(Z) d{\bf W}_{t_1}^M \right) dt_2,
\end{equation}

\begin{equation}
\label{kkk2}
I_2[B(Z)]_{T,t}^M=\int\limits_{t}^{T}B''(Z) \left(
\int\limits_{t}^{t_2}B(Z) d{\bf W}_{t_1}^M,
\int\limits_{t}^{t_2}B(Z) d{\bf W}_{t_1}^M 
\right) d{\bf W}_{t_2}^M.\
\end{equation}

\vspace{5mm}

Let conditions 1 and 2 of Theorem 4 be fulfilled.
Let $B''(v)(B(v),B(v))$ 
be a 3-linear Hilbert--Schmidt operator mapping from
$U_0\times U_0\times U_0$ to $H$ for all $v\in H$.

Then we have w.\ p.\ 1 (see (\ref{xx605}))

\vspace{1mm}
\begin{equation}
\label{abc2}
I_0[B(Z),F(Z)]_{T,t}^M=
\sum_{r_1\in J_M}B'(Z)F(Z)e_{r_1}\sqrt{\lambda_{r_1}}I_{(01)T,t}^{(0 r_1)},
\end{equation}

\vspace{2mm}
\begin{equation}
\label{da0}
I_1[B(Z),F(Z)]_{T,t}^M=
\sum_{r_1\in J_M}F'(Z)(B(Z)e_{r_1})\sqrt{\lambda_{r_1}}I_{(10)T,t}^{(r_1 0)},
\end{equation}

\vspace{4mm}
$$
I_2[B(Z)]_{T,t}^M=
\sum_{r_1,r_2,r_3\in J_M}B''(Z)\left(B(Z)e_{r_1}, B(Z)e_{r_2}\right)e_{r_3}
\sqrt{\lambda_{r_1}\lambda_{r_2}\lambda_{r_3}}\times
$$

\begin{equation}
\label{jjj1}
\times \int\limits_t^T
\left(\int\limits_t^s d{\bf w}_{\tau}^{(r_1)}   
\int\limits_t^s d{\bf w}_{\tau}^{(r_2)}\right)
d{\bf w}_s^{(r_3)}.
\end{equation}

\vspace{4mm}

Using the It\^{o} formula, we obtain

\vspace{1mm}
\begin{equation}
\label{da12}
\int\limits_t^s d{\bf w}_{\tau}^{(r_1)}   
\int\limits_t^s d{\bf w}_{\tau}^{(r_2)}
=
I_{(11)s,t}^{(r_1 r_2)}+
I_{(11)s,t}^{(r_2 r_1)}+
{\bf 1}_{\{r_1=r_2\}}(s-t)\ \ \ \hbox{w.\ p.\ 1.}
\end{equation}

\vspace{2mm}

From (\ref{da12}) we have

\vspace{1mm}
\begin{equation}
\label{jjj2}
\int\limits_t^T
\left(\int\limits_t^s d{\bf w}_{\tau}^{(r_1)}   
\int\limits_t^s d{\bf w}_{\tau}^{(r_2)}\right)
d{\bf w}_s^{(r_3)}=
I_{(111)T,t}^{(r_1 r_2 r_3)}+
I_{(111)T,t}^{(r_2 r_1 r_3)}+
{\bf 1}_{\{r_1=r_2\}}I_{(01)T,t}^{(0 r_3)}\ \ \ \hbox{w.\ p.\ 1.}
\end{equation}

\vspace{4mm}

Note that in 
(\ref{abc2}), (\ref{da0}), (\ref{da12}), and (\ref{jjj2}) 
we use the notations from Sect.~4.

After substituting (\ref{jjj2}) into (\ref{jjj1}), we have

\vspace{5mm}
$$
I_2[B(Z)]_{T,t}^M=
\sum_{r_1,r_2,r_3\in J_M}B''(Z)\left(B(Z)e_{r_1}, B(Z)e_{r_2}\right)e_{r_3}
\sqrt{\lambda_{r_1}\lambda_{r_2}\lambda_{r_3}} \times
$$

\vspace{-1mm}
\begin{equation}
\label{da1}
\times
\left(I_{(111)T,t}^{(r_1 r_2 r_3)}+
I_{(111)T,t}^{(r_2 r_1 r_3)}+
{\bf 1}_{\{r_1=r_2\}}I_{(01)T,t}^{(0 r_3)}\right)\ \ \ \hbox{w.\ p.\ 1.}
\end{equation}

\vspace{7mm}

Taking into account (\ref{opp1}) and (\ref{opp2}), we put for $q=1$

\vspace{4mm}
\begin{equation}
\label{opp3}
I_{(01)T,t}^{(0 r_3)q}=
I_{(01)T,t}^{(0 r_3)}=\frac{(T-t)^{3/2}}{2}\biggl(\zeta_0^{(r_3)}+
\frac{1}{\sqrt{3}}\zeta_1^{(r_3)}\biggr)\ \ \
(q=1)\ \ \ \hbox{w.\ p.\ 1},
\end{equation}

\vspace{2mm}
\begin{equation}
\label{opp4}
I_{(10)T,t}^{(r_1 0)q}=
I_{(10)T,t}^{(r_1 0)}=\frac{(T-t)^{3/2}}{2}\biggl(\zeta_0^{(r_1)}-
\frac{1}{\sqrt{3}}\zeta_1^{(r_1)}\biggr)\ \ \ 
(q=1)\ \ \ \hbox{w.\ p.\ 1},
\end{equation}

\vspace{7mm}
\noindent
where 
$I_{(01)T,t}^{(0 r_3)q},$
$I_{(10)T,t}^{(r_1 0)q}$ denote the approximations of
corresponding
iterated It\^{o} stochastic integrals.

Denote by
$I_0[B(Z),F(Z)]_{T,t}^{M,q},$ $I_1[B(Z),F(Z)]_{T,t}^{M,q},$ 
$I_2[B(Z)]_{T,t}^{M,q}$ the approximations of iterated stochastic integrals
(\ref{abc2}), (\ref{da0}), (\ref{da1})

\vspace{2mm}

\begin{equation}
\label{abc3}
I_0[B(Z),F(Z)]_{T,t}^{M,q}=
\sum_{r_1\in J_M}B'(Z)F(Z)e_{r_1}\sqrt{\lambda_{r_1}}I_{(01)T,t}^{(0 r_1)q},
\end{equation}

\vspace{2mm}

\begin{equation}
\label{da3}
I_1[B(Z),F(Z)]_{T,t}^{M,q}=
\sum_{r_1\in J_M}F'(Z)(B(Z)e_{r_1})\sqrt{\lambda_{r_1}}I_{(10)T,t}^{(r_1 0)q},
\end{equation}

\vspace{4mm}

$$
I_2[B(Z)]_{T,t}^{M,q}=
\sum_{r_1,r_2,r_3\in J_M}B''(Z)\left(B(Z)e_{r_1}, B(Z)e_{r_2}\right)e_{r_3}
\sqrt{\lambda_{r_1}\lambda_{r_2}\lambda_{r_3}} \times
$$

\begin{equation}
\label{da4}
\times
\left(I_{(111)T,t}^{(r_1 r_2 r_3)q}+
I_{(111)T,t}^{(r_2 r_1 r_3)q}+
{\bf 1}_{\{r_1=r_2\}}I_{(01)T,t}^{(0 r_3)q}\right),
\end{equation}

\vspace{7mm}
\noindent
where $q\ge 1,$ and the approximations $I_{(111)T,t}^{(r_1 r_2 r_3)q},$
$I_{(111)T,t}^{(r_2 r_1 r_3)q}$ are defined by 
(\ref{kr1}).

From (\ref{abc2}), (\ref{da0}), (\ref{da1}),
(\ref{abc3})--(\ref{da4}) it follows that

\vspace{4mm}

$$
I_0[B(Z),F(Z)]_{T,t}^{M}-I_0[B(Z),F(Z)]_{T,t}^{M,q}=0\ \ \ \hbox{w.\ p.\ 1,}
$$

\vspace{2mm}

$$
I_1[B(Z),F(Z)]_{T,t}^{M}-I_1[B(Z),F(Z)]_{T,t}^{M,q}=0\ \ \ \hbox{w.\ p.\ 1,}
$$

\vspace{4mm}

$$
I_2[B(Z)]_{T,t}^{M}-I_2[B(Z)]_{T,t}^{M,q}=
\sum_{r_1,r_2,r_3\in J_M}B''(Z)\left(B(Z)e_{r_1}, B(Z)e_{r_2}\right)e_{r_3}
\sqrt{\lambda_{r_1}\lambda_{r_2}\lambda_{r_3}} \times
$$

\vspace{1mm}
$$
\times
\left(\left(I_{(111)T,t}^{(r_1 r_2 r_3)}-I_{(111)T,t}^{(r_1 r_2 r_3)q}\right)+
\left(I_{(111)T,t}^{(r_2 r_1 r_3)}-I_{(111)T,t}^{(r_2 r_1 r_3)q}
\right)
\right)\ \ \ \hbox{w.\ p.\ 1.}
$$

\vspace{9mm}

Repeating with an insignificant 
modification the proof 
of Theorem 4 for the case $k=3$, we obtain

\vspace{3mm}
$$
\hspace{-55mm}
{\sf M}\left\{\biggl\Vert
I_2[B(Z)]_{T,t}^{M}-I_2[B(Z)]_{T,t}^{M,q}\biggr\Vert_H^2\right\}\le
$$

\vspace{2mm}
$$
~~~~~~~~~~~~~~~~~~~~~~~~~~~~~~~~~~~~~~~~~~~~~\le 4C(3!)^2
\left({\rm tr}\ Q\right)^3
\Biggl(\frac{(T-t)^{3}}{6}-\sum_{j_1,j_2,j_3=0}^{q}
C_{j_3j_2j_1}^2\Biggr),
$$

\vspace{7mm}
\noindent
where here and further constant $C$ has the same meaning 
as constant $L_k$ in Theorem 4 
($k$ is the multiplicity of the iterated stochastic integral),
and

\vspace{7mm}
$$
C_{j_3j_2j_1}=\frac{\sqrt{(2j_1+1)(2j_2+1)(2j_3+1)}(T-t)^{3/2}}{8}\bar
C_{j_3j_2j_1},
$$

\vspace{3mm}
$$
\bar C_{j_3j_2j_1}=\int\limits_{-1}^{1}P_{j_3}(z)
\int\limits_{-1}^{z}P_{j_2}(y)
\int\limits_{-1}^{y}
P_{j_1}(x)dx dy dz,
$$

\vspace{5mm}
\noindent
where $P_j(x)$ is the Legendre polynomial.

\vspace{5mm}

\section{Approximation of Some Iterated Stochastic Integrals
of Third and Fourth Miltiplicity with
Respect to the $Q$-Wiener Process}

\vspace{5mm}

In this section, we consider the approximation of iterated
stochastic integrals of the following form (see Sect.~1)

\vspace{2mm}
$$
I_3[B(Z)]_{T,t}^M=\int\limits_{t}^{T}B'''(Z) \left(
\int\limits_{t}^{t_2}B(Z) d{\bf W}_{t_1}^M,
\int\limits_{t}^{t_2}B(Z) d{\bf W}_{t_1}^M,
\int\limits_{t}^{t_2}B(Z) d{\bf W}_{t_1}^M 
\right) d{\bf W}_{t_2}^M,\
$$

\vspace{3mm}
$$
I_4[B(Z)]_{T,t}^M=\int\limits_{t}^{T}B'(Z) \left(
\int\limits_{t}^{t_3}B''(Z) \left(
\int\limits_{t}^{t_2}B(Z) d{\bf W}_{t_1}^M,
\int\limits_{t}^{t_2}B(Z) d{\bf W}_{t_1}^M 
\right) d{\bf W}_{t_2}^M\right)d{\bf W}_{t_3}^M,\
$$

\vspace{3mm}
$$
I_5[B(Z)]_{T,t}^M=\int\limits_{t}^{T}B''(Z) \left(
\int\limits_{t}^{t_3}B(Z)d{\bf W}_{t_1}^M,
\int\limits_{t}^{t_3}B'(Z)\left(
\int\limits_{t}^{t_2}B(Z) d{\bf W}_{t_1}^M
\right) d{\bf W}_{t_2}^M\right)d{\bf W}_{t_3}^M,\
$$

\vspace{3mm}
$$
I_6[B(Z),F(Z)]_{T,t}^M=\int\limits_{t}^{T}F'(Z)\left(
\int\limits_{t}^{t_3}B'(Z) \left(
\int\limits_{t}^{t_2}B(Z) d{\bf W}_{t_1}^M
\right) d{\bf W}_{t_2}^M 
\right) dt_3,
$$

\vspace{3mm}
$$
I_7[B(Z),F(Z)]_{T,t}^M=\int\limits_{t}^{T}F''(Z) \left(
\int\limits_{t}^{t_2}B(Z) d{\bf W}_{t_1}^M,
\int\limits_{t}^{t_2}B(Z) d{\bf W}_{t_1}^M 
\right) dt_2,
$$

\vspace{3mm}
$$
I_8[B(Z),F(Z)]_{T,t}^M=\int\limits_{t}^{T}B''(Z) \left(
\int\limits_{t}^{t_2}F(Z) dt_1,
\int\limits_{t}^{t_2}B(Z) d{\bf W}_{t_1}^M 
\right) d{\bf W}_{t_2}^M.
$$

\vspace{7mm}

Consider the stochastic integral $I_3[B(Z)]_{T,t}^M.$
Let conditions 1 and 2 of Theorem 4 be fulfilled.
Let $B'''(v)(B(v),B(v),B(v))$
be a 4-linear Hilbert--Schmidt operator 
mapping from
$U_0\times U_0\times U_0\times U_0$
to $H$
for all $v\in H$.

We have (see (\ref{xx605}))

\vspace{2mm}
$$
I_3[B(Z)]_{T,t}^M=
\sum_{r_1,r_2,r_3,r_4\in J_M}
B'''(Z)\left(B(Z)e_{r_1}, B(Z)e_{r_2},B(Z)e_{r_3}\right)e_{r_4}
\sqrt{\lambda_{r_1}\lambda_{r_2}\lambda_{r_3}\lambda_{r_4}}\times
$$

\begin{equation}
\label{jjjk1}
\times\int\limits_t^T
\left(\int\limits_t^s d{\bf w}_{\tau}^{(r_1)}   
\int\limits_t^s d{\bf w}_{\tau}^{(r_2)}
\int\limits_t^s d{\bf w}_{\tau}^{(r_3)}\right)
d{\bf w}_s^{(r_4)}\ \ \ \hbox{w.\ p.\ 1.}
\end{equation}

\vspace{7mm}

From \cite{37} (pp.~A.438--A.439) (also see \cite{37a}-\cite{37aaxx})
or using the It\^{o} formula
we obtain

\vspace{2mm}
$$
I_{(1)s,t}^{(r_1)}
I_{(1)s,t}^{(r_2)}
I_{(1)s,t}^{(r_3)}=
$$

\vspace{1mm}

$$
=
I_{(111)s,t}^{(r_1r_2r_3)}+
I_{(111)s,t}^{(r_1r_3r_2)}+
I_{(111)s,t}^{(r_2r_1r_3)}+
I_{(111)s,t}^{(r_2r_3r_1)}+
I_{(111)s,t}^{(r_3r_1r_2)}+
I_{(111)s,t}^{(r_3r_2r_1)}+
$$

\vspace{1mm}
$$
+{\bf 1}_{\{r_1=r_2\}}\left(
I_{(10)s,t}^{(r_3 0)}+I_{(01)s,t}^{(0 r_3)}\right)+
{\bf 1}_{\{r_1=r_3\}}\left(
I_{(10)s,t}^{(r_2 0)}+I_{(01)s,t}^{(0 r_2)}\right)+
$$

\vspace{1mm}
$$
+{\bf 1}_{\{r_2=r_3\}}\left(
I_{(10)s,t}^{(r_1 0)}+I_{(01)s,t}^{(0 r_1)}\right)=
$$

\begin{equation}
\label{q12}
=\sum\limits_{(r_1,r_2,r_3)}
I_{(111)s,t}^{(r_1r_2r_3)}+
(s-t)\left({\bf 1}_{\{r_2=r_3\}}I_{(1)s,t}^{(r_1)}+
{\bf 1}_{\{r_1=r_3\}}I_{(1)s,t}^{(r_2)}+
{\bf 1}_{\{r_1=r_2\}}I_{(1)s,t}^{(r_3)}\right)\ \ \ \hbox{w.\ p.\ 1,}
\end{equation}

\vspace{5mm}
\noindent
where 
$$
\sum\limits_{(r_1,r_2,r_3)}
$$ 

\vspace{2mm}
\noindent
means the sum with respect to all
possible permutations
$(r_1,r_2,r_3)$ and we use the notations from Sect.~4.

After substituting (\ref{q12}) into (\ref{jjjk1}), we obtain

\vspace{5mm}

$$
I_3[B(Z)]_{T,t}^M=
\sum_{r_1,r_2,r_3,r_4\in J_M}
B'''(Z)\left(B(Z)e_{r_1}, B(Z)e_{r_2},B(Z)e_{r_3}\right)e_{r_4}
\sqrt{\lambda_{r_1}\lambda_{r_2}\lambda_{r_3}\lambda_{r_4}}\times
$$

\begin{equation}
\label{jjjk2}
\times 
\Biggl(\sum\limits_{(r_1,r_2,r_3)}
I_{(1111)T,t}^{(r_1r_2r_3r_4)}-
{\bf 1}_{\{r_1=r_2\}}J_{(01)T,t}^{(r_3r_4)}-
{\bf 1}_{\{r_1=r_3\}}J_{(01)T,t}^{(r_2r_4)}-
{\bf 1}_{\{r_2=r_3\}}J_{(01)T,t}^{(r_1r_4)}\Biggr)\ \ \ \hbox{w.\ p.\ 1,}
\end{equation}

\vspace{5mm}
\noindent
where 

\vspace{-2mm}
\begin{equation}
\label{ros1}
J_{(01)T,t}^{(r_1r_2)}=
\int\limits_t^T
(t-s)\int\limits_t^s d{\bf w}_{\tau}^{(r_1)}   
d{\bf w}_s^{(r_2)}.
\end{equation}

\vspace{4mm}

Denote by
$I_3[B(Z)]_{T,t}^{M,q}$ the approximation of the iterated stochastic integral
(\ref{jjjk2}),
which has the following form

\vspace{5mm}

$$
I_3[B(Z)]_{T,t}^{M,q}=
\sum_{r_1,r_2,r_3,r_4\in J_M}
B'''(Z)\left(B(Z)e_{r_1}, B(Z)e_{r_2},B(Z)e_{r_3}\right)e_{r_4}
\sqrt{\lambda_{r_1}\lambda_{r_2}\lambda_{r_3}\lambda_{r_4}}\times
$$

\begin{equation}
\label{12346}
\times 
\left(\sum\limits_{(r_1,r_2,r_3)}
I_{(1111)T,t}^{(r_1r_2r_3r_4)q}-
{\bf 1}_{\{r_1=r_2\}}J_{(01)T,t}^{(r_3r_4)q}-
{\bf 1}_{\{r_1=r_3\}}J_{(01)T,t}^{(r_2r_4)q}-
{\bf 1}_{\{r_2=r_3\}}J_{(01)T,t}^{(r_1r_4)q}\right),
\end{equation}

\vspace{7mm}
\noindent
where 
the approximations $I_{(1111)T,t}^{(r_1r_2r_3r_4)q},$
$J_{(01)T,t}^{(r_1r_2)q}$
are based on Theorems 1, 2 and Le\-gen\-dre polynomials.

The approximation $J_{(01)T,t}^{(r_1 r_2)q}$ of the stochastic integral
$J_{(01)T,t}^{(r_1 r_2)}$ $(r_1,r_2=1,\ldots,M),$
which is based on Theorems 1, 2 and 
Le\-gendre polynomials has the following form (see 
\cite{37} (formula (6.91), p.~A.544) or 
\cite{34} (formula (5.7), p.~A.249))

\vspace{3mm}
$$
J_{(01)T,t}^{(r_1 r_2)q}=
-\frac{T-t}{2}
I_{(11)T,t}^{(r_1 r_2)q}
-\frac{(T-t)^2}{4}\Biggl(
\frac{1}{\sqrt{3}}\zeta_0^{(r_1)}\zeta_1^{(r_2)}+\Biggr.
$$

\vspace{1mm}
\begin{equation}
\label{vini0}
+\Biggl.\sum_{i=0}^{q}\Biggl(
\frac{(i+2)\zeta_i^{(r_1)}\zeta_{i+2}^{(r_2)}
-(i+1)\zeta_{i+2}^{(r_1)}\zeta_{i}^{(r_2)}}
{\sqrt{(2i+1)(2i+5)}(2i+3)}-
\frac{\zeta_i^{(r_1)}\zeta_{i}^{(r_2)}}{(2i-1)(2i+3)}\Biggr)\Biggr),
\end{equation}

\vspace{8mm}

\begin{equation}
\label{vini}
I_{(11)T,t}^{(r_1 r_2)q}=
\frac{T-t}{2}\left(\zeta_0^{(r_1)}\zeta_0^{(r_2)}+\sum_{i=1}^{q}
\frac{1}{\sqrt{4i^2-1}}\biggl(
\zeta_{i-1}^{(r_1)}\zeta_{i}^{(r_2)}-
\zeta_i^{(r_1)}\zeta_{i-1}^{(r_2)}\biggr)-{\bf 1}_{\{r_1=r_2\}}
\right),
\end{equation}

\vspace{8mm}
\noindent
where notations are the same as in Theorems 1, 2.

Moreover (see \cite{37} (formula (6.106), p.~A.551)
or \cite{34} (formula (5.19), p.~A.252--A.253)),

\vspace{3mm}

$$
{\sf M}\left\{\left(J_{(01)T,t}^{(r_1 r_2)}-
J_{(01)T,t}^{(r_1 r_2)q}\right)^2\right\}=\frac{(T-t)^4}{16}\times
$$

\begin{equation}
\label{987}
\times\Biggl(\frac{5}{9}-
2\sum_{i=2}^q\frac{1}{4i^2-1}-
\sum_{i=1}^q
\frac{1}{(2i-1)^2(2i+3)^2}-
\sum_{i=0}^q\frac{(i+2)^2+(i+1)^2}{(2i+1)(2i+5)(2i+3)^2}
\Biggr)\ \ \ (r_1\ne r_2).
\end{equation}

\vspace{7mm}

From (\ref{star00011}), (\ref{12345}) we obtain

\vspace{3mm}

$$
{\sf M}\left\{\left(J_{(01)T,t}^{(r_1 r_2)}-
J_{(01)T,t}^{(r_1 r_2)q}\right)^2\right\}\le 
$$

$$
\le\frac{(T-t)^4}{8}\Biggl(\frac{5}{9}-
2\sum_{i=2}^q\frac{1}{4i^2-1}-
\sum_{i=1}^q
\frac{1}{(2i-1)^2(2i+3)^2}-
\sum_{i=0}^q\frac{(i+2)^2+(i+1)^2}{(2i+1)(2i+5)(2i+3)^2}
\Biggr),
$$

\vspace{7mm}
\noindent
where $r_1,r_2=1,\ldots,M.$

From (\ref{jjjk2}), (\ref{12346})
it follows that

\vspace{3mm}

$$
I_3[B(Z)]_{T,t}^{M}-
I_3[B(Z)]_{T,t}^{M,q}=
$$

\vspace{2mm}
$$
=
\sum_{r_1,r_2,r_3,r_4\in J_M}
B'''(Z)\left(B(Z)e_{r_1}, B(Z)e_{r_2},B(Z)e_{r_3}\right)e_{r_4}
\sqrt{\lambda_{r_1}\lambda_{r_2}\lambda_{r_3}\lambda_{r_4}}\times
$$

\vspace{2mm}
$$
\times 
\Biggl(\sum\limits_{(r_1,r_2,r_3)}
\left(I_{(1111)T,t}^{(r_1r_2r_3r_4)}-I_{(1111)T,t}^{(r_1r_2r_3r_4)q}\right)-
{\bf 1}_{\{r_1=r_2\}}\left(J_{(01)T,t}^{(r_3r_4)}-J_{(01)T,t}^{(r_3r_4)q}
\right)-\Biggr.
$$

\begin{equation}
\label{12341}
\Biggl.-{\bf 1}_{\{r_1=r_3\}}
\left(J_{(01)T,t}^{(r_2r_4)}-J_{(01)T,t}^{(r_2r_4)q}
\right)-
{\bf 1}_{\{r_2=r_3\}}\left(J_{(01)T,t}^{(r_1r_4)}-
J_{(01)T,t}^{(r_1r_4)q}\right)\Biggr)\ \ \ \hbox{w.\ p.\ 1.}
\end{equation}

\vspace{7mm}

Repeating with an insignificant 
modification the proof 
of Theorem 4 for the cases $k=2, 4$, we obtain

\vspace{3mm}
$$
\hspace{-70mm}{\sf M}\Biggl\{\biggl\Vert
I_3[B(Z)]_{T,t}^{M}-I_3[B(Z)]_{T,t}^{M,q}\biggr\Vert_H^2\Biggr\}\le
$$

\vspace{2mm}
$$
~~~~~~~~~~~~~~~~~~~~~~~~~~~~~~~~~\le C\left({\rm tr}\ Q\right)^4
\Biggl(6^2(4!)^2
\Biggl(\frac{(T-t)^{4}}{24}-\sum_{j_1,j_2,j_3,j_4=0}^{q}
C_{j_4j_3j_2j_1}^2\Biggr)+3^2(2!)^2E_q\Biggr),
$$

\vspace{7mm}
\noindent
where 
$E_q$ is the right-hand side of (\ref{987}), and

\vspace{4mm}

\begin{equation}
\label{tq1}
C_{j_4j_3j_2j_1}
=\frac{\sqrt{(2j_1+1)(2j_2+1)(2j_3+1)(2j_4+1)}(T-t)^{2}}{16}\bar
C_{j_4j_3j_2j_1},
\end{equation}

\vspace{3mm}
$$
\bar C_{j_4j_3j_2j_1}=\int\limits_{-1}^{1}P_{j_4}(u)
\int\limits_{-1}^{u}P_{j_3}(z)
\int\limits_{-1}^{z}P_{j_2}(y)
\int\limits_{-1}^{y}
P_{j_1}(x)dx dy dz du,
$$

\vspace{6mm}
\noindent
where $P_j(x)$ is the Legendre polynomial.

Consider the stochastic integral $I_4[B(Z)]_{T,t}^M.$
Let conditions 1 and 2 of Theorem 4 be fulfilled.
Let $B'(v)(B''(v)(B(v),B(v)))$
be a 4-linear Hilbert--Schmidt operator 
mapping from
$U_0\times U_0\times U_0\times U_0$ to $H$
for all $v\in H$.

We have (see (\ref{xx605}))

\vspace{2mm}
$$
I_4[B(Z)]_{T,t}^M=
\sum_{r_1,r_2,r_3,r_4\in J_M}
B'(Z)\left(B''(Z)
\left(B(Z)e_{r_1}, B(Z)e_{r_2}\right)e_{r_3}\right)e_{r_4}
\sqrt{\lambda_{r_1}\lambda_{r_2}\lambda_{r_3}\lambda_{r_4}}\times
$$

\vspace{1mm}
\begin{equation}
\label{jjjk5}
\times\int\limits_t^T
\int\limits_t^s \left(
\int\limits_t^{\tau} d{\bf w}_{u}^{(r_1)}
\int\limits_t^{\tau} d{\bf w}_{u}^{(r_2)}\right)
d{\bf w}_{\tau}^{(r_3)}
d{\bf w}_s^{(r_4)}\ \ \ \hbox{w.\ p.\ 1.}
\end{equation}

\vspace{5mm}

From (\ref{jjj2}) and (\ref{jjjk5}) we obtain

\vspace{5mm}

$$
I_4[B(Z)]_{T,t}^M=
\sum_{r_1,r_2,r_3,r_4\in J_M}
B'(Z)\left(B''(Z)
\left(B(Z)e_{r_1}, B(Z)e_{r_2}\right)e_{r_3}\right)e_{r_4}
\sqrt{\lambda_{r_1}\lambda_{r_2}\lambda_{r_3}\lambda_{r_4}}\times
$$

\begin{equation}
\label{jjjk6}
\times 
\Biggl(I_{(1111)T,t}^{(r_1r_2r_3r_4)}+
I_{(1111)T,t}^{(r_2r_1r_3r_4)}-
{\bf 1}_{\{r_1=r_2\}}J_{(10)T,t}^{(r_3r_4)}\Biggr)
\ \ \ \hbox{w.\ p.\ 1,}
\end{equation}

\vspace{5mm}
\noindent
where 

\vspace{-2mm}
\begin{equation}
\label{ros2}
J_{(10)T,t}^{(r_3r_4)}=
\int\limits_t^T
\int\limits_t^s (t-\tau)d{\bf w}_{\tau}^{(r_3)}   
d{\bf w}_s^{(r_4)}.
\end{equation}

\vspace{4mm}

Denote by
$I_4[B(Z)]_{T,t}^{M,q}$ the approximation of the iterated stochastic integral
(\ref{jjjk6}), 
which has the following form

\vspace{3mm}

$$
I_4[B(Z)]_{T,t}^{M,q}=
\sum_{r_1,r_2,r_3,r_4\in J_M}
B'(Z)\left(B''(Z)
\left(B(Z)e_{r_1}, B(Z)e_{r_2}\right)e_{r_3}\right)e_{r_4}
\sqrt{\lambda_{r_1}\lambda_{r_2}\lambda_{r_3}\lambda_{r_4}}\times
$$

\vspace{1mm}
\begin{equation}
\label{jjjk7}
\times 
\Biggl(I_{(1111)T,t}^{(r_1r_2r_3r_4)q}+
I_{(1111)T,t}^{(r_2r_1r_3r_4)q}
-{\bf 1}_{\{r_1=r_2\}}J_{(10)T,t}^{(r_3r_4)q}\Biggr)
\ \ \ \hbox{w.\ p.\ 1,}
\end{equation}

\vspace{7mm}
\noindent
where 
the approximations $I_{(1111)T,t}^{(r_1r_2r_3r_4)q},$
$J_{(10)T,t}^{(r_1r_2)q}$
are based on Theorems 1, 2 and Legendre polynomials.

The approximation $J_{(10)T,t}^{(r_1 r_2)q}$ of the stochastic integral
$J_{(10)T,t}^{(r_1 r_2)}$ $(r_1,r_2=1,\ldots,M),$
which is based on Theorems 1, 2 and 
Legendre polynomials has the following form (see 
\cite{37} (formula (6.92), p.~A.544) or \cite{34} (formula (5.8), p.~A.249))

\vspace{2mm}

$$
J_{(10)T,t}^{(r_1 r_2)q}=
-\frac{T-t}{2}I_{(11)T,t}^{(r_1 r_2)q}
-\frac{(T-t)^2}{4}\Biggl(
\frac{1}{\sqrt{3}}\zeta_0^{(r_2)}\zeta_1^{(r_1)}+\Biggr.
$$

\vspace{2mm}
\begin{equation}
\label{4006x}
+\Biggl.\sum_{i=0}^{q}\Biggl(
\frac{(i+1)\zeta_{i+2}^{(r_2)}\zeta_{i}^{(r_1)}
-(i+2)\zeta_{i}^{(r_2)}\zeta_{i+2}^{(r_1)}}
{\sqrt{(2i+1)(2i+5)}(2i+3)}+
\frac{\zeta_i^{(r_1)}\zeta_{i}^{(r_2)}}{(2i-1)(2i+3)}\Biggr)\Biggr),
\end{equation}

\vspace{8mm}
\noindent
where the approximation $I_{(11)T,t}^{(r_1 r_2)q}$
is defined by (\ref{vini}).

Moreover,

\vspace{1mm}
\begin{equation}
\label{roza1}
{\sf M}\left\{\left(J_{(10)T,t}^{(r_1 r_2)}-
J_{(10)T,t}^{(r_1 r_2)q}\right)^2\right\}=E_q\ \ (r_1\ne r_2),
\end{equation}

\vspace{6mm}
\noindent
where 
$E_q$ is the right-hand side of (\ref{987}) 
(see 
\cite{37} (formula (6.106), p.~A.551) or
\cite{34} (formula (5.19), p.~A.252--A.253)).

From (\ref{jjjk6}), (\ref{jjjk7})
it follows that

\vspace{3mm}

$$
I_4[B(Z)]_{T,t}^{M}-
I_4[B(Z)]_{T,t}^{M,q}
=
$$

\vspace{3mm}
$$
=\sum_{r_1,r_2,r_3,r_4\in J_M}
B'(Z)\left(B''(Z)
\left(B(Z)e_{r_1}, B(Z)e_{r_2}\right)e_{r_3}\right)e_{r_4}
\sqrt{\lambda_{r_1}\lambda_{r_2}\lambda_{r_3}\lambda_{r_4}}\times
$$

\vspace{1mm}
$$
\times 
\Biggl(\left(I_{(1111)T,t}^{(r_1r_2r_3r_4)}
-I_{(1111)T,t}^{(r_1r_2r_3r_4)q}\right)
+
\left(I_{(1111)T,t}^{(r_2r_1r_3r_4)}
-I_{(1111)T,t}^{(r_2r_1r_3r_4)q}\right)-\Biggr.
$$

\vspace{1mm}
$$
\Biggl.-
{\bf 1}_{\{r_1=r_2\}}\left(J_{(10)T,t}^{(r_3r_4)}
-J_{(10)T,t}^{(r_3r_4)q}\right)\Biggr)
\ \ \ \hbox{w.\ p.\ 1.}
$$

\vspace{7mm}

Repeating with an insignificant 
modification the proof 
of Theorem 4 for the cases $k=2, 4$, we obtain

\vspace{2mm}
$$
\hspace{-70mm}{\sf M}\Biggl\{\biggl\Vert
I_4[B(Z)]_{T,t}^{M}-I_4[B(Z)]_{T,t}^{M,q}\biggr\Vert_H^2\Biggr\}\le
$$

\vspace{1mm}
$$
~~~~~~~~~~~~~~~~~~~~~~~~~~~~~\le
C
\left({\rm tr}\ Q\right)^4
\Biggl(2^2(4!)^2
\Biggl(\frac{(T-t)^{4}}{24}-\sum_{j_1,j_2,j_3,j_4=0}^{q}
C_{j_4j_3j_2j_1}^2\Biggr) +(2!)^2E_q\Biggr),
$$

\vspace{7mm}
\noindent
where 
$E_q$ is the  right-hand side of (\ref{987}), and
$C_{j_4j_3j_2j_1}$ is defined by (\ref{tq1}).

Consider the stochastic integral $I_5[B(Z)]_{T,t}^M.$
Let conditions 1 and 2 of Theorem 4 be fulfilled.
Let $B''(v)(B(v),B'(v)(B(v)))$
be a 4-linear Hilbert--Schmidt operator 
mapping from
$U_0\times U_0\times U_0\times U_0$ to $H$
for all $v\in H$.

We have (see (\ref{xx605}))

\vspace{2mm}
$$
I_5[B(Z)]_{T,t}^M=
\sum_{r_1,r_2,r_3,r_4\in J_M}
B''(Z)(B(Z)e_{r_3},
B'(Z)(B(Z)e_{r_2})e_{r_1})e_{r_4}
\sqrt{\lambda_{r_1}\lambda_{r_2}\lambda_{r_3}\lambda_{r_4}}\times
$$

\vspace{1mm}
\begin{equation}
\label{jjj500}
\times\int\limits_t^T\left(
\int\limits_t^s d{\bf w}_{\tau}^{(r_3)}
\int\limits_t^{s} 
\int\limits_t^{\tau} d{\bf w}_{u}^{(r_2)}
d{\bf w}_{\tau}^{(r_1)}
\right)
d{\bf w}_s^{(r_4)}\ \ \ \hbox{w.\ p.\ 1.}
\end{equation}

\vspace{5mm}

Using the theorem on the integration order replacement in iterated 
It\^{o} stochastic integrals (see 
\cite{37} (p.~A.150, p.~A.163), \cite{37a}-\cite{37aaxx}, 
\cite{integ}) or the It\^{o} formula, we obtain

\vspace{3mm}
$$
\int\limits_t^T\left(
\int\limits_t^s d{\bf w}_{\tau}^{(r_3)}
\int\limits_t^{s} 
\int\limits_t^{\tau} d{\bf w}_{u}^{(r_2)}
d{\bf w}_{\tau}^{(r_1)}
\right)
d{\bf w}_s^{(r_4)}=
$$

\vspace{2mm}
$$
=
I_{(1111)T,t}^{(r_2r_1r_3r_4)}+
I_{(1111)T,t}^{(r_2r_3r_1r_4)}+
I_{(1111)T,t}^{(r_3r_2r_1r_4)}+
$$

\begin{equation}
\label{jjj501}
+
{\bf 1}_{\{r_1=r_3\}}\left(J_{(10)T,t}^{(r_2r_4)}-
J_{(01)T,t}^{(r_2r_4)}\right)-
{\bf 1}_{\{r_2=r_3\}}J_{(10)T,t}^{(r_1r_4)}\ \ \ \hbox{w.\ p.\ 1,}
\end{equation}

\vspace{7mm}
\noindent
where 
we use the notations from Sect.~4,
and $J_{(01)T,t}^{(r_1r_2)}$, $J_{(10)T,t}^{(r_1r_2)}$
are defined by (\ref{ros1}), (\ref{ros2}).

After substituting (\ref{jjj501}) into (\ref{jjj500}), we obtain

\vspace{5mm}
$$
I_5[B(Z)]_{T,t}^M=
\sum_{r_1,r_2,r_3,r_4\in J_M}
B''(Z)(B(Z)e_{r_3},
B'(Z)(B(Z)e_{r_2})e_{r_1})e_{r_4}
\sqrt{\lambda_{r_1}\lambda_{r_2}\lambda_{r_3}\lambda_{r_4}}\times
$$

\vspace{2mm}
$$
\times\Biggl(I_{(1111)T,t}^{(r_2r_1r_3r_4)}+
I_{(1111)T,t}^{(r_2r_3r_1r_4)}+
I_{(1111)T,t}^{(r_3r_2r_1r_4)}+\Biggr.
$$

\begin{equation}
\label{ros3}
\Biggl.+
{\bf 1}_{\{r_1=r_3\}}\left(J_{(10)T,t}^{(r_2r_4)}-
J_{(01)T,t}^{(r_2r_4)}\right)-
{\bf 1}_{\{r_2=r_3\}}J_{(10)T,t}^{(r_1r_4)}\Biggr)\ \ \ \hbox{w.\ p.\ 1.}
\end{equation}

\vspace{5mm}

Denote by
$I_5[B(Z)]_{T,t}^{M,q}$ the approximation of the iterated stochastic integral
(\ref{ros3}), 
which has the following form

\vspace{5mm}
$$
I_5[B(Z)]_{T,t}^{M,q}=
\sum_{r_1,r_2,r_3,r_4\in J_M}
B''(Z)(B(Z)e_{r_3},
B'(Z)(B(Z)e_{r_2})e_{r_1})e_{r_4}
\sqrt{\lambda_{r_1}\lambda_{r_2}\lambda_{r_3}\lambda_{r_4}}\times
$$

\vspace{2mm}
$$
\times\Biggl(I_{(1111)T,t}^{(r_2r_1r_3r_4)q}+
I_{(1111)T,t}^{(r_2r_3r_1r_4)q}+
I_{(1111)T,t}^{(r_3r_2r_1r_4)q}+\Biggr.
$$

\begin{equation}
\label{ros4}
\Biggl.+
{\bf 1}_{\{r_1=r_3\}}\left(J_{(10)T,t}^{(r_2r_4)q}-
J_{(01)T,t}^{(r_2r_4)q}\right)-
{\bf 1}_{\{r_2=r_3\}}J_{(10)T,t}^{(r_1r_4)q}\Biggr)\ \ \ \hbox{w.\ p.\ 1,}
\end{equation}

\vspace{7mm}
\noindent
where 
the approximations $I_{(1111)T,t}^{(r_1r_2r_3r_4)q},$
$J_{(01)T,t}^{(r_1r_2)q},$ and $J_{(10)T,t}^{(r_1r_2)q}$
are based on Theorems 1, 2 and Legendre polynomials.

From (\ref{ros3}), (\ref{ros4}) it follows that

\vspace{3mm}
$$
I_5[B(Z)]_{T,t}^{M}-I_5[B(Z)]_{T,t}^{M,q}=
\sum_{r_1,r_2,r_3,r_4\in J_M}
B''(Z)(B(Z)e_{r_3},
B'(Z)(B(Z)e_{r_2})e_{r_1})e_{r_4}
\sqrt{\lambda_{r_1}\lambda_{r_2}\lambda_{r_3}\lambda_{r_4}}\times
$$

\vspace{1mm}
$$
\times\Biggl(\left(I_{(1111)T,t}^{(r_2r_1r_3r_4)}-
I_{(1111)T,t}^{(r_2r_1r_3r_4)q}\right)+
\left(I_{(1111)T,t}^{(r_2r_3r_1r_4)}-I_{(1111)T,t}^{(r_2r_3r_1r_4)q}\right)+
\left(I_{(1111)T,t}^{(r_3r_2r_1r_4)}-I_{(1111)T,t}^{(r_3r_2r_1r_4)q}\right)
+\Biggr.
$$

\vspace{1mm}
$$
\Biggl.+
{\bf 1}_{\{r_1=r_3\}}\left(\left(J_{(10)T,t}^{(r_2r_4)}-
J_{(10)T,t}^{(r_2r_4)q}\right)-
\left(J_{(01)T,t}^{(r_2r_4)}-J_{(01)T,t}^{(r_2r_4)q}\right)\right)-
{\bf 1}_{\{r_2=r_3\}}
\left(J_{(10)T,t}^{(r_1r_4)}-J_{(10)T,t}^{(r_1r_4)q}\right)
\Biggr)\ \ \ \hbox{w.\ p.\ 1.}
$$

\vspace{7mm}

Repeating with an insignificant 
modification the proof 
of Theorem 4 for the cases $k=2, 4$ and taking into account
(\ref{roza1}),
we obtain

\vspace{2mm}
$$
\hspace{-70mm}{\sf M}\Biggl\{\biggl\Vert
I_5[B(Z)]_{T,t}^{M}-I_5[B(Z)]_{T,t}^{M,q}\biggr\Vert_H^2\Biggr\}\le
$$

\vspace{1mm}
$$
~~~~~~~~~~~~~~~~~~~~~~~~~~~~~\le
C
\left({\rm tr}\ Q\right)^4
\Biggl(3^2(4!)^2
\Biggl(\frac{(T-t)^{4}}{24}-\sum_{j_1,j_2,j_3,j_4=0}^{q}
C_{j_4j_3j_2j_1}^2\Biggr) +3^2(2!)^2E_q\Biggr),
$$

\vspace{7mm}
\noindent
where 
$E_q$ is the right-hand side of (\ref{987}), and
$C_{j_4j_3j_2j_1}$ is defined by (\ref{tq1}).

Consider the stochastic integral $I_6[B(Z),F(Z)]_{T,t}^M.$
Let conditions 1 and 2 of Theorem 4 be fulfilled.
We have (see (\ref{xx605}))

\vspace{1mm}
$$
I_6[B(Z),F(Z)]_{T,t}^M=
\sum_{r_1,r_2\in J_M}
F'(Z)(B'(Z)(B(Z)e_{r_1})e_{r_2})
\sqrt{\lambda_{r_1}\lambda_{r_2}}\times
$$

\begin{equation}
\label{jjjk51}
\times\int\limits_t^T
\int\limits_t^s 
\int\limits_{t}^{\tau} d{\bf w}_{u}^{(r_1)}
d{\bf w}_{\tau}^{(r_2)}
ds\ \ \ \hbox{w.\ p.\ 1.}
\end{equation}

\vspace{5mm}

Using the theorem on the integration order replacement in iterated 
It\^{o} stochastic integrals (see 
\cite{37} (p.~A.150, p.~A.163), \cite{37a}-\cite{37aaxx}, 
\cite{integ}) or the It\^{o} formula, we obtain

\vspace{2mm}
\begin{equation}
\label{ff11}
\int\limits_t^T
\int\limits_t^s 
\int\limits_{t}^{\tau} d{\bf w}_{u}^{(r_1)}
d{\bf w}_{\tau}^{(r_2)}
ds=(T-t)I_{(11)T,t}^{(r_1r_2)}+J_{(01)T,t}^{(r_1r_2)}\ \ \ \hbox{w.\ p.\ 1.}
\end{equation}

\vspace{5mm}

After substituting (\ref{ff11}) into (\ref{jjjk51})
we have

\vspace{2mm}
$$
I_6[B(Z),F(Z)]_{T,t}^M=
\sum_{r_1,r_2\in J_M}
F'(Z)(B'(Z)(B(Z)e_{r_1})e_{r_2})
\sqrt{\lambda_{r_1}\lambda_{r_2}}\times
$$

\vspace{2mm}
\begin{equation}
\label{jjjk101}
\times\left(
(T-t)I_{(11)T,t}^{(r_1r_2)}+J_{(01)T,t}^{(r_1r_2)}
\right)\ \ \ \hbox{w.\ p.\ 1.}
\end{equation}

\vspace{7mm}

Denote by
$I_6[B(Z),F(Z)]_{T,t}^{M,q}$ the approximation of the
iterated stochastic integral
(\ref{jjjk101}), 
which has the following form

\vspace{3mm}
$$
I_6[B(Z),F(Z)]_{T,t}^{M,q}=                     
\sum_{r_1,r_2\in J_M}
F'(Z)(B'(Z)(B(Z)e_{r_1})e_{r_2})
\sqrt{\lambda_{r_1}\lambda_{r_2}}\times
$$

\vspace{2mm}
\begin{equation}
\label{jjjk102}
\times\left(
(T-t)I_{(11)T,t}^{(r_1r_2)q}+J_{(01)T,t}^{(r_1r_2)q}
\right),
\end{equation}

\vspace{7mm}
\noindent
where 
the approximations $J_{(01)T,t}^{(r_1r_2)q}$, $I_{(11)T,t}^{(r_1r_2)q}$
are defined by (\ref{vini0}), (\ref{vini}).

From (\ref{jjjk101}), (\ref{jjjk102}) it follows that

\vspace{3mm}
$$
I_6[B(Z),F(Z)]_{T,t}^M-I_6[B(Z),F(Z)]_{T,t}^{M,q}=                     
\sum_{r_1,r_2\in J_M}
F'(Z)(B'(Z)(B(Z)e_{r_1})e_{r_2})
\sqrt{\lambda_{r_1}\lambda_{r_2}}\times
$$

\vspace{2mm}
$$
\times\Biggl(
(T-t)\left(I_{(11)T,t}^{(r_1r_2)}-I_{(11)T,t}^{(r_1r_2)q}\right)
+\left(J_{(01)T,t}^{(r_1r_2)}-J_{(01)T,t}^{(r_1r_2)q}\right)
\Biggr)\ \ \ \hbox{w.\ p.\ 1.}
$$

\vspace{7mm}

Repeating with an insignificant 
modification the proof 
of Theorem 4 for the case $k=2$, we obtain

\vspace{2mm}
$$
\hspace{-40mm}
{\sf M}\left\{\biggl\Vert
I_6[B(Z),F(Z)]_{T,t}^M-I_6[B(Z),F(Z)]_{T,t}^{M,q}\biggr\Vert_H^2\right\}\le
$$

\vspace{2mm}
$$
~~~~~~~~~~~~~~~~~~~~~~~~~~~~~~~~~~~~~~~~~~~~~~~~~~~~~~~\le
2C(2!)^2\left({\rm tr}\ Q\right)^2
\biggl((T-t)^2 G_q+E_q\biggr),
$$

\vspace{7mm}
\noindent
where $G_q$ and $E_q$ are the right-hand sides of (\ref{909}) and
(\ref{987}) correspondingly.

Consider the stochastic integral $I_7[B(Z),F(Z)]_{T,t}^M.$
Let conditions 1 and 2 of Theorem 4 be fulfilled.
Then we have (see (\ref{xx605})) w.\ p.\ 1

\vspace{2mm}
$$
I_7[B(Z),F(Z)]_{T,t}^M=
\sum_{r_1,r_2\in J_M}F''(Z)\left(B(Z)e_{r_1}, B(Z)e_{r_2}\right)
\sqrt{\lambda_{r_1}\lambda_{r_2}}\times
$$

\begin{equation}
\label{da11}
\times \int\limits_t^T
\left(\int\limits_t^s d{\bf w}_{\tau}^{(r_1)}   
\int\limits_t^s d{\bf w}_{\tau}^{(r_2)}\right)
ds.
\end{equation}

\vspace{4mm}

Using the It\^{o} formula, we obtain

\begin{equation}
\label{da12w}
\int\limits_t^s d{\bf w}_{\tau}^{(r_1)}   
\int\limits_t^s d{\bf w}_{\tau}^{(r_2)}
=
I_{(11)s,t}^{(r_1 r_2)}+
I_{(11)s,t}^{(r_2 r_1)}+
{\bf 1}_{\{r_1=r_2\}}(s-t)\ \ \ \hbox{w.\ p.\ 1,}
\end{equation}

\vspace{3mm}
\noindent
where we use the notations from Sect.~4.

From (\ref{da12w}) and (\ref{ff11}) we have

\vspace{4mm}
$$
\int\limits_t^T
\left(\int\limits_t^s d{\bf w}_{\tau}^{(r_1)}   
\int\limits_t^s d{\bf w}_{\tau}^{(r_2)}\right)ds=
$$

\vspace{1mm}
$$
=
\int\limits_t^T
I_{(11)s,t}^{(r_1 r_2)}ds+
\int\limits_t^T
I_{(11)s,t}^{(r_2 r_1)}ds+{\bf 1}_{\{r_1=r_2\}}\frac{(T-t)^2}{2}=
$$

\vspace{2mm}
$$
=(T-t)\left(I_{(11)T,t}^{(r_1 r_2)}+
I_{(11)T,t}^{(r_2 r_1)}\right)+
J_{(01)T,t}^{(r_1 r_2)}+J_{(01)T,t}^{(r_2 r_1)}+
{\bf 1}_{\{r_1=r_2\}}\frac{(T-t)^2}{2}=
$$

\vspace{2mm}
$$
=
(T-t)\left(I_{(1)T,t}^{(r_1)}I_{(1)T,t}^{(r_2)}-
{\bf 1}_{\{r_1=r_2\}}(T-t)\right)+
J_{(01)T,t}^{(r_1 r_2)}+J_{(01)T,t}^{(r_2 r_1)}+
{\bf 1}_{\{r_1=r_2\}}\frac{(T-t)^2}{2}=
$$

\vspace{2mm}
$$
=
(T-t)I_{(1)T,t}^{(r_1)}I_{(1)T,t}^{(r_2)}+
J_{(01)T,t}^{(r_1 r_2)}+J_{(01)T,t}^{(r_2 r_1)}-
$$

\vspace{2mm}
\begin{equation}
\label{da13}
-
{\bf 1}_{\{r_1=r_2\}}\frac{(T-t)^2}{2}\ \ \ \hbox{w.\ p.\ 1.}
\end{equation}

\vspace{9mm}

After substituting (\ref{da13}) into (\ref{da11}) we obtain

\vspace{5mm}
$$
I_7[B(Z),F(Z)]_{T,t}^M=
\sum_{r_1,r_2\in J_M}F''(Z)\left(B(Z)e_{r_1}, B(Z)e_{r_2}\right)
\sqrt{\lambda_{r_1}\lambda_{r_2}}\times
$$

\vspace{1mm}
\begin{equation}
\label{da14}
\times
\left(
(T-t)I_{(1)T,t}^{(r_1)}I_{(1)T,t}^{(r_2)}+
J_{(01)T,t}^{(r_1 r_2)}+J_{(01)T,t}^{(r_2 r_1)}-
{\bf 1}_{\{r_1=r_2\}}\frac{(T-t)^2}{2}\right)\ \ \ \hbox{w.\ p.\ 1.}
\end{equation}

\vspace{8mm}

Denote by
$I_7[B(Z),F(Z)]_{T,t}^{M,q}$ the approximation of 
the
iterated stochastic integral
(\ref{da14}),
which has the following form

\vspace{4mm}
$$
I_7[B(Z),F(Z)]_{T,t}^{M,q}=                     
\sum_{r_1,r_2\in J_M}F''(Z)\left(B(Z)e_{r_1}, B(Z)e_{r_2}\right)
\sqrt{\lambda_{r_1}\lambda_{r_2}}\times
$$

\vspace{2mm}
\begin{equation}
\label{da20}
\times
\left(
(T-t)I_{(1)T,t}^{(r_1)}I_{(1)T,t}^{(r_2)}+
J_{(01)T,t}^{(r_1 r_2)q}+J_{(01)T,t}^{(r_2 r_1)q}-
{\bf 1}_{\{r_1=r_2\}}\frac{(T-t)^2}{2}\right),
\end{equation}

\vspace{8mm}
\noindent
where 
the approximation $J_{(01)T,t}^{(r_1r_2)q}$
is defined by (\ref{vini0}).

From (\ref{da14}), (\ref{da20}) it follows that

\vspace{4mm}
$$
I_7[B(Z),F(Z)]_{T,t}^M-I_7[B(Z),F(Z)]_{T,t}^{M,q}=                     
\sum_{r_1,r_2\in J_M}F''(Z)\left(B(Z)e_{r_1}, B(Z)e_{r_2}\right)
\sqrt{\lambda_{r_1}\lambda_{r_2}}\times
$$

\vspace{2mm}
$$
\times
\Biggl(
\left(J_{(01)T,t}^{(r_1 r_2)}-J_{(01)T,t}^{(r_1 r_2)q}\right)
+\left(J_{(01)T,t}^{(r_2 r_1)}-J_{(01)T,t}^{(r_2 r_1)q}\right)
\Biggr)\ \ \ \hbox{w.\ p.\ 1.}
$$

\vspace{8mm}

Repeating with an insignificant 
modification the proof 
of Theorem 4 for the case $k=2$, we obtain

\vspace{3mm}
$$
\hspace{-30mm}
{\sf M}\left\{\biggl\Vert
I_7[B(Z),F(Z)]_{T,t}^M-I_7[B(Z),F(Z)]_{T,t}^{M,q}\biggr\Vert_H^2\right\}\le
$$

\vspace{2mm}
$$
~~~~~~~~~~~~~~~~~~~~~~~~~~~~~~~~~~~~~~~~~~~~~~~~~~~~\le
2^2 C (2!)^2\left({\rm tr}\ Q\right)^2
E_q,
$$

\vspace{7mm}
\noindent
where $E_q$ is the right-hand side of (\ref{987}).

Consider the stochastic integral $I_8[B(Z),F(Z)]_{T,t}^M.$
Let conditions 1 and 2 of Theorem 4 be fulfilled.
Then we have (see (\ref{xx605})) w.\ p.\ 1

\vspace{3mm}
\begin{equation}
\label{da110}
I_8[B(Z),F(Z)]_{T,t}^M=-
\sum_{r_1,r_2\in J_M}B''(Z)\left(F(Z), B(Z)e_{r_1}\right)e_{r_2}
\sqrt{\lambda_{r_1}\lambda_{r_2}}J_{(01)T,t}^{(r_1r_2)}.
\end{equation}

\vspace{5mm}

Denote by
$I_8[B(Z),F(Z)]_{T,t}^{M,q}$ the approximation of the
iterated stochastic integral
(\ref{da110}),
which has the following form

\vspace{3mm}
\begin{equation}
\label{da200}
I_8[B(Z),F(Z)]_{T,t}^{M,q}=                     
-\sum_{r_1,r_2\in J_M}B''(Z)\left(F(Z), B(Z)e_{r_1}\right)e_{r_2}
\sqrt{\lambda_{r_1}\lambda_{r_2}}J_{(01)T,t}^{(r_1r_2)q},
\end{equation}

\vspace{6mm}
\noindent
where 
the approximation $J_{(01)T,t}^{(r_1r_2)q}$
is defined by (\ref{vini0}).

From (\ref{da110}), (\ref{da200}) it follows that

\vspace{4mm}
$$
I_8[B(Z),F(Z)]_{T,t}^M-I_8[B(Z),F(Z)]_{T,t}^{M,q}=
$$

\vspace{2mm}
$$
=-
\sum_{r_1,r_2\in J_M}B''(Z)\left(F(Z), B(Z)e_{r_1}\right)e_{r_2}
\sqrt{\lambda_{r_1}\lambda_{r_2}}\times
$$

\vspace{2mm}
$$
\times\left(J_{(01)T,t}^{(r_1r_2)}-
J_{(01)T,t}^{(r_1r_2)q}\right)\ \ \ \hbox{w.\ p.\ 1.}
$$

\vspace{8mm}

Repeating with an insignificant 
modification the proof 
of Theorem 4 for the case $k=2$, we obtain

\vspace{5mm}
$$
{\sf M}\left\{\biggl\Vert
I_8[B(Z),F(Z)]_{T,t}^M-I_8[B(Z),F(Z)]_{T,t}^{M,q}\biggr\Vert_H^2\right\}\le
C (2!)^2\left({\rm tr}\ Q\right)^2
E_q,
$$

\vspace{6mm}
\noindent
where $E_q$ is the right-hand side of (\ref{987}).

Using computational experiments
it was shown in \cite{1002}, \cite{1003} (also see \cite{37a}, Sect.~5.4) that
we can neglect the 
multiplier factor $k!$ in the estimate (\ref{star00011}).
As a result, the computational costs for the approximation 
of iterated It\^{o} stochastic integrals are significantly reduced.
For the same reason, we can replace the multiplier factor
$(k!)^2$ by $k!$ in the formula (\ref{zzz1}) in practical calculations.

\vspace{8mm}
\noindent
{\bf Acknowledgement}. I 
would like to thank Leonid Makarovsky for his help in 
translation this article into English and Konstantin Rybakov 
for useful discussion of some presented results.

\end{document}